\documentclass[a4paper,11pt]{article}
\usepackage[noblocks]{authblk}
\usepackage{todonotes}
\usepackage{float}
\usepackage{amsmath}
\usepackage{amssymb}
\usepackage{amsthm}
\usepackage{pifont}
\usepackage{graphicx}
\usepackage{xcolor}
\usepackage{fullpage}
\usepackage[percent]{overpic}

\newcommand{\nodes}{{\mathcal A}}

\graphicspath{{./Figures/}}

\usepackage{xcolor}
\definecolor{newcolor}{rgb}{.8,.349,.1}
\usepackage{mathrsfs} 

\usepackage{amsmath}
\usepackage[noblocks]{authblk}
\usepackage{amsthm}
\usepackage{algpseudocode}
\usepackage{algorithm}
\usepackage{amsfonts}
\usepackage{graphicx}
\usepackage{geometry}
\usepackage{bm}
\usepackage{todonotes}
\usepackage{mathtools}
\usepackage{comment}
\usepackage{epstopdf}
\usepackage{float}
\usepackage{placeins}
\newcommand{\vx}{\vec{x}}
\newcommand{\tpar}{\widetilde{\partial}}
\usepackage{xcolor}
\usepackage{amsmath}
\usepackage{amssymb}
\usepackage{amsthm}
\usepackage{amsfonts}
\usepackage{graphicx}
\usepackage{geometry}
\usepackage{bm}
\usepackage{todonotes}

\usepackage{mathtools}
\newtheorem{pro}{Proposition}

\newtheorem{remark}{Remark}
\usepackage[percent]{overpic}
\usepackage[latin1]{inputenc}
\usepackage{epstopdf}
\usepackage{float}
\usepackage{placeins}
\graphicspath{ {./Figures/} }
\newcommand{\pad}[2]{\frac{\partial #1}{\partial #2}}

\newcommand{\td}[2] {\frac{ {\rm d} #1}{ {\rm d} #2}}
\newcommand{\BB}{\text{\Large$\mathbb{B}^\varepsilon$} }

\title{Asymptotic Preserving and Accurate scheme for Multiscale Poisson-Nernst-Planck (MPNP) system}

\author{Clarissa Astuto}
\author{Giovanni Russo}
\affil{Department of Mathematics and Computational Science, University of Catania, Italy}

\begin{document}
\maketitle

\begin{abstract}
In this paper, we propose and validate a two-species Multiscale model for a Poisson-Nernst-Planck (PNP) system, focusing on the correlated motion of positive and negative ions under the influence of a  trap. Specifically, we aim to model surface traps whose attraction range, of length $\delta$, is much smaller then the scale of the problem. The physical setup we refer to is an anchored gas drop (bubble) surrounded by a flow of charged surfactants {(composed by positive and negative ions) that diffuses in water. When the diffusing
surfactants reach the surface of the trap, the negative ions are adsorbed because of their hydrophobic tail that is attracted by the air bubble}. As in our previous works \cite{astuto2023multiscale,ASTUTO2023111880,astuto2025time,astuto2024high}, the effect of the attractive potential is replaced by a suitable boundary condition derived by mass conservation and asymptotic analysis. The novelty of this work is the extension of the model proposed in \cite{astuto2023multiscale}, now incorporating the influence of both carriers -- positive and negative ions -- simultaneously, which is often neglected in traditional approaches that treat ion species independently. {The two carriers interact through the Coulomb potential, that is computed by a Poisson equation. This leads to a multiscale model with two additional equations compared to the initial problem formulated in \cite{astuto2023multiscale}.}
In the second part of the paper, we address the treatment of the Coulomb interaction. When the Debye length $\lambda_D$ ({related} to a small parameter $\varepsilon$) is very small, one can adopt the so-called Quasi-Neutral limit, which significantly simplifies the system, reducing it to a {single} diffusion equation for {the sum of the two carriers} with effective diffusion coefficient \cite{jungel,CiCP-31-707}. {While this approach significantly simplifies the mathematical model, by reducing the system from three equations to a single one and eliminating the stiffness induced by a vanishingly small Debye length, it fails to capture the effects arising for non-negligible values of $\varepsilon$. In the regime where the Debye length is small but not negligible, capturing small deviations from the quasi-neutral limit may become computationally expensive when using standard methods in the literature.} 
One of the objectives of this work is to develop an \textit{Asymptotic Preserving} (AP) second order numerical scheme that works for all Debye lengths and becomes a consistent discretization of the Quasi Neutral limit as $\varepsilon \to 0$, with no stability restriction on the time step.
Furthermore, the numerical scheme we propose is also \textit{Asymptotic Accurate} (AA), which means that it preserves second order accuracy in the Quasi-Neutral limit.
\end{abstract}

\section{Introduction}
In this work, we are interested in modeling the chemical trapping of heterogeneous substances, such us surfactants, since the existing methods often rely on physical techniques that monitor only one component at a time \cite{osti_6724746,fernandez2016existence}, limiting their scope. We propose a model to simultaneously measure the interfacial concentrations and distributions of different ions, such as positive and negative ions. 

Aqueous surfactants are important in diverse applications, including biological and biochemical processes. They influence foam properties \cite{brown1990foam}, wettability, coating flows \cite{valentini1991role}, and {are widely used as spray to increase efficacy of foliar-applied agrochemicals, enhancing pesticide penetration into foliage of a wide range of plant species \cite{knoche1991performance}.} Surfactants are also applied to pulmonary mechanics, specifically in the context of the human lungs and alveoli \cite{notter1975pulmonary}. The surfactant molecules present in the alveolar lining layer play a crucial role in controlling the surface tension at the liquid-air interface within the alveoli. The surfactants help stabilize the lungs during respiration. In \cite{tomlinson2023unsteady,tomlinson2023laminar}, the authors investigate the effect of soluble surfactant fluctuations on drag reduction in superhydrophobic channels, focusing on how surfactant transport and adsorption modify flow dynamics; using asymptotic modeling, they derive coupled equations for surfactant and fluid behaviors under laminar, pressure-driven flow. The study predicts that nonlinear wave phenomena, including shock formation in surfactant flux, can degrade slip and drag-reduction performance. Results highlight the interplay between surfactant concentration gradients, adsorption effects, and effective slip length, providing insights into optimizing superhydrophobic surfaces for fluid applications.

\begin{figure}[h]
	\centering
	\centering
\begin{overpic}[abs,width=0.6\textwidth,unit=1mm,scale=.25]{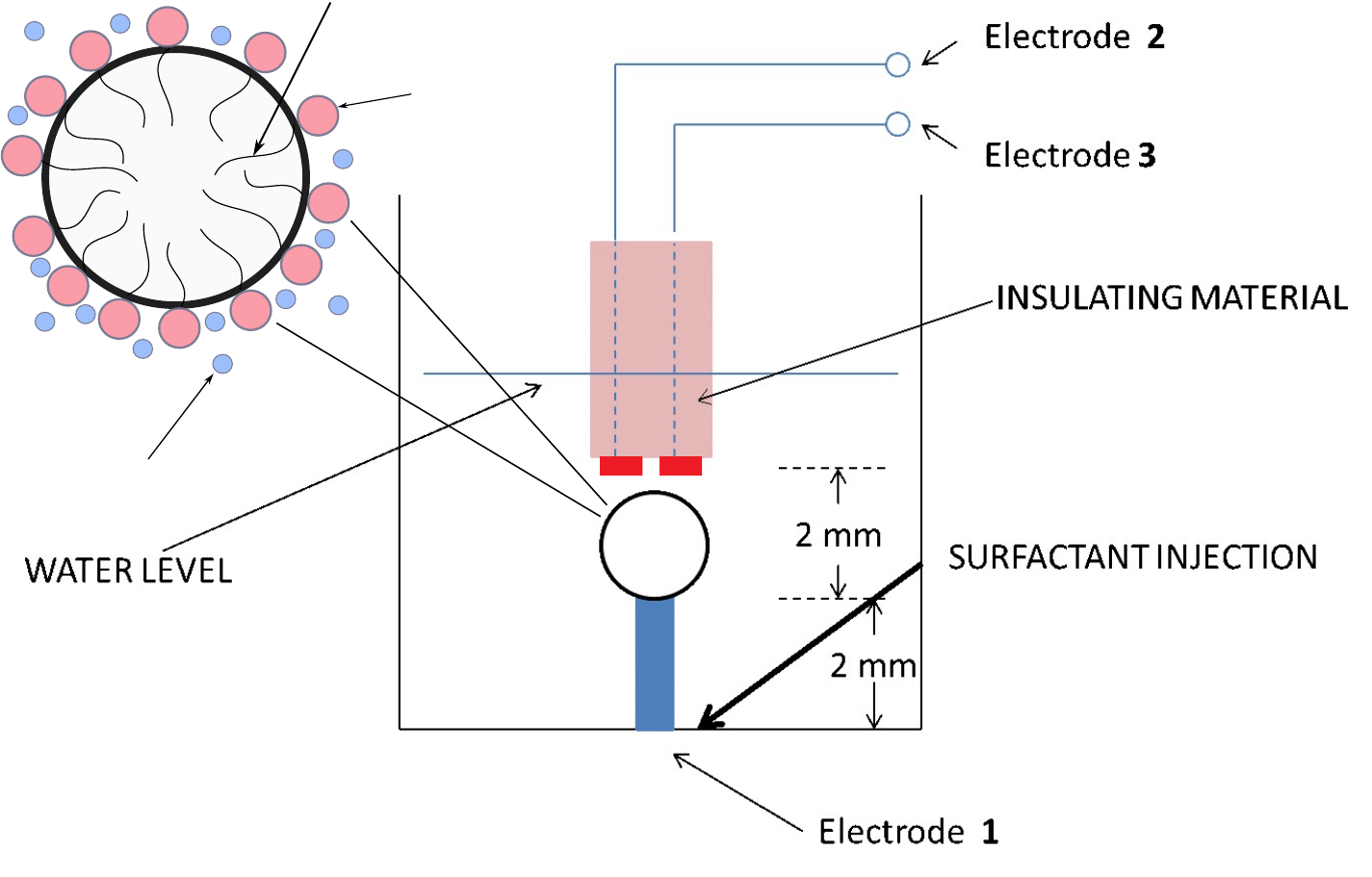}
\put(13,62.5){\footnotesize hydrophobic tail}
\put(22,56.5){\footnotesize hydrophilic}
\put(26.5,53){\footnotesize head}
\put(8,27.5){\footnotesize cation}
\end{overpic}
	\caption{\textit{Scheme of the experimental setup.  On the top left there is a zoom in of the anions and cations behavior at the air-surface of the bubble: the cations (blue) are composed by hydrophilic heads;  the anions (red) have hydrophobic tails inside the air bubble, and hydrophilic heads on the surface.}}
\label{fig_setup}
\end{figure}
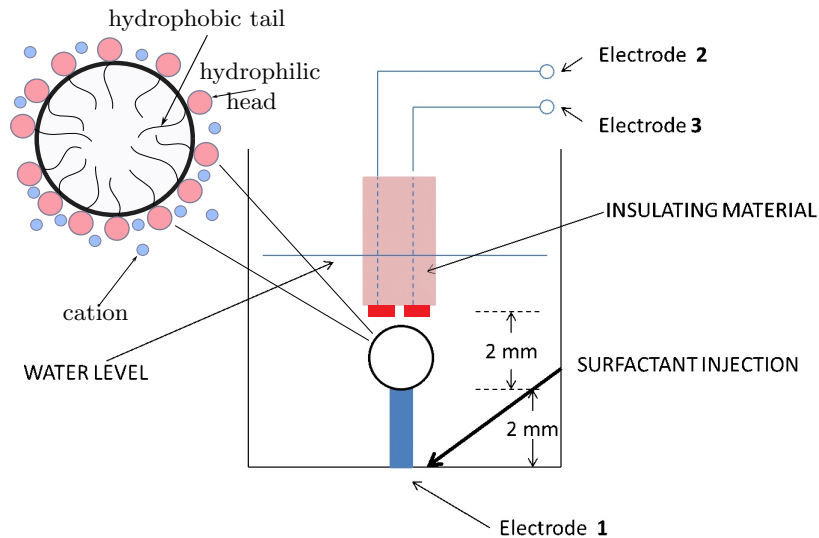

{In this paper, we consider a biomimetic experimental model system designed to mimic the capture rates (chemoreception) of a diffusing substance \cite{Raudino20168574}. In particular, we investigate the role of surface traps acting on species freely diffusing. When surfactants reach the surface of the bubble, they undergo reversible adsorption, allowing for an accurate measurement of their local concentration. More generally, this setting provides a prototypical example of a phenomenon that requires the coupling of multiple length scales, ranging from the molecular level up to the assembly scale that may even extend to the millimeter range \cite{van2007self}. The system considered in \cite{Raudino20168574} involves} an anchored gas droplet subjected to a diffusive flow of charged surfactants, where conductivity measurements are used to detect surfactant concentration beyond the oscillating bubble, see Fig.~\ref{fig_setup}. Pulsating bubble devices have been widely used to measure the dynamic surface tension because of their simplicity \cite{lu1991shape,piedfert2018numerical,neergaard1929new,notter1975pulmonary,fernandez2016existence,LIAO2006183}. 

In this work, we consider the PNP model for the diffusion of the two carriers \cite{EISENBERG2007,LU20112475,wang2017modeling}. {The drift term is given by the gradient of a potential which includes both the effect of the bubble and the Coulomb interaction between the carriers, described by $V_\pm$ and $\varphi$, respectively. The latter term is obtained from the solution of a self-consistent Poisson equation.}  We begin with a {$d$-}dimensional model, and consider two different regions of the domain, one of which represents the bubble and the other one the fluid. We continue with a simple one dimensional problem to deduce and validate new boundary conditions at the surface of the bubble (that corresponds to a point in one dimension), which mimics the attractive-repulsive interaction of the bubble on the surfactants. Following the strategy introduced in \cite{astuto2023multiscale} for one single carrier, { we extend the framework to a coupled setting for two carriers. In particular, the system is augmented by two additional equations: one governing the diffusion of positive ions, $c_+$, and one for the electrostatic potential $\varphi$, which describes the coupling between anions and cations. This extension leads to the derivation of three novel time-dependent boundary conditions at the bubble surface for the three variables: the two ionic concentrations and the electrostatic potential.} 

{After deducing the new Multiscale PNP model (MPNP), we develop an AP numerical scheme that has no restriction on the time step for all values of vanishingly small Debye lengths $\varepsilon\geq 0$, and that works also in the so called \textit{Quasi--Neutral Limit} (QNL) regime, $\varepsilon\to 0$, as seen, for instance, in \cite{jungel,CiCP-31-707}. Although in the present application, characterized by very low concentrations, the quasineutrality regime is frequently observed, there exist many biologically relevant situations in which electrostatic interactions play a fundamental role and cannot be neglected. In particular, processes such as protein misfolding and aggregation, which are associated with a wide range of pathological conditions including Alzheimer and Parkinson diseases and type II diabetes, involve interactions occurring at scales comparable to the Debye length \cite{gsponer2006theoretical,van2007self}. }

We introduce here the local {molar} concentrations of cations $c_+= c_+(\vec{x},t)$ and of anions $c_-= c_-(\vec{x},t)$. Their time evolution in a fluid is governed by the following conservation laws
\begin{equation}
	\displaystyle \frac{\partial c_\pm\left(\vec{x},t\right)}{\partial t} + \nabla\cdot J_\pm(\vec{x},t) = 0.
	\label{equation_flux}
\end{equation}
where $J_\pm$ are the fluxes associated to $c_\pm$, respectively, $\vec{x} \in {\Omega \subset} \mathbb R^d, \, d\geq 1$ and $t>0$. In simple diffusion, the fluxes consist of the gradient of the concentrations. However, in our case, the fluxes are augmented by the presence of additional potentials, which significantly influence the system behavior. These potentials stem from two distinct contributions: an external potential, describing the effect of the bubble on the ions, which has been extensively analyzed in previous works \cite{astuto2023multiscale,ASTUTO2023111880,astuto2025time,SIAM1_fernandez2016existence,SIAM3_morgan2015mathematical}, and the electrostatic interaction between the ions \cite{CiCP-31-707}. These combined effects introduce a more complex dynamic, altering the standard diffusion process by incorporating both external forces and ion-to-ion interactions. The dynamics have been extensively studied using various numerical approaches: an arbitrary Lagrangian-Eulerian finite-element method \cite{SIAM2_ganesan2012arbitrary}, a finite-volume method for bulk diffusion combined with a Voronoi decomposition for surface diffusion \cite{novak2007diffusion}, a CutFEM for coupled bulk-surface problems on time-dependent domains, where a level-set function describes the time evolution of the interface \cite{hansbo2016cut}, and a free-boundary formulation incorporating a mathematical model for the evolving interface \cite{SIAM3_morgan2015mathematical}.

Eqs.~\eqref{equation_flux} are coupled self-consistently to the Poisson equation for the electrostatic potential, $\varphi$, between ions, as follows 
\begin{equation}
	\displaystyle -\epsilon_0\epsilon_r{\Delta\varphi}=q\left(\tilde n_+-\tilde n_-\right)
	\label{equation_poisson_phi}
\end{equation}
where $\epsilon_0$ is the vacuum permittivity, $\epsilon_r$ is the relative permittivity,  $q $ is the (positive) electron charge and $\tilde n_\pm$ are the ion charge density which are proportional to the ion concentrations ${c_\pm}$ by the relation:
\begin{equation}
	{ \tilde n_\pm=\frac{{c_\pm} N_A\rho}{\tilde{m}} }
	\label{eq_n_pm}
\end{equation}
where $N_A$ is the Avogadro's number, {$\rho$ is the water mass density, $\tilde m = m_0 m_{H_2O}$  is the mass of one mole of water (expressed in Kg/mol) and $m_{H_2O} = 18$ is the molecular mass of water.}
Multiplying Eq.~\eqref{equation_poisson_phi} by $q/(\epsilon_0 \epsilon_r)$ and replacing $\tilde n_\pm$ with Eq.~\eqref{eq_n_pm}, we obtain
\begin{equation}
	{	\displaystyle -  {\Delta} (q\varphi) = \frac{q^2 N_A \rho}{\epsilon_0\epsilon_r \tilde{m}}\left(c_+-c_-\right). }
		\label{phiNond}
\end{equation}
Dividing by $k_B T$ (where $k_B$ is the Boltzmann's constant, $T$ is the absolute temperature, assumed to be constant), Eq.~\eqref{phiNond} becomes
\begin{equation}	
	{	 -  {\Delta} \Phi = {K}\left(c_+ - c_-\right) }
	\label{phiNond2}
\end{equation}
where {$\displaystyle \Phi = \frac{q\varphi}{k_B T}$ is dimensionless and  $\displaystyle {K}=\frac{q^2N_A \rho}{k_B T\epsilon_0\epsilon_r \tilde m}$. 
Eq.~\eqref{phiNond2} can be written as 
\begin{equation}
\label{eq:npm}
	-\varepsilon {\Delta}\Phi =c_+ - c_-,	
\end{equation} 
where $\varepsilon :=  \lambda_D^2 =K^{-1}$ has the dimension of the square of a length, and its value comes from the  parameters shown in Table~\ref{table_parameters}.} 

At the end, the system reads 
\begin{subequations}
\label{eq_full_model_adim}
\begin{align}
		\displaystyle \frac{\partial c_\pm}{\partial t} &= - \nabla \cdot J_\pm, \quad {\rm in }\, \Omega\\
		\displaystyle J_\pm &= - D_\pm \left(\nabla c_\pm + c_\pm \nabla (U_\pm \pm \Phi)\right), \quad {\rm in }\, \Omega\\ 
		-\varepsilon \Delta \Phi &= {c_+ - c_-}, \quad {\rm in }\, \Omega\\		
		\displaystyle J_\pm \cdot \widehat n &= 0, \quad {\rm on }\, \partial\Omega \\
        \nabla \Phi \cdot \widehat n &= 0, \quad {\rm on }\, \partial\Omega,
\end{align}
\end{subequations}
where $D_\pm$ are the diffusion coefficients, respectively, for $c_\pm$, and $U_\pm = V_\pm/k_BT$ the suitable potential functions that model the \textit{attractive-repulsive} behavior of the bubble with the ions. 
 Regarding the external potential for cations, $V_+$, the bubble behavior is always repulsive, acting like a wall along its surface. With the negative ions, the potential $V_-$ acts in a different way. In particular, when an anion is close to the bubble, this acts as a trap and the particle feels attraction towards the surface. On the contrary, when the anion is at a very short distance from the surface of the bubble, the potential is designed to repulse the particle in order to mimic impermeability; see {Fig.~\ref{fig:Lpotential} (a). In this paper, we investigate two dimensional cases. We start with $\Omega \subset \mathbb R^d,\, d =1$, where we set the center of the external potentials at $x = 0$. Furthermore, we set $d = 2, $ a more realistic domain where the bubble is firstly represented by a potential with a radial symmetry and then replaced by as a circular hole in the multiscale setup.}

In this paper, we first deduce new boundary conditions for the two-carrier MPNP model, and secondly, we design and study a two dimensional numerical scheme for the Poisson-Nernst-Planck system, applicable for any values of $\varepsilon$. These asymptotic limiting processes have been widely studied. To mention a few, in \cite{alves2024zero}, the authors examine a system of bipolar Euler-Poisson system, focusing on two asymptotic limiting processes. The first involves the limit of zero electron mass. In the second step, they explore the simultaneous application of both the zero-electron-mass limit and the Quasi-Neutral limit. In \cite{MarkowichQuasineutral}, the authors investigate the classical time-dependent drift-diffusion model for semiconductors, focusing on scenarios where the Debye length is small, treating it as a singular perturbation parameter.

The numerical scheme we employ converges for any value of $\varepsilon \geq 0$ and remains stable at the Quasi-Neutral limit as $\varepsilon\to 0$. For the time discretization, we employ a second-order Asymptotic Preserving {(AP)} numerical scheme. This approach utilizes an IMplicit-EXplicit (IMEX) strategy \cite{pareschi2000implicit,boscarino2016high,dimarco2013asymptotic,pareschi2000asymptotic,boscarino2024implicit}, where the stiff terms of the equations are treated implicitly to ensure stability, while the non-stiff components are handled explicitly to preserve computational efficiency. 
{The concept of an AP method has been introduced in \cite{jin1999efficient}. A problem $\mathcal{P}^\varepsilon$, that depends on a (small) parameter $\varepsilon$, is considered. As the parameter vanishes, the problem becomes the limit problem  $\mathcal{P}^0$. Examples from the literature include compressible Euler equations that relax to the incompressible Euler system as the Mach number vanishes \cite{degond2009asymptotic,filbet2010class}, or the Boltzmann equation of rarefied gas dynamics, which relaxes to the compressible Euler equations as the Knudsen number vanishes \cite{filbet2011asymptotic}. Then, a discretized version of the problem is considered, which is denoted by $\mathcal{P}^\varepsilon_h$, and $h$ is the discretization parameter. As $h\to 0$, the discrete problem converges to the continuous one. A numerical scheme $\mathcal P^\varepsilon_h$ for problem $\mathcal{P}^\varepsilon$ with discretization parameter $h$ is called AP if it becomes a consistent discretization of the limit problem $\mathcal{P}^0$ as $\varepsilon\to 0$.}

In space, we consider a ghost nodal finite element method, recently developed in \cite{astuto2024nodal} and further applied in \cite{astuto2024self,astuto2024comparison,dilip2025multigrid}. Since the numerical method does not require the use of a mesh fitted to the domain over which we are solving the MPNP system, it belongs to the realm of the ``unfitted'' finite element methods. Other examples of unfitted FEM are the so called CutFEM \cite{burman2012fictitious,burman2015cutfem,lehrenfeld2016high,lehrenfeld2018analysis,burman2022cutfem}, or the AgFEM \cite{BADIA202360}. There are other numerical schemes based on finite volume in space, such as \cite{bessemoulin2014study}, where the authors perform a numerical approximation of the classical time-dependent drift-diffusion system near quasi-neutrality, with a fully implicit time discretization combined with a finite volume method in space, approximating the convection-diffusion fluxes using Scharfetter-Gummel fluxes. 
 In \cite{belaouar2009asymptotically}, the authors develop a semi-Lagrangian scheme for the Vlasov-Poisson equation in the Quasi-Neutral regime. The key is a reformulation of the Poisson equation that allows for asymptotically stable simulations and the advantage is that this approach has no restriction on the time step as the Debye length and plasma period approach zero. \cite{brull2011asymptotic} is another paper based on AP semi-discretization in time for the simulation of a strongly magnetized plasma considered as a mixture of an ion fluid and an electron fluid, described by Euler equations. Regarding high order numerical schemes, in 
\cite{crouseilles2024high}, the authors develop IMEX finite volume methods for simulating plasmas in quasineutral regimes; to overcome stability challenges at small scales near the Debye length, the authors propose a class of penalized IMEX Runge-Kutta methods tailored for the Euler-Poisson system. 

{The paper is organized as follows: in Section~\ref{sec:MPNP}, we first deduce the multiscale model in one-dimension and then its extension in higher (two and three) dimensions. Before introducing the quasi-neutrality, we ensure the validity of the MPNP model, first describing the numerical scheme in one dimension, in Section~\ref{sec:discr1D}, and then showing its numerical results in Section~\ref{section_accuracy_test_1D}. In Section~\ref{sec:QNL}, we derive the quasi-neutral limit of the MPNP model. To this end, we first introduce two new variables (the sum of the concentrations  and their difference divided by $\varepsilon$), and then deduce the corresponding time-dependent boundary conditions. Finally, we perform the limit $\varepsilon \to 0$ on the resulting system. In Section~\ref{sec_2D}, we describe the discretization of the system in two dimensions. We first introduce its variational formulation, then we apply a ghost-FEM scheme for the space discretization, and finally present an Asymptotic-Preserving IMEX numerical scheme for the time discretization.  Section~\ref{section_numerics_2D} is devoted to the numerical results in two dimensions. Finally, in Section~\ref{sec:conclusions}, we draw some Conclusions and future work.}

\section{Multiscale Poisson-Nernst-Planck system (MPNP)}
\label{sec:MPNP}
\section{Multiscale Poisson-Nernst-Planck (MPNP) system }
\label{sec:MPNP}
In this section, we first present the one-dimensional version of the {Multiscale PNP model (MPNP), and then extend the procedure for higher dimensions}. {We follow the strategy proposed in \cite{astuto2023multiscale}, considering this time two additional equations: one describing the diffusion of positive ions, $c_+$, and one for the electrostatic potential $\varphi$, which couples anions and cations. As a result, we derive three new time-dependent boundary conditions for the system.} We then extend the model to higher dimensions, discussing the additional complexities that arise from the geometry.
\subsection{One-dimensional MPNP model}
\label{sec:MPNP_1D}
\begin{figure}
    \centering 
\begin{minipage}[b]
		{.49\textwidth}
		\centering    
\begin{overpic}[abs,width=\textwidth,unit=1mm,scale=.25]{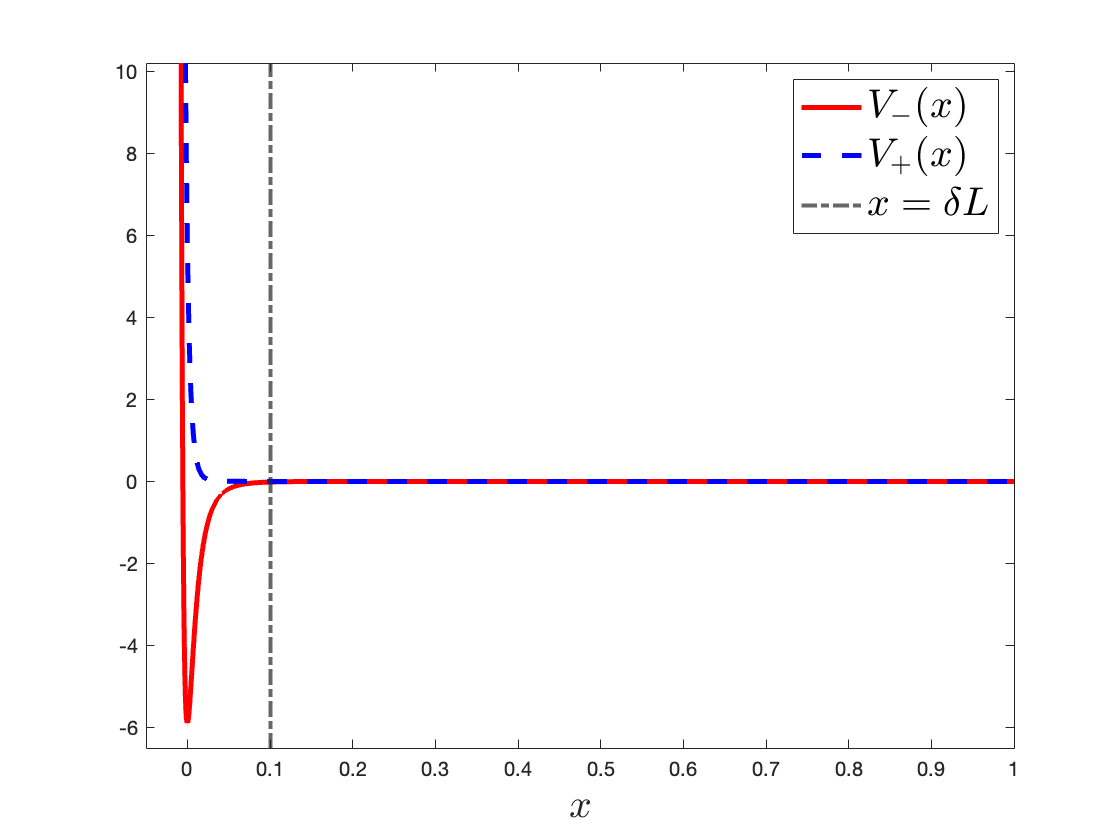}
\put(4,54){(a)}
\put(12,57){$\Omega^\delta_{\rm b}$}
\put(60,57){$\Omega^\delta_{\rm f}$}
\put(18,2){$\delta L$}
\end{overpic}
\end{minipage}
\begin{minipage}[b]
		{.49\textwidth}
		\centering    
\begin{overpic}[abs,width=\textwidth,unit=1mm,scale=.25]{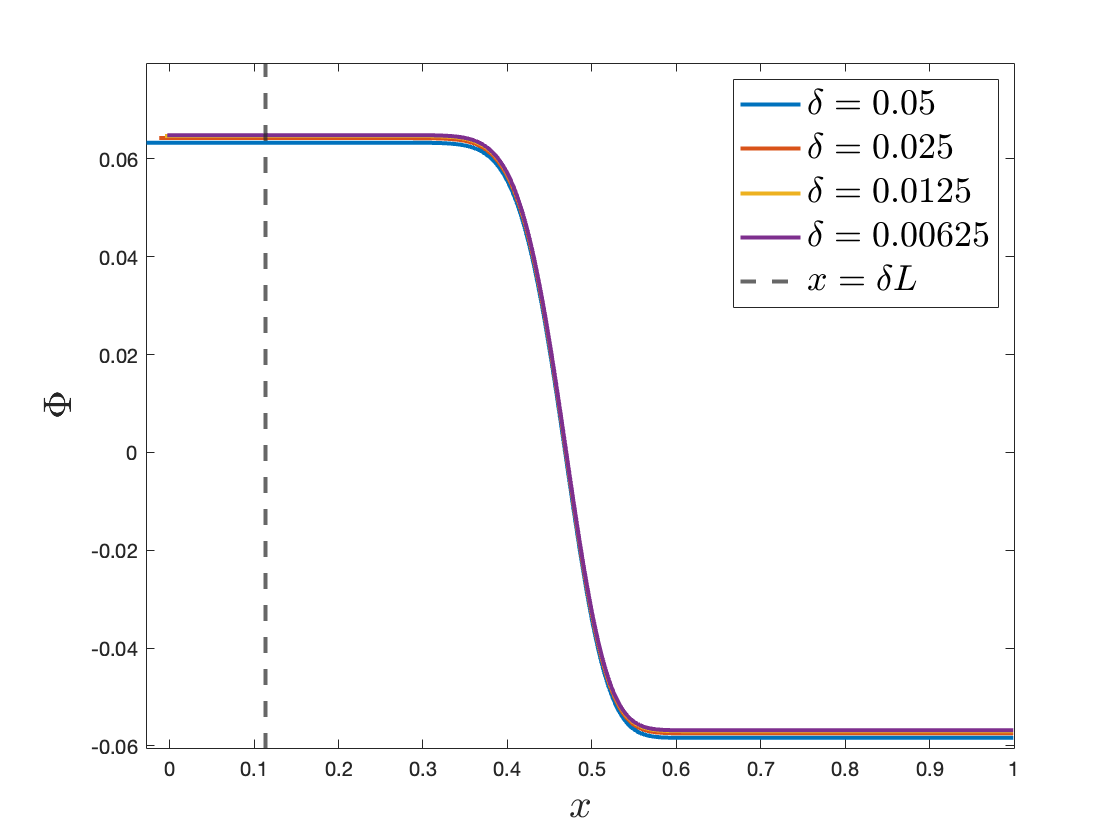}
\put(12,57){$\Omega^\delta_{\rm b}$}
\put(60,57){$\Omega^\delta_{\rm f}$}
\put(18,2){$\delta L$}
\put(4,54){(b)}
\end{overpic}
\end{minipage}
\caption{\textit{{(a): Representation of the external potentials $V_\pm(x)$, for $\delta = 5\times 10^{-2}$. After $x = \delta L$, the contribution of the potentials is negligible. For this reason, the dashed line $x=\delta L$ denotes the right boundary of $\Omega_{\rm b}^\delta$, the domain affected by the potentials. Numerically, we notice that a reasonable value for this quantity is $L=2$. (b) Representation of the Coulumb potential $\Phi$, at final time $t = 1.5$. We emphasize the dashed line $x=\delta L$, for $\delta = 0.05$, to show that $\Phi$ is approximately constant in $\Omega_{\rm b}^\delta$.} }}
    \label{fig:Lpotential}
\end{figure}

\begin{figure}
    \centering 
\includegraphics[width=0.6\textwidth]{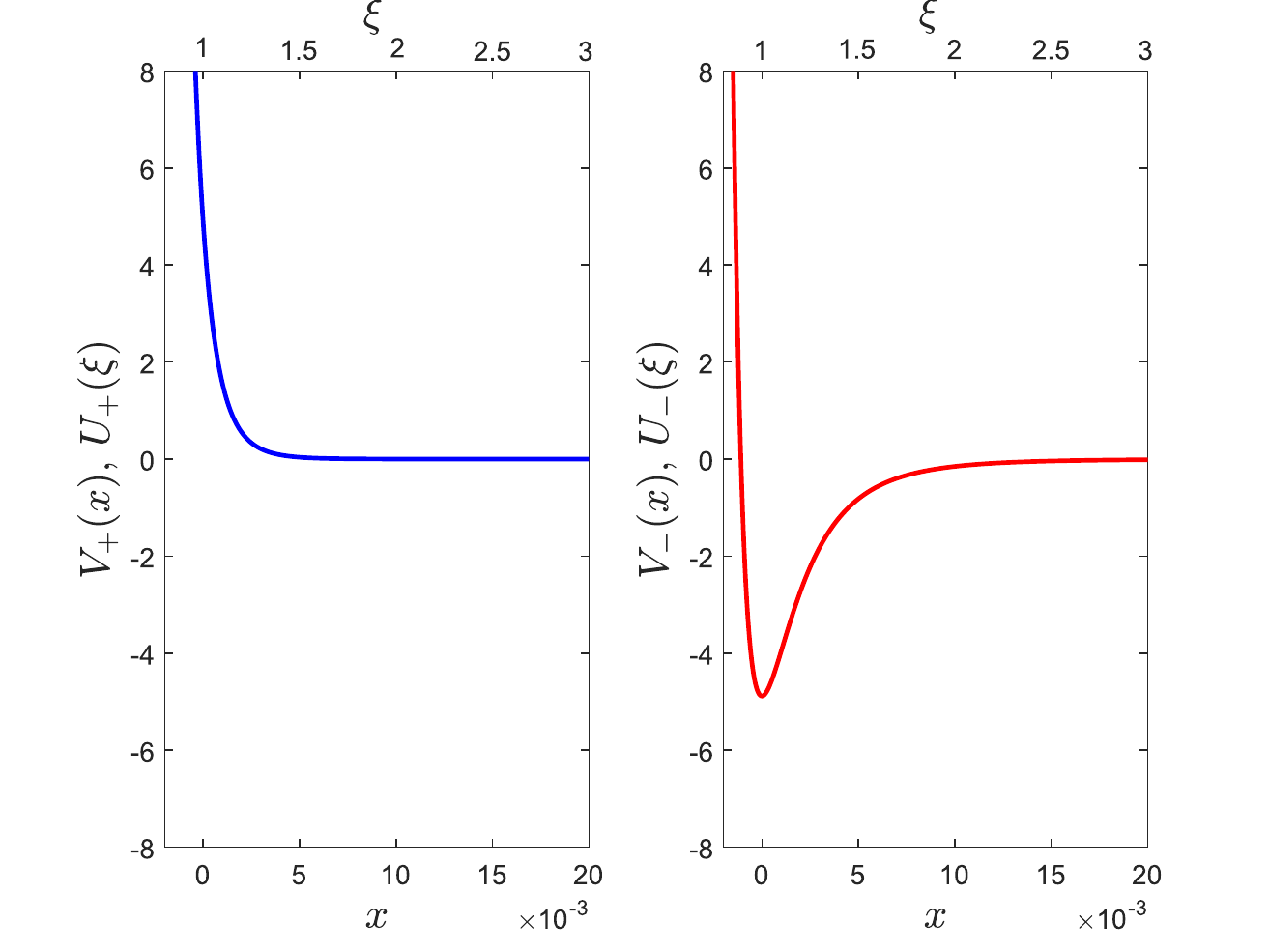}        
    \caption{\textit{Example of the potentials $V_\pm(x)$ and, after a change of variable, the corresponding $U_\pm(\xi)$ for $\delta = 10^{-2}$.}}
    \label{figure_potential_V_1D}
\end{figure}

In this section, we introduce the multiscale model in one dimension. The domain is divided in two regions: $\Omega^\delta = [- \delta ,1] =\Omega^\delta _{\rm b} \cup \Omega^\delta _{\rm f}$. $\,\Omega^\delta_{\rm b} = [- \delta ,\delta L]$ is the region which is affected by the bubble, described by the external potentials $V_\pm$, while $\, \Omega^\delta _{\rm f} = [\delta L,1]$ is the region where the influence of the external potentials is negligible, {say} $V_\pm = 0$; { see Fig.~\ref{fig:Lpotential} (a) for more details. Fig.~\ref{figure_potential_V_1D} displays the potentials restricted to the portion of the domain where their contribution is non negligible. Now
 we write the equations starting from the interval where the influence of the potentials is most relevant, that is} $\Omega^\delta_{\rm b} = [-\delta,\delta L]$
\begin{subequations}
\label{sys:bc_0}
\begin{align}
\label{eq_continuity_1delta}
		\displaystyle \frac{\partial c_\pm}{\partial t} &= -\frac{\partial J_\pm}{\partial x}, \quad {\rm in }  \,\Omega^\delta_{\rm b}\\
		\displaystyle J_\pm &= - D_\pm \left(\frac{\partial c_\pm}{\partial x} + c_\pm \frac{\partial (U_\pm \pm \Phi)}{\partial x}\right), \quad {\rm in }  \,\Omega^\delta_{\rm b} \\ \label{eq_phi2}
		-\varepsilon\frac{\partial^2 \Phi }{\partial x^2} & {= c_+ - c_-}, \quad {\rm in }  \,\Omega^\delta_{\rm b}\\
		\displaystyle \left. J_\pm \right|_{x=-\delta} &= 0\\
		\displaystyle \left. \frac{\partial \Phi}{\partial x}\right|_{x=-\delta} &= 0.
	\label{system_1_delta}
\end{align}
\end{subequations}
In fact we observe that no boundary conditions are necessary in $x=-\delta$ since the potentials barrier prevents from reaching that boundary {(see the behaviour of $V_\pm$ at the left boundary in Fig.~\ref{fig:Lpotential} (a))}.
{As a prototype potential {for the negative ions,} we consider the Lennard-Jones (LJ) potential, which describes attraction at long distances and repulsion at short distances due to Van der Waals and Pauli terms, respectively.} The potential associated with the anions, $V_-(x)$, is given by the following expression
\begin{equation}
\label{eq:Vminus}
	V_-(x) = E\left( \left(\frac{x+\delta}{\delta}\right)^{-12} - 2\left(\frac{x+\delta}{\delta}\right)^{-6} \right),
\end{equation}
where $\delta$ denotes the range of the potential and $E$ represents the depth of the well. For the potential $V_+(x)$, which simulates only a repulsive behavior {with the bubble}, we choose the following expression
\begin{equation*}
	V_+(x) = E\left(\frac{x+\delta}{\delta}\right)^{-12}.
\end{equation*}
Now we write the system for the interval $\Omega^\delta_{\rm f} = [\delta L,1]${, where the effect of the potentials is negligible}
\begin{subequations}
\label{sys:bc_1}
\begin{align}
		\displaystyle \frac{\partial c_\pm}{\partial t } &= -\frac{\partial J_\pm}{\partial x},\quad {\rm in }  \,\Omega^\delta_{\rm f}\\
		\displaystyle J_\pm &= -D_\pm\left(\frac{\partial c_\pm}{\partial x} \pm c_\pm\frac{\partial \Phi}{\partial x}\right),\quad {\rm in }  \,\Omega^\delta_{\rm f}\\
		\displaystyle -\varepsilon\frac{\partial^2 \Phi }{\partial x^2} &= {c_+ - c_-},\quad {\rm in }  \,\Omega^\delta_{\rm f}\\
		\displaystyle \left. J_\pm \right|_{x=1} &= 0\\		
		\displaystyle \left. \frac{\partial \Phi}{\partial x}\right|_{x=1} &= 0. 
\end{align}  
\end{subequations}
The goal of this section is to replace the effect of the bubble region by a suitable boundary condition on the left boundary of $\Omega^\delta_{\rm f}$ {(that coincides with the right boundary of $\Omega^\delta_{\rm b}$)}, as $\delta \to 0$.
If we define $c^{\rm b}_\pm, \Phi^{\rm b}_\pm$ the unknowns inside the bubble region $\Omega_{\rm b}^\delta$ and $c^{\rm f}_\pm,\Phi^{\rm f}_\pm$ the ones in the fluid region $\Omega_{\rm f}^\delta$, the boundary conditions  that we should impose  at the interface are 
{
\begin{subequations}
    \begin{align}
        c^{\rm b}_\pm(L\delta) &= c^{\rm f}_\pm(L\delta) \\ 
        \Phi^{\rm b}_\pm(L\delta) &= \Phi^{\rm f}_\pm(L\delta).
    \end{align}
\end{subequations}}
$L$ is a scaled distance beyond which the potentials are negligible. {In all our calculations we used $L=2$; see Fig.~\ref{fig:Lpotential} (a).}
We follow the strategy that we introduced in \cite{astuto2023multiscale}. Let us assume that the potentials $V_\pm(x) = V_{\pm,\delta}(x)$ depend on $\delta$ but maintain the same functional form, i.e.~there exist two pairs of functions: $\widetilde{U}_\pm\colon[0,L+1] \rightarrow \mathbb{R}$  and $\mathcal{U}_\pm\colon\mathbb{R} \rightarrow \mathbb{R}$ that do not depend on $\delta$ and such that
\begin{align}
\displaystyle \mathcal{U}_\pm(\xi)=\left\{
\begin{array}{l}
+\infty \quad \xi\leq 0\\
\widetilde U_\pm(\xi) \quad \xi \in [0,L+1] \\
0 \quad \quad \xi > L+1, 
\end{array}
\right.
\end{align}
where we adopt a non dimensional form of the potentials, expressed as a function of the rescaled variable 
\begin{equation}
\label{eq:xi_def}
\xi = 1 + x/\delta \in [0,L+1],
\end{equation}
and such that
\begin{subequations}
\begin{align}
\label{expr_U_LJ}
    \widetilde U_-(\xi) &= \nu \left( \xi^{-12} - 2\xi^{-6} \right) \\ \label{eq_U+}
   \widetilde U_+(\xi) &= \nu \, \xi^{-12},
\end{align}
\end{subequations}
$\displaystyle \nu = {E}/{k_BT}$ represents the ratio between the depth of the potential well of the bubble, $E$, and $k_B T$ {(see Fig.~\ref{figure_potential_V_1D} for a representation of the two  external potentials)}. Typical values of $\nu$ should range between 10 and 20, \cite{everett1976adsorption}. {The notation $\widetilde \cdot$ indicates that the functions depend on $\xi$.}

It is clear that the solutions $ c_{\pm,\delta},\Phi_{\pm,\delta}$ and the flux $J_{\pm,\delta}$ depend on $\delta$ as well, thus Eqs.(\ref{eq_continuity_1delta}-\ref{eq_phi2}) can be written as
\begin{subequations}
\begin{align}
\label{eq_flux_eps}
\displaystyle \frac{\partial  c_{\pm,\delta}}{\partial t} + \frac{\partial J_{\pm,\delta}}{\partial x} &=0\quad \rm{in } \, \Omega^\delta_{\rm b}\\
\displaystyle J_{\pm,\delta}&= -D\left(\frac{\partial  c_{\pm,\delta}}{\partial x} +  c_{\pm,\delta}  \, \frac{\partial}{\partial x} \left(U_{\pm,\delta}\pm \Phi\right)\right) \\ 
-\varepsilon \frac{\partial^2 \Phi}{\partial x^2} &= {c_+ - c_-}
\end{align}
\end{subequations}
For $x \in \Omega^\delta_{\rm b}$, we use the scaled variable $\displaystyle\xi \in [0,L+1]$ and consequently, the flux in $\Omega^\delta_{\rm b}$ can be expressed as
\begin{equation} \label{newflux}
\widetilde J_{\pm,\delta} = -D_\pm \frac{1}{\delta}\left(\frac{\partial \widetilde  c_{\pm,\delta}}{\partial \xi}+ \widetilde c_{\pm,\delta} \, \frac{\partial}{\partial \xi} \left(\widetilde U_{\pm,\delta}\pm \widetilde \Phi\right)  \right).  
\end{equation} 
 {Rewriting \eqref{eq_flux_eps} in $\xi$, and using Eq.~\eqref{newflux}, we obtain
\begin{equation}
    \label{c_eps2}
    \frac{\partial  \widetilde  c_{\pm,\delta}}{\partial t} = -\frac{1}{\delta}\pad{\widetilde  J_{\pm,\delta}}{\xi} = \frac{D_\pm}{\delta^2}\frac{\partial}
    {\partial \xi} \left(\pad{\widetilde  c_{\pm,\delta}}{\xi} + \widetilde  c_{\pm,\delta} \frac{\partial}{\partial \xi} \left(\widetilde U_{\pm,\delta}\pm \widetilde \Phi\right) \right),
\end{equation}
}

The range of the scaled variable $\xi$ does not depend on $\delta$, allowing us to assume that $\widetilde c_{\pm,\delta}(\xi,t)$ has the following asymptotic expansion in $\Omega^\delta_{\rm b}$:
\begin{equation}\label{exprho}
{\widetilde c_{\pm,\delta}(\xi,t) =  \widetilde c_\pm^{\,(0)}(\xi,t)+\delta^2 \widetilde c_\pm^{\,(1)}(\xi,t)+
\delta^4 \widetilde c_\pm^{\,(2)}(\xi,t)+\cdots.}
\end{equation} 
 
Inserting the expansion~\eqref{exprho} into Eq.~\eqref{newflux}, we obtain the following expansion for the flux:
\begin{equation}
    \delta \widetilde J_{\pm,\delta}(\xi,t) = \widetilde J_\pm^{\,(0)}(\xi,t) + \delta^2 \widetilde J_\pm^{\,(1)}(\xi,t)  + \delta^4 \widetilde J_\pm^{\,(2)}(\xi,t) + \cdots,
    \label{Jexpand}
\end{equation}
with 
\[
   \widetilde J_\pm^{\,(k)} = -D_\pm\left(\pad{ \widetilde c_\pm^{\, (k)}}{\xi} + \widetilde c_\pm^{\, (k)}\frac{\partial}{\partial \xi} \left(\widetilde U_{\pm,\delta}\pm \widetilde \Phi\right)  \right), \quad k\geq 0.
\]
Using expansion \eqref{Jexpand}  in Eq.~\eqref{c_eps2}, we obtain, to the various orders in $\delta$:
\begin{eqnarray}
    O(\delta^{-2})  :& \displaystyle  \phantom{\quad\frac{\partial  c_\pm^{(k)}}{\partial t} + } \pad{\widetilde J_\pm^{\, (0)}}{\xi} & = 0 ,
    \label{eq:eps0}\\
    O(\delta^{2k})  : & \displaystyle \quad\frac{\partial  \widetilde c_\pm^{\, (k)}}{\partial t} + \pad{\widetilde J_\pm^{\, (k)}}{\xi}
     &  =  0, \quad k\ge 0.\label{eq:epsk}
\end{eqnarray}
Eq.~\eqref{eq:eps0} states that  the lowest order flux $\widetilde J_\pm^{\, (0)}$ is constant. Since we have zero-flux boundary conditions, we deduce that $\widetilde J^{\, (0)}_\pm = 0$, and from its definition we have
\begin{align}
\label{eq:int_xi}
		\displaystyle\frac{1}{\widetilde c^{\, (0)}_\pm}\frac{\partial \widetilde c^{\, (0)}_\pm}{\partial \xi} & = -\frac{\partial}{\partial \xi} \left(\widetilde U_{\pm,\delta}\pm \widetilde \Phi\right).
\end{align}
Integrating Eq.~\eqref{eq:int_xi} between $\xi$ and $L+1$, we have 
\begin{align}
		\ln\left(\frac{\widetilde c^{\, (0)}_\pm(\xi)}{\widetilde c^{\, (0)}_\pm({\xi = L+1})}\right) &= -\left(\widetilde U_\pm(\xi) \pm \widetilde \Phi(\xi) \mp \widetilde \Phi(\xi = L+1)\right)  
\end{align}
whose solution is 
\begin{equation} 
\label{eq_cpm_exp}
    \displaystyle \widetilde c^{\, (0)}_\pm(\xi) = \widetilde c^{\, (0)}_\pm(L+1)  \exp\left( \pm \left( \widetilde \Phi( L+1) - \widetilde \Phi(\xi)\right)\right)\exp\left( -\widetilde U_\pm(\xi)\right)
\end{equation}

{Now we substitute the expressions in Eq.~\eqref{eq_cpm_exp} in Eq.~\eqref{eq_phi2}, and integrate in the interval $[-\delta,\delta L]$, obtaining
\begin{equation}
\label{eq:cB1}
\begin{split}
    -\varepsilon\left.\frac{\partial \Phi }{\partial x}\right|_{x=\delta L} = \frac{1}{m_+} \int_{-\delta}^{\delta L} c_+(\delta L) \exp\left( \left( \Phi(\delta L) -  \Phi(x)\right)\right)\exp\left( - U_+(x)\right) - \\ \frac{1}{m_-} \int_{-\delta}^{\delta L} c_-(\delta L) \exp\left(- \left( \Phi(\delta L) -  \Phi(x)\right)\right)\exp\left( - U_-(x)\right)
\end{split}
\end{equation}
where we used that ${\partial \Phi }/{\partial x} = 0$ if $x = -\delta$ and we omit to specify the zeroth order in the apex to simplify the notation. Defining the following quantities:
\begin{subequations}
\label{eq:Mpm_cB}
\begin{align}
\label{eq:Mpm}
M_\pm & = \delta \int_{0}^{L+1}\exp\left( -\left(\widetilde U_\pm \pm \widetilde \Phi\right) \right)d\xi \\ 
c_\pm^B &= c_\pm(\delta L)\exp(\pm  \Phi(\delta L)),
\end{align}
\end{subequations}
{where the integrand of Eq.~\eqref{eq:Mpm} denotes the classical Boltzmann factor.}
Eq.~\eqref{eq:cB1} can be rewritten as
\begin{align}
        -\varepsilon\left.\frac{\partial \Phi }{\partial x}\right|_{x=\delta L} &  = { M_+{c_+^B} -  M_- {c_-^B}}. 
	\label{phi_2_delta}
\end{align}
}

Here we make some assumptions about the quantities in Eqs.~\eqref{eq:Mpm_cB} to simplify Eq.~\eqref{phi_2_delta}. 
From the expression of $U_+$ in Eq.~\eqref{eq_U+}, we observe that $\displaystyle \exp\left( -\widetilde U_+ \right)$ is bounded by 1, while assuming that $\delta \ll 1$. {Moreover, for sufficiently small values of $\delta$,  the quantity $\exp\left( \mp\widetilde\Phi \right)$ can be considered as approximately constant in the region $x \in \Omega^\delta_{\rm b} = [-\delta,\delta L]$ (see for example Fig.~\ref{fig:Lpotential} (b)). Thus, we obtain} that the term $ M_+  $ is of the order $\delta$, and more precisely 
\begin{equation}
\label{eq:Mp_order_delta}
M_+ \approx  \exp\left( -\widetilde\Phi \right) \delta \int_{0}^{L+1} \exp\left( -\widetilde U_+ \right) \, d\xi < \exp\left(-\widetilde \Phi\right)\,\delta\,(L+1) \to 0 \quad {\rm as } \> \delta \to 0.
\end{equation}
Analogously, for $ M_- $, we can write 
\begin{equation}
\label{eq_Mm}
 M_- \approx  \exp\left( \widetilde\Phi \right) \underbrace{\delta \int_{0}^{L+1} \exp\left( -\widetilde U_- \right) \, d\xi}_{M},
\end{equation}
and we assume that $M = \delta \int_0^{L+1}\exp(-\widetilde U_-) d\xi $ takes a finite value as $\delta \to 0$, as in \cite{astuto2023multiscale}. {In Remark~\ref{remark:Mfixed}, we explain the details of our strategy.}

At the end, we rewrite Eq.~\eqref{phi_2_delta} as follows
\begin{equation}
  -\varepsilon \left. \frac{\partial  \Phi }{\partial x}\right|_{x=L\delta} = O(\delta) - {  M {c_-^B}}
\end{equation}
and if we perform the limit as $\delta \to 0$, with fixed $M$, we obtain the boundary condition for ${ \Phi }$ in $x = 0^+$
\[ \varepsilon \left. \frac{\partial  \Phi }{\partial x}\right|_{x=0}  = { M {c_-^B}  \left( =  M {c_-(x = 0^+)}  \right)}. \]

Now, we look for the boundary conditions on the concentrations $c_\pm$. We start from Eq.~\eqref{eq_continuity_1delta}, integrate in the interval $[-\delta,L\delta]$ and substitute the concentrations with the expressions found in \eqref{eq:Mpm_cB}, obtaining
\begin{equation}
\label{eq_Mpm}
  \frac{\partial }{\partial t} M_\pm c_\pm^B = D_\pm \left(\left.\frac{\partial c_\pm}{\partial x}\right|_{x=L\delta} \pm \left. c_\pm \frac{\partial \Phi}{\partial x}\right|_{x=L\delta}\right)
\end{equation}
where the right hand side of the last equation coincides with the value of the flux $J_\pm$ in $x = L\delta$, and we obtain it using the boundary condition $J(x=-\delta) = 0$.

Regarding the equation for \( c_+ \), we make the same assumptions as before, where \( M_+ \) is of order \( \delta \). Taking the limit as \( \delta \to 0 \), we obtain
\begin{equation}
\label{eq_BC_cp}
    \left.\frac{\partial c_+}{\partial x}\right|_{x=0} + \left. c_+ \frac{\partial \Phi}{\partial x}\right|_{x=0} = 0.
\end{equation}
Now we derive the expression for $c_-$. Substituting Eq.~\eqref{eq_Mm} in Eq.~\eqref{eq_Mpm}, we obtain
\begin{equation}
    M\left.\frac{\partial c_-}{\partial t} \right|_{x=0} = D_- \left(\left.\frac{\partial c_-}{\partial x}\right|_{x=0} - \left. c_-\frac{\partial \Phi}{\partial x}\right|_{x=0}\right).
\end{equation}

Summarizing, in the limit of $\delta \to 0$, with $M$ finite, we obtain the system
\begin{subequations}
\label{eq_final_system_1D}
\begin{align}
\label{eq_final_system_cpm}
		\displaystyle \frac{\partial c_\pm}{\partial t} &= -\frac{\partial J_\pm}{\partial x}, \quad [0,1] 
 \\
  \label{eq_final_system_phi}
		-\varepsilon\frac{\partial^2 \Phi }{\partial x^2} &= { c_+ - c_-}, \quad [0,1]
\end{align}
\end{subequations}
where 
\begin{equation}
      \label{eq_final_system_flux}
		\displaystyle J_\pm = - D_\pm \left(\frac{\partial c_\pm}{\partial x} \pm  c_\pm \frac{\partial \Phi}{\partial x}\right),
\end{equation}
and boundary conditions
\begin{subequations}
\label{eq_final_system_1D_bc}
\begin{align}
 \label{eq_final_system_bc_cp}
	 0 & = \left.\frac{\partial c_+}{\partial x}\right|_{x=0}  \\ 
 \label{eq_final_system_bc_cm}
 \displaystyle \left.M\frac{\partial c_-}{\partial t}\right|_{x=0} & = D_-\left(\left.\frac{\partial c_-}{\partial x}\right|_{x=0} - \left. c_-\frac{\partial \Phi}{\partial x}\right|_{x=0}\right) \\
 \varepsilon \left.\frac{\partial \Phi}{\partial x}\right|_{x=0} & = \left.M c_-\right|_{x=0} \label{system_1_delta_2}  \\ \label{system_1_delta_3}		
		\displaystyle \left. J_\pm \right|_{x=1} &= 0\\
		\displaystyle \left. \frac{\partial \Phi}{\partial x}\right|_{x=1} &= 0
 \label{eq_final_system_bc_phi}
\end{align}
\end{subequations}
In the multiscale model, the boundary conditions that we obtain are independent of the shape of the potentials $V_\pm$ and depend solely on the constant $M$, that contains the relevant physical properties of the system at different scales.

Note that this model does not take into account saturation effects, which are due to the fact that the concentrations are limited by 1, therefore the charge accumulated at the bubble cannot exceed $\delta$. Such effect has been taken into account in the original multiscale model for one carrier \cite{astuto2023multiscale}, and will be considered in a future work {(e.g., \cite{astuto2026gradient})}.

\begin{remark}
\label{remark:Mfixed}
{Here, we briefly recall our strategy to compute the constant $M$. We observe that, if the potential $\widetilde U_-(\xi)$ does not depend on $\delta$, $M \rightarrow 0$ as $\delta \rightarrow 0$ and then the condition \eqref{eq_Mm} reduces to a zero Neumann boundary condition. Therefore, in validating the model, we shall assume that $\delta \to 0$ with $M$ fixed and we explain how to perform such a limit.}

{There are several ways to choose the values of the parameters $M$, $\delta$, and of the function $\widetilde U_-(\xi)$. As we said before, since $M \to 0 $ such as $\delta \to 0$, we may consider a potential $\widetilde U_-(\xi)$ which scales in such a way that $M$ is finite as $\delta$ approaches $\delta_0$, where $\delta_0$ is a typical length of the attractive-repulsive core of the bubble and it is prohibitively small for a detailed calculation. In practice, we can assume that $M$ is non negligible when $\delta\ll l$, where $l$ is a typical length scale of the physical setup. Here we show how to choose the parameters such that $M$ is constant when $\delta \to \delta_0$.} 

{Once we define $\widetilde U_-(\xi)$ in Eq.~\eqref{expr_U_LJ}, the resulting expression for $M$ from Eq.~\eqref{eq_Mm} is
\begin{align}
	M = \delta \mathcal{I}_L(\nu)
	\label{eq_choiseparam}
\end{align}
where 
\begin{equation}
\mathcal{I}_L(\nu) = \int_0^{L+1} \exp\left(-\nu\left( \xi^{-12} - 2\xi^{-6} \right) \right) d\xi.
\label{eq_Iphi_mult}
\end{equation}
Notice that, since the integrand in \eqref{eq_Iphi_mult} is exponentially large in $\nu$ over a finite interval, then 
we can choose $(\delta,\nu)$ such that $M$ is constant by solving the non linear equation for $\nu$, $\mathcal{I}(\nu) = M/\varepsilon$.}
\end{remark}
\subsection{Extension to higher dimensions  MPNP model}

{In this section, we provide a unified treatment of two- and three-dimensional cases, making use of the projection operator on the line  (2D case) or on the surface  (3D case), which delimits the bubble ${\mathcal B}$. 

We assume that the bubble $\mathcal{B}$ is implicitly defined by a level set function $\phi(\vx)$, $\vx\in\mathbb{R}^d$, $d=2,3$, and assume $\phi(\vx)>0$ inside the bubble and $\phi(\vx)<0$ in $\Omega$. For example, $\phi(\vx)$ may be the signed distance from $\Gamma_\mathcal{B}$  (see \cite{Russo2000}), where $\Gamma_\mathcal{B}$ denotes the boundary of the bubble $\mathcal{B}$. The  unit outer normal $\widehat n$ to $\Gamma_\mathcal{B}$ is given by $\widehat n(\vx) = \nabla\phi(\vx)/|\nabla\phi(\vx)|$, 
$\forall \vx$ on $\Gamma_\mathcal{B}$. Note that the $\widehat n$ is naturally defined everywhere in $\Omega\cup\mathcal{B}$. We denote by $n_i, i=1,\ldots,d$ the Cartesian components of $\widehat n$, and we shall denote by $h_{ij} = \delta_{ij}-n_i n_j$ the projection operator on the plane  tangent to $\Gamma_\mathcal{B}$.  Finally, we denote by $\tpar_i \equiv h_{ij}\partial_j$ the gradient operator on the tangent plane to $\Gamma_\mathcal{B}$. We adopt standard Einstein's convention of summation over repeated indices. Furthermore, we assume that the surface $\Gamma_\mathcal{B}$ is smooth, its curvature $\kappa$ times the range of the potential {$\delta$} is much smaller than one, and in the region where the potential $V_\pm$ is non zero, its gradient is parallel to the normal $\widehat{n}$.

As a further simplification of the model, we assume that the bubble potential depends only on the coordinate normal to the surface, and not on the transversal coordinates. This is justified because we assumed $\kappa\delta \ll 1$.

For points near the boundary, we decompose the flux vector $J_{\pm,i},i=1,\ldots,d$ into a normal and tangential component, i.e.
\begin{equation}
\label{eq_components_flux}
        J_{\pm,i} = J^n_{\pm,i} + J_{\pm,i}^{\tau}
\end{equation}
where $J_{\pm,i}^n = n_iJ^n_{\pm} = n_in_jJ_{\pm,j}$ and $J_{\pm,i}^{\tau} = h_{ij} J_{\pm,j}$.

The divergence of the flux is given by:
\begin{equation}
\label{eq_flux_i}
    \nabla \cdot \vec{J}_\pm = \partial_i J_{\pm,i}
\end{equation}
If we substitute Eq.~\eqref{eq_components_flux} in Eq.~\eqref{eq_flux_i} we obtain
\begin{equation}
    \nabla \cdot \vec{J}_\pm = \partial_i \left( n_i J^n_{\pm} \right) + \partial_i \left( h_{i j} J_{\pm,j} \right).
\end{equation}

Now the first term in the last equation can be written as
\begin{equation}
\label{eq_n_i}
    \partial_i \left( n_i J^n_{\pm} \right) = n_i \partial_i J^n_{\pm} + J^n_{\pm} \underbrace{\partial_i n_i}_{\chi_{ii}} 
\end{equation}
where $\chi_{ij} = \partial_in_j$ denotes the second fundamental form of the surface  $\Gamma_\mathcal{B}$, and its trace $\chi_{ii}$ gives the curvature of the surface $\Gamma_\mathcal{B}$:
\begin{equation}
    \chi_{ii} = \kappa = \left\{
\begin{array}{rl} \displaystyle
1/R &  {\rm in}\, 2D \\
\kappa_1 + \kappa_2 &  {\rm in} \, 3D, 
\end{array}
\right.
\end{equation}
where $R$ is the local radius of curvature,  $\kappa_1$ and $\kappa_2$ denote the Gauss principal curvatures.

Since the gradient of the potential $V_\pm$ is parallel to the normal, then one has
\begin{equation}
    \label{eq:Jtau}
    J_{\pm,i}^{\tau} = -D_\pm \left( \tpar_i c_\pm \pm c_\pm \tpar_i \Phi \right) = -D_\pm \left( h_{ij}\partial_j c_\pm \pm c_\pm  h_{ij}\partial_j \Phi \right)
\end{equation}
At this point, we substitute the quantities (\ref{eq_n_i}-\ref{eq:Jtau}) in Eq.~\eqref{eq_flux_i} and the expression for the divergence of the flux becomes:
\begin{equation}
    \nabla \cdot \vec{J}_\pm = n_i \partial_i J^n_{\pm} + J^n_{\pm}\chi_{ii} + \tpar_i J_\pm^\tau 
\end{equation}
where 
\[ \tpar_i J^\tau_\pm = - D_\pm \left( \Delta_\perp c_\pm \pm \tpar_i \left( c_\pm \tpar_i \Phi \right)\right) \]
and $\Delta_\perp c_\pm= \tpar_i \tpar_i c_\pm= h_{ij}\partial_i\partial_j c_\pm- \chi_{ii}n_j\partial_j c_\pm$, 
that denotes the surface Laplacian of the concentrations. 

The evolution equation for the concentration near the bubble surface is then given by
\begin{equation}
\label{eq:dcpm_dt}
    \pad{c_\pm}{t} = -\pad{J^n_{\pm}}{n} - J^n_{\pm} \chi_{ii} - \tpar_i {J^\tau_\pm} 
\end{equation}
We follow the procedure adopted in one dimension, but we start integrating the drift-diffusion equations because of the complexity of these equations. Let $\vec{x}$ be a point in $\Gamma_\mathcal{B}$ and integrate Eq.~\eqref{eq:dcpm_dt} along the normal direction:
\begin{align} \nonumber
	\frac{\partial}{\partial t}\int_{-\delta}^{\delta L} c_\pm(\vec{x}+r\widehat n)\,dr = & -J^n_{\pm}(\vec{x}+L\delta\widehat n) 
	- \int_{-\delta}^{\delta L} \chi_{ii}  J^n_{\pm} (\vec{x}+r\widehat n) \, dr  \\  
	& +  D_\pm \Delta_\perp \int_{-\delta}^{\delta L}  c_\pm(\vec{x}+ r \, \widehat n)\,dr  \pm  D_\pm \int_{-\delta}^{\delta L}  \tpar_i \left( c_\pm(\vec{x}+r \, \widehat n) \tpar_i \Phi(\vec{x}+r \, \widehat n) \right) \,dr 
\end{align}
{ 
where we dropped the terms $J^n_{\pm}(\vec{x}-\delta\,\widehat n), J^\tau_{\pm}(\vec{x}-\delta\,\widehat n)$ which vanish because of the repulsive core. 

Notice that the term 
$\int_{-\delta}^{\delta L} \chi_{ii}  J^n_{\pm} (\vec{x}+r\widehat n) \, dr  = (1+L)\,\delta\,\langle \chi_{ii} J^n_{\pm} \rangle$ can be neglected since we assume that the range of the potential is much smaller than the radius of curvature (2D) or to the inverse of the mean curvature (3D).
}

Now, we consider a change of variable from $r$ to $\xi = r/\delta$ in the integration intervals, and obtain 
\begin{align} \nonumber
	\delta\frac{\partial}{\partial t}\int_{-1}^{L} c_\pm(\vec{x}+\delta\xi \,\widehat n)\,d\xi 
    = & -J^n_{\pm}(\vec{x}+\delta L\,\widehat n) + \delta D_\pm \Delta_\perp \int_{-1}^{L}  c_\pm(\vec{x}+\delta \xi \, \widehat n)\,d\xi \\ \label{eq:delta_Jn_Jtau}
    &  \pm  \delta D_\pm  \int_{-1}^{L}  \tpar_i \left( c_\pm(\vec{x}+\delta \xi \, \widehat n) \tpar_i \Phi(\vec{x}+\delta \xi \, \widehat n) \right) \,d\xi.
\end{align}
Let us consider {  the normal components of the flux $J^n_{\pm}$}:
\begin{equation}
\label{eq_J_n}
	J^n_{\pm} = -D_\pm\left(\frac{\partial  c_\pm}{\partial {n}}+\frac{\partial \left(U_\pm \mp \Phi\right)}{\partial {n}} c_\pm\right).
\end{equation}
{As in the one-dimensional case, we notice that} there is also a dependence on $\xi $ of the solution inside the layer, such that $ c_\pm= c_\pm(\vec{x} + \delta \xi \, \widehat n,t) = \widetilde{c}_\pm(\xi,t)$, thus we can put $1/\delta$ as common factor:
\begin{equation}
	J^n_{\pm} = -\frac{D_\pm}{\delta}\left(\frac{\partial \widetilde  c_\pm}{\partial \xi}+ \frac{\partial (\widetilde U_\pm \mp \widetilde \Phi)}{\partial \xi} \widetilde c_\pm\right).
\end{equation}
Following the same argument that we used in one dimension to derive Eq.~\eqref{exprho}, we perform a formal expansion in $\delta$ of the solution. To the lowest order in $\delta$, {  for the normal component of the flux, we have $\displaystyle {J_{\pm}^{n,(0)}} = 0$, } which implies 
\begin{equation}
\label{bcbubble}
	\frac{\partial \widetilde c_\pm^{\, (0)}}{\partial \xi} + \frac{\partial (\widetilde U_\pm \mp \widetilde \Phi)}{\partial \xi} \widetilde c_\pm^{\, (0)} = 0.
\end{equation}
{Therefore, we obtain {  the expression for the concentration to order zero in $\delta$:}
\begin{equation}
\label{eq:cpm_limit}
	c_\pm^{(0)}(\vec{x}+\delta \xi \,\widehat n)= c^B_\pm  \exp\left(-U_\pm(\vec{x}+\delta \xi \,\widehat n)\mp \Phi(\vec{x}+\delta \xi \,\widehat n) \right),
\end{equation}
where, as in one dimension, 
\[  c^B_\pm  = c_\pm^{(0)}(\vec{x}+\delta L\, \widehat n)\exp\left( \pm \Phi(\vec{x}+\delta L\, \widehat n) \right).\]
Here, we integrate Eq.~\eqref {eq:cpm_limit} in the bubble layer, to calculate the quantity of entrapped ions at the surface of the trap, as follows
\begin{equation}
\label{eq:Mm_2D}
    \delta \int_{-1}^{L} c_\pm^{(0)}(\vec{x}+\delta \xi \,\widehat n) \,d\xi \approx c^B_\pm M_\pm
\end{equation}
where, for anions, it holds the following expression for $M_-$
\begin{equation}
\label{eq:M_2D}
    M_- \approx \exp\left(\Phi(\vec{x}+\delta L \,\widehat n) \right) M
\end{equation}
and $M = \delta \int_{-1}^{L}   \exp\left(-U_-(\vec{x}+\delta \xi \,\widehat n) \right) d\xi.$ For cations, we apply the same considerations as in the one-dimensional case, concluding that $M_+ \to 0$ as $\delta $ goes to zero in Eq.~\eqref{eq:Mp_order_delta}. 
}

At this point, we analyze each term in Eq.~\eqref{eq:delta_Jn_Jtau} (starting with the one coming from $J_\pm^\tau$), and  conduct distinct analyses for the concentration of positive $c_+$ and negative $c_-$ ions. For very small values of $\delta$, we adopt approximations in Eqs.~(\ref{eq:Mm_2D}-\ref{eq:M_2D}) in the left hand side of Eq.~\eqref{eq:delta_Jn_Jtau}. For the cations, we obtain
\begin{align} 
\label{eq:c_Mpm}
 \delta \pad{}{t} \int_{-1}^{L}  c_+(\vec{x}+\delta \xi \, \widehat n)\,d\xi \approx   \pad{}{t}\left(M_+ c^B_+\right) \to 0\quad  \mbox{\rm because of 
    Eq.~\eqref{eq:Mp_order_delta}.}
\end{align}
Likewise, the second and last terms on the right hand side in Eq.~\eqref{eq:delta_Jn_Jtau} vanish as $\delta\to 0$, yielding the following boundary condition for $c_+$ at the bubble surface:
\[  
  \frac{\partial  c_+}{\partial n}  = 0 \quad \rm{ on }\, \Gamma_\mathcal{B}.
\]
For the negative ions, Eq.~\eqref{eq:delta_Jn_Jtau} is more complex, and we need a more sophisticated approximation. 
First, the left hand side of Eq.~\eqref{eq:delta_Jn_Jtau} becomes  
\begin{align}
\label{eq:appr_dt}
    \delta \pad{}{t} \int_{-1}^{L}  c_-(\vec{x}+\delta \xi \, \widehat n)\,d\xi \approx   \pad{}{t}\left(M_- c^B_-\right) \approx M \pad{}{t}\left(\exp{(\Phi(\vec{x}+\delta L \,\widehat n))} \,c^B_-\right) = M \pad{c_-}{t},
\end{align}
where, the last approximation comes from Eqs.~\eqref{eq:cpm_limit},\eqref{eq:M_2D}.

Secondly, we note that \( M \), as defined in \eqref{eq:M_2D}, is independent of the location on the surface. This follows from our earlier assumption in the subsection that the potential \( U_- \), which defines the bubble, depends only on \( r \).
In this way, the second term of the right hand side of Eq.~\eqref{eq:delta_Jn_Jtau} becomes  
\begin{equation}
\label{eq:appr_Delta}
\delta D_- \Delta_\perp \int_{-1}^{L}  c_-(\vec{x}+\delta \xi \, \widehat n)\,d\xi  \approx  D_- M \Delta_\perp c_- 
\end{equation}
Now, we approximate the last term of the right hand side of Eq.~\eqref{eq:delta_Jn_Jtau}, as follows
\begin{eqnarray}
\label{eq:int_tparPhi}
    \delta D_-  \int_{-1}^{L}  \tpar_i \left( c_-(\vec{x}+\delta \xi \, \widehat n) \tpar_i \Phi(\vec{x}+\delta \xi \, \widehat n) \right) \,d\xi \approx \delta D_-  \tpar_i \left( \tpar_i \Phi(\vec{x}+\delta L \,\widehat n)  \int_{-1}^{L}   c_-(\vec{x}+\delta \xi \, \widehat n)  \,d\xi \right)
\end{eqnarray}
 where we consider that $\tpar_i$ is a differentiation in the direction orthogonal to the normal $\widehat n$, and
 since $\delta $ is very small, we assume $\tpar_i\Phi(\vec{x}+\delta \xi \, \widehat n)$ to be essentially constant inside the integral. 
 Making use of approximations in Eqs.~(\ref{eq:Mm_2D}-\ref{eq:M_2D}), we write 
\begin{eqnarray}
\label{eq:int_tparPhi3}
 \delta D_-  \tpar_i \left( \tpar_i \Phi(\vec{x}+\delta \xi \,\widehat n)  \int_{0}^{L+1}   c_-(\vec{x}+\delta \xi \, \widehat n)  \,d\xi \right) \approx D_- \tpar_i\left(M_- c_-^B\, \tpar_i \Phi  \right)\approx D_- M \tpar_i\left( c_- \,  \tpar_i \Phi \right),
\end{eqnarray}
where, the last approximation comes from Eq~\eqref{eq:M_2D}.

Now, making use of Eqs.~\eqref{eq:appr_dt}, \eqref{eq:appr_Delta}, \eqref{eq:int_tparPhi3} in  Eq.~\eqref{eq:delta_Jn_Jtau} in the limit of $\delta \to 0$, we obtain the boundary condition for $c_-$ at the bubble surface
\begin{align} 
	M \pad{}{t} c_-  =  D_- \left( \frac{\partial c_-}{\partial n} +
	  M \left( \Delta_\perp  c_-   + \tpar_i\left( c_- \, \tpar_i \Phi \right) \right) \right). 
\end{align}
At the end, in two- and three-dimensions, the system of equations is
\begin{subequations}
\label{system_multiscale_all}
\begin{align}
	\label{system_multiscale}
	\displaystyle \frac{\partial c_\pm}{\partial t} &= D_\pm \Delta c_\pm \pm \nabla\cdot \left( c_\pm \nabla \Phi \right) \quad \rm{in} \> \Omega \\ 
    -\varepsilon \Delta \Phi &= {  {c_+} - {c_-}} \quad \rm{in} \> \Omega
\end{align}
\end{subequations}
together with the following boundary conditions for $c_\pm$ and $\Phi$:
\begin{subequations}
\label{system_multiscale_bc}
\begin{align}
	\displaystyle \nabla c_\pm \cdot \widehat n &= 0 \quad \rm{ on } \, \partial \Omega\setminus \Gamma_\mathcal{B} \\ \label{bc_multiscale_cp}
	  \displaystyle \nabla c_+ \cdot \widehat n &= 0 \quad \rm{ on }\, \Gamma_\mathcal{B}
 \\ \label{bc_multiscale_cm}
	\displaystyle M\frac{\partial c_-}{\partial t} &= D_-\left(\frac{\partial  c_-}{\partial n} + M \left( \Delta_\perp c_- +  \tpar_i\left( c_-   \tpar_i \Phi \right)\right) \right) \quad \rm{ on }\, \Gamma_\mathcal{B} \\ \nonumber
 \nabla \Phi \cdot \widehat n & = 0 \quad \rm{ on } \, \partial \Omega\setminus \Gamma_\mathcal{B} \\ 
 \varepsilon \nabla \Phi \cdot \widehat n & = Mc_- \quad \rm{ on } \, \Gamma_\mathcal{B}.
\end{align}
\end{subequations}
Notice that in two-dimensions the Laplace-Beltrami operator reduces to the second derivative with respect to the arclength of the boundary\footnote{Here by $\partial/\partial \tau$ 
we denote the derivative on $\Gamma_\mathcal{B}$, i.e.\ the derivative along the arclength that parametrizes the curve, likewise
$\partial^2/\partial \tau^2$ denotes the second derivative along $\Gamma_\mathcal{B}$, not the second derivative along the tangent direction. }
\begin{equation}
    \Delta_\perp c_\pm= \frac{\partial^2 c_\pm}{\partial \tau^2},
    \label{eq:BL2D}
\end{equation}
while $\tpar_i$ reduces to $\partial/\partial \tau$.

\section{One-dimensional discretization and validation for $\protect \delta \to 0$}
\label{sec:discr1D}
In this section, we discretize the one-dimensional problem with a second order accurate discretization in space and time. We further validate the two-species MPNP system in the limit $\delta \to 0$, demonstrating that the difference between the two approaches is an infinitesimal of order one in $\delta$. 
\subsection{Space discretization in one dimension}
The one dimensional domain is $\Omega^\delta = [- \delta ,1].$ 
The computational domain $\Omega^\delta_h$ is a discretization of $\Omega^\delta$ obtained by a uniform Cartesian mesh with spatial step $ h:   h  N_x = 1 +  \delta,\, N_x\in\mathbb N$. The concentrations $c_{\pm,i}$ are defined at the center of the cells with $x_i=-\delta  + (i-1/2) h \in\{1,\dots,N_x\} $. We choose a cell centered discretization in order to guarantee the exact conservation of the  
{total integral of the solute $v=\sum_i c_{\pm,i} h$}, which is a consequence of the zero boundary condition for the flux. The scheme is second order accurate and it is stable even in presence of a drift term, provided the following condition for the so called {mesh P\'{e}clet number}, pec, is satisfied \cite{Wesseling2023600}:
\begin{equation}
\label{eq_pec}
	{\rm pec} := \max_x\left|{\partial_x U_\pm}\right| h < {2}.
\end{equation} 
The full one dimensional problem~\eqref{eq_full_model_adim} is then discretized in space, leading to the following system of Differential-Algebraic equations:
\begin{subequations}
\label{eq_semi_discrete_full}
\begin{eqnarray}
\label{eq_semi_discrete_eq_c}
	\td{c_{\pm,h}}{t} &=& \left(\mathbb L^{\rm 1D} + \mathbb D^{\rm 1D}_\pm[U_\pm] \right)  c_{\pm,h} \pm \mathbb D^{\rm 1D}_\pm[c_{\pm,h}] \Phi_h, \quad c^0_{\pm,h} = c_{\pm,h}(t = 0) ,\\ \label{eq_semi_discrete_eq_phi}
 -\varepsilon \mathbb L^{\rm 1D}\, \Phi_h &=& {  {c_{+,h}} - {c_{-,h}}}
\end{eqnarray}
\end{subequations}
where $\mathbb L^{\rm 1D},\mathbb D^{\rm 1D}\,\left[U_{\pm,h}\right] $ and $\mathbb D^{\rm 1D}\,\left[ c_{\pm,h}\right] $ are ${N_x} \times {N_x}$ matrices representing the discretization of the space derivatives, as follows
\begin{eqnarray} 
\label{eq_Lh}
\mathbb L^{\rm 1D}\, c_{\pm,h}\Big|_{j} &=& D_\pm \left( \frac{c_{\pm,j+1}  + c_{\pm,j-1} - 2c_{\pm,j}}{h^2} \right)
\\ \nonumber \mathbb D^{\rm 1D}_\pm\,\left[ U_{\pm,h}\right] \,c_{\pm,h}\Big|_{j} & = & D_\pm \frac{(U_{\pm,j+1}-U_{\pm,j})\,c_{\pm,j+1} + (U_{\pm,j+1}-2U_{\pm,j}+U_{\pm,j-1})\,c_{\pm,j} + (U_{\pm,j-1}-U_{\pm,j})\,c_{\pm,j-1}}{2h^2}
\\ \nonumber \mathbb D^{\rm 1D}_\pm\,\left[ c_{\pm,h}\right] \,\Phi_{h}\Big|_{j} & = & D_\pm 
\frac{(c_{\pm,j+1}+c_{\pm,j})\,\Phi_{j+1} - (c_{\pm,j+1}+2c_{\pm,j}+c_{\pm,j-1})\,\Phi_{j} + (c_{\pm,j}+c_{\pm,j-1})\,\Phi_{j-1}}{2h^2}
.
\end{eqnarray} 
The numerical solution at time $t$ is represented by the vector $c_{\pm,h}$, whose components $c_{\pm,i}(t)$ are approximations of the exact solution on the grid points of $\Omega_h$,  i.e.\
$c_{\pm,i}(t)\approx c_\pm(x_i,t)$. The quantities $ c^0_{\pm,h}$ and $ U_{\pm,h}$ are analogously defined. 

For the multiscale model (\ref{eq_final_system_1D}-\ref{eq_final_system_1D_bc}), the domain is $\Omega^0 = [0,1]$, while the discrete computational domain is $\Omega^0_h$, that we discretize as follows: $x_i=(i-1/2) h,\, i \in\{0,\dots,N\},  h  N = 1$. The semidiscrete numerical scheme is the following
\begin{subequations}
\label{eq_semi_discrete_mul}
\begin{eqnarray}
\label{eq_semi_discrete_eq_c_mul}
	\td{c_{\pm,h}}{t} &=& \mathbb L^{\rm 1D} c_{\pm,h} \pm \mathbb D^{\rm 1D}_\pm[c_{\pm,h}] \Phi_h, \quad c^0_{\pm,h} = c_{\pm,h}(t = 0),\\ \label{eq_semi_discrete_eq_phi_mul}
 -\varepsilon \mathbb L^{\rm 1D}\, \Phi_h &=& { {c_{+,h}} - {c_{-,h}}}.
\end{eqnarray}
\end{subequations}

To close the system we need to discretize the boundary conditions. 
To this purpose we add two ghost cells, one to the left of the boundary $x=0$, and call $c_{\pm,0}$ the value at the center of such a cell, and one to the right of the right boundary, and denote by $c_{\pm,N+1}$ their values.  To second order accuracy, the value of $c_-$ at $x=0$ is given by $c_{-,1/2}\approx (c_{-,0}+c_{-,1})/2$. Notice that boundary conditions 
\eqref{eq_final_system_bc_cp}, and (\ref{system_1_delta_3}-\ref{eq_final_system_bc_phi}) impose that  $c_{+,0} = c_{+,1}, \, c_{+,N+1} = c_{+,N}, \, c_{-,N+1} = c_{-,N}$ and $\Phi_{N+1} = \Phi_N$.

The discretization of \eqref{eq_semi_discrete_eq_c_mul} at the left boundary reads
\begin{equation}\label{eq:c_ghost}
    \frac M 2 \left(\td{c_{-,0}}{t} + \td{c_{-,1}}{t} \right)= D_-\left( \frac{c_{-,1}-c_{-,0}}{h} - \frac{c_{-,1} + c_{-,0}}{2} \frac{\Phi_1-\Phi_0}{h}  \right).
\end{equation}
Now we discretize the boundary condition for $\Phi$, that at the left boundary becomes 
\begin{equation}
\label{eq:Phi_ghost}
    \varepsilon \frac{\Phi_1 - \Phi_0}{h}  =  M \frac{c_{-,0} + c_{-,1} }{2}.
\end{equation}
Eqs.~(\ref{eq_semi_discrete_mul}-\ref{eq:Phi_ghost}) constitute a system of $3N+2$ equations  for the $3N+2$ unknowns $c_{+,i},\,  i = 1,\cdots,N, \> c_{-,i}, \, i = 0,\cdots,N$ and $\Phi_i, \, i = 0,\cdots,N$.

\begin{figure}
    \centering
\begin{minipage}[b]
		{.32\textwidth}
		\centering
	\begin{overpic}[abs,width=\textwidth,unit=1mm,scale=.25]{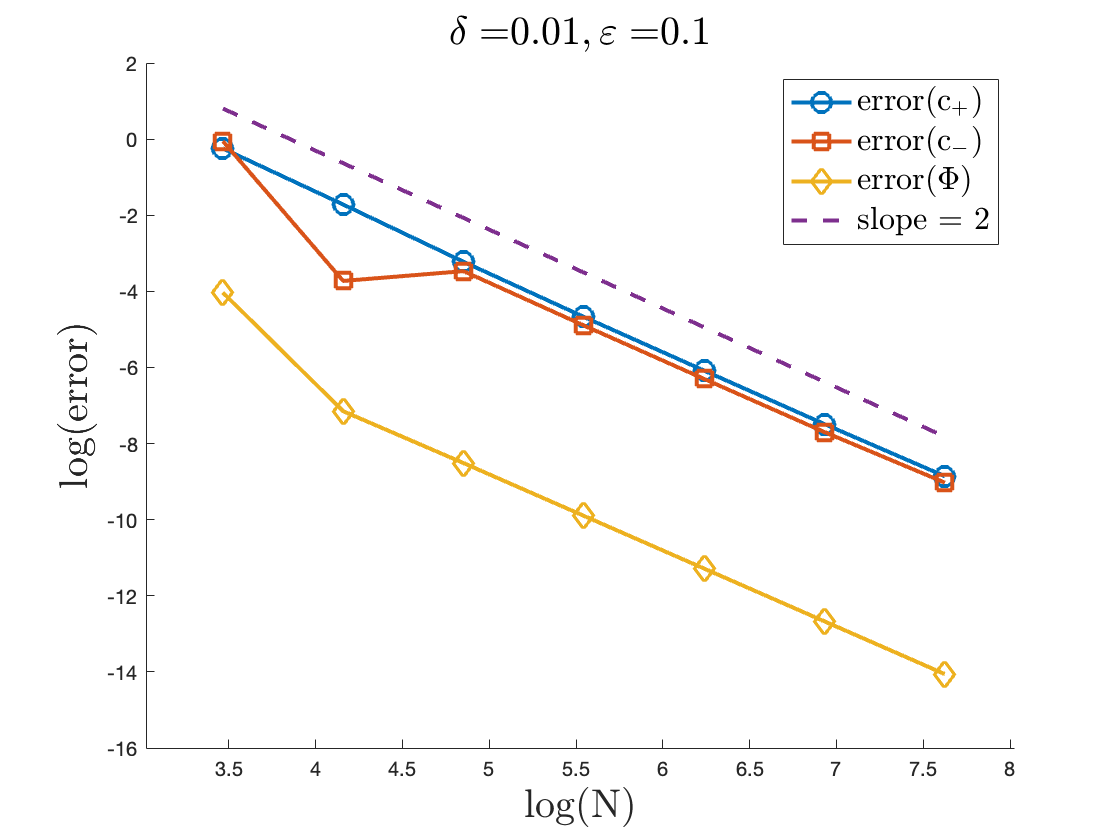}
\put(1,37){(a)}
\end{overpic}
\end{minipage}         
\begin{minipage}[b]
		{.32\textwidth}
		\centering
	\begin{overpic}[abs,width=\textwidth,unit=1mm,scale=.25]{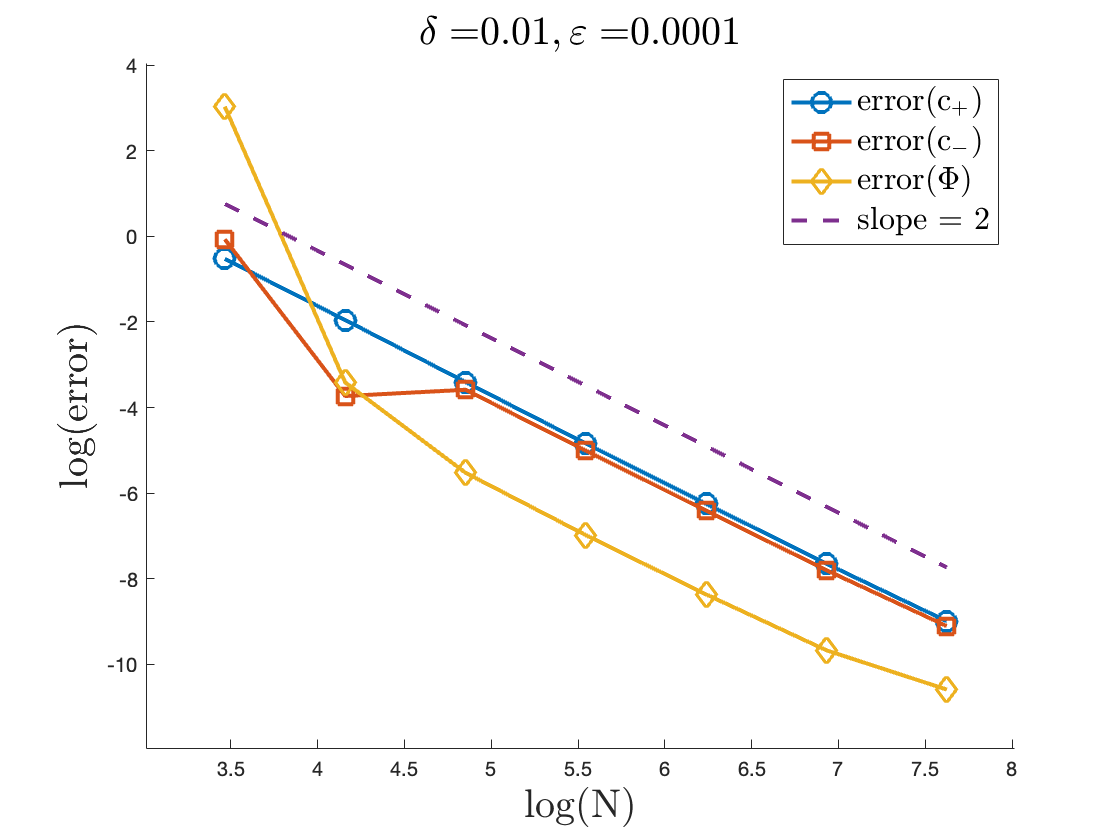}
\put(1,37){(c)}
\end{overpic}
\end{minipage}         
\begin{minipage}[b]
		{.32\textwidth}
		\centering
	\begin{overpic}[abs,width=\textwidth,unit=1mm,scale=.25]{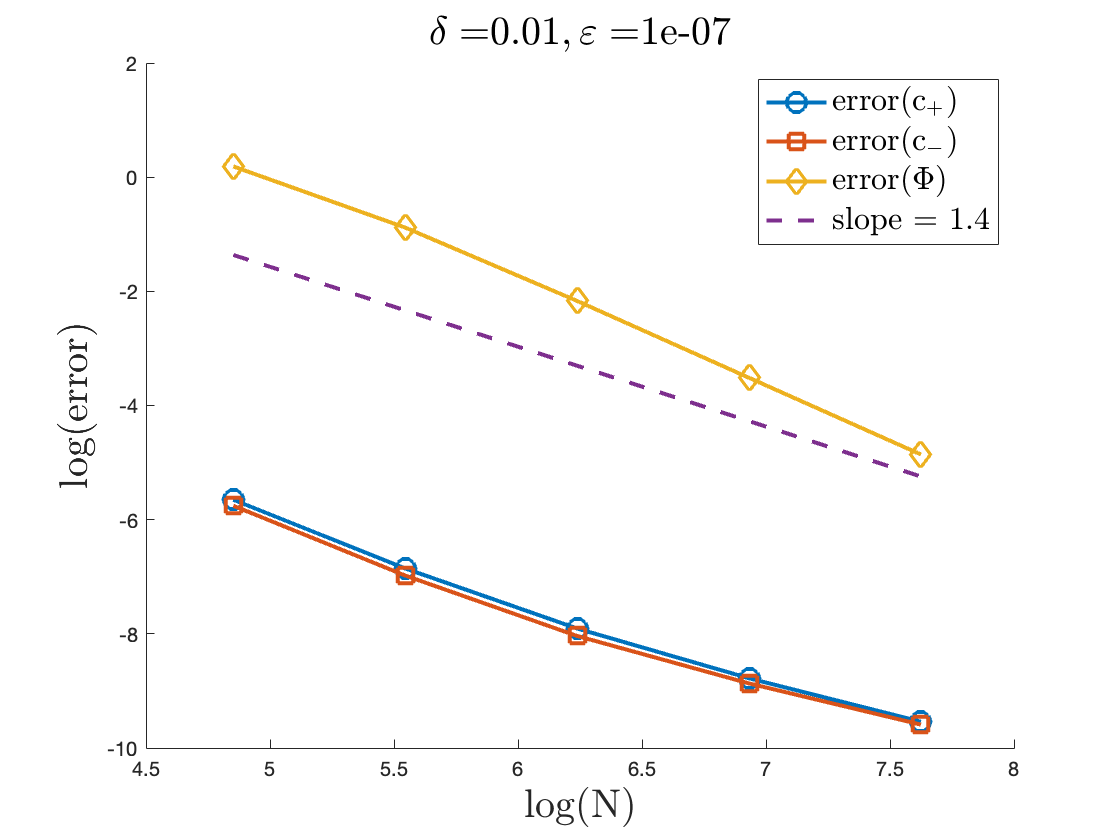}
\put(1,37){(e)}
\end{overpic}
\end{minipage} 
\begin{minipage}[b]
		{.32\textwidth}
		\centering
	\begin{overpic}[abs,width=\textwidth,unit=1mm,scale=.25]{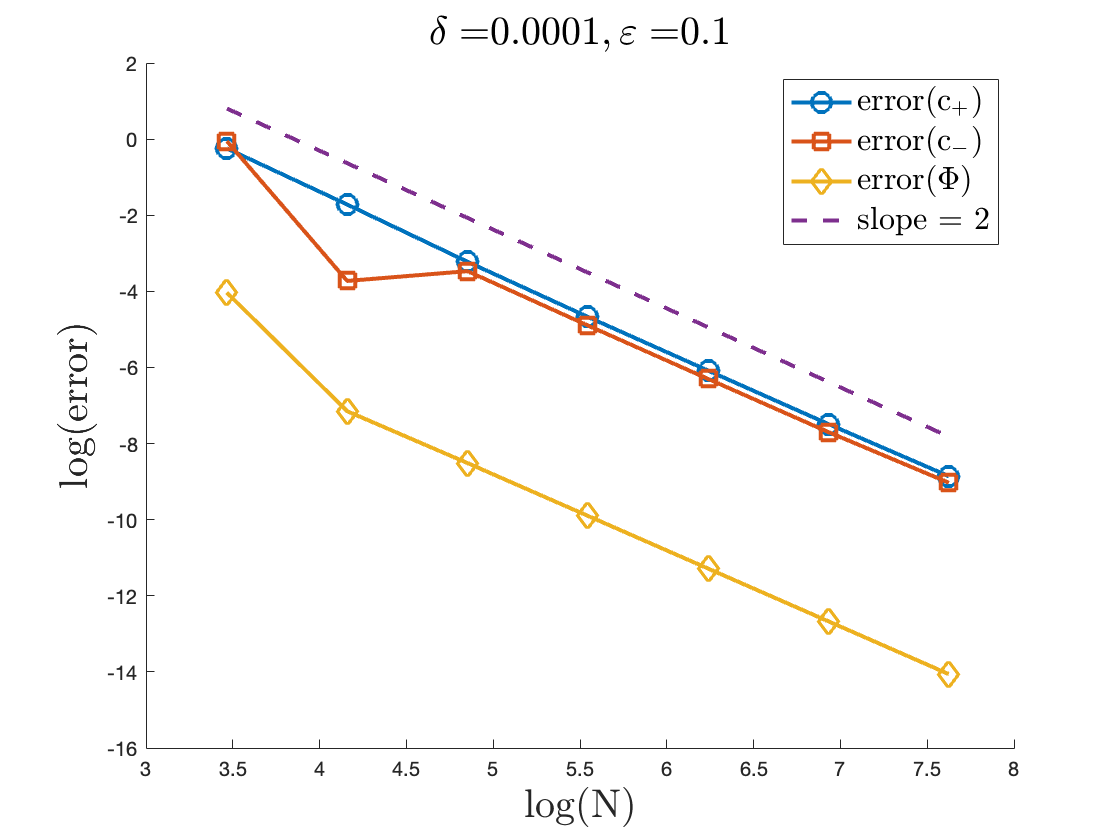}
\put(1,37){(b)}
\end{overpic}
\end{minipage}        
\begin{minipage}[b]
		{.32\textwidth}
		\centering
	\begin{overpic}[abs,width=\textwidth,unit=1mm,scale=.25]{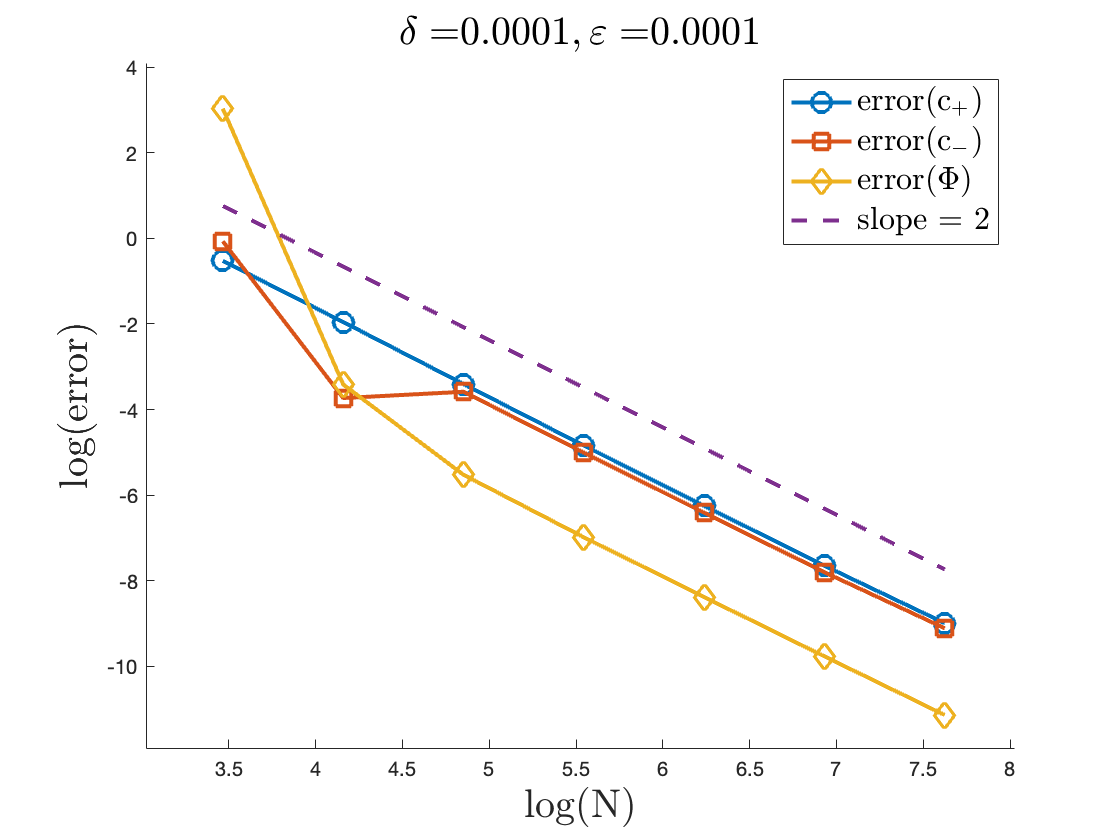}
\put(1,37){(d)}
\end{overpic}
\end{minipage}             
\begin{minipage}[b]
		{.32\textwidth}
		\centering
	\begin{overpic}[abs,width=\textwidth,unit=1mm,scale=.25]{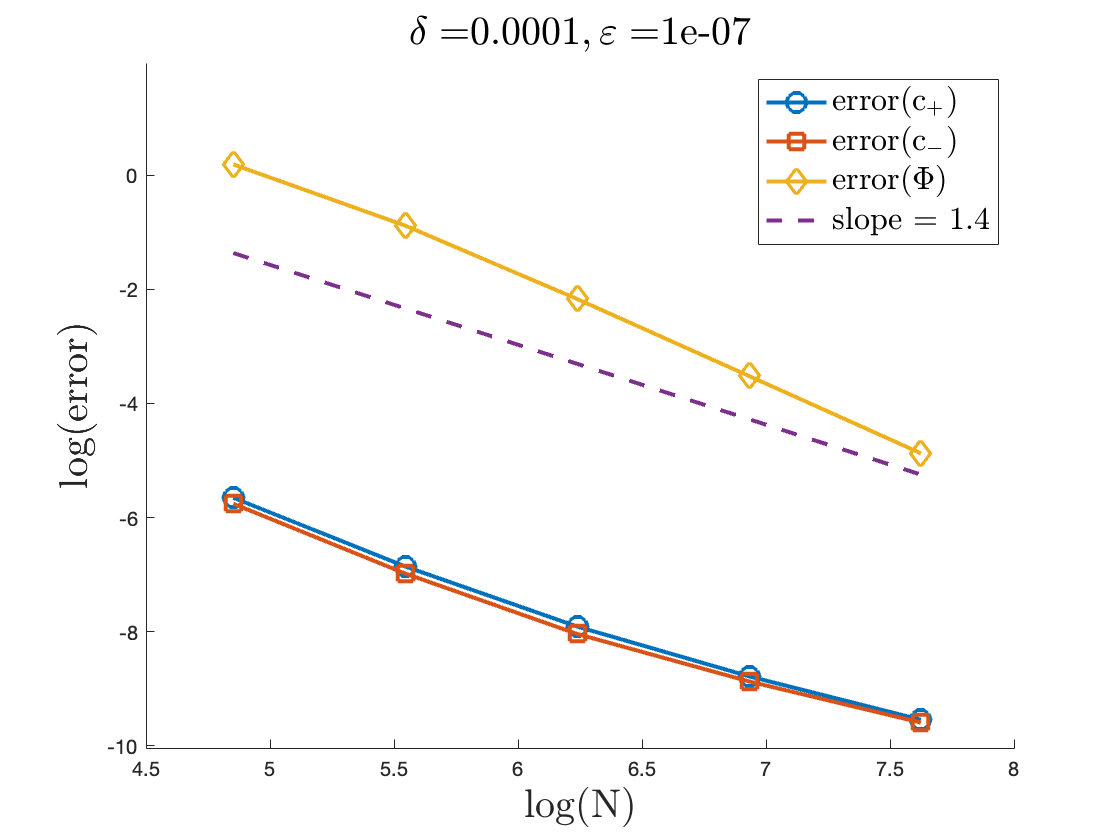}
\put(1,37){(f)}
\end{overpic}
\end{minipage}        
    \caption{\textit{Space and time accuracy orders of the $\delta$-model at final time $t = 0.1$, for different values of $\delta$ and $\varepsilon$. Simulations details are provided in Section~\ref{section_accuracy_test_1D}.}}
    \label{fig:accuracy_1D}
\end{figure}

\subsection{Time discretization in one dimension}
\label{sec:time_discretization_1D}
In this subsection, we describe the time discretization for the one dimensional problem. 
If we denote with $c_\pm^\delta$ and $c_\pm^0$ the solutions of the full {or $\delta$-model} and {MPNP (or 0-model)} models, respectively, the problems~\eqref{eq_semi_discrete_full} and \eqref{eq_semi_discrete_mul} can be summarized as follows
\begin{equation}\label{eq:Theta_eps}
	\td{ \textbf{q}^\delta }{t} = \Theta^\delta [ \textbf{q}^\delta ]\textbf{q}^\delta,
\end{equation}
where $\textbf{q}^\delta = [c_{+,h},c_{-,h},\Phi_h]^\top$, and we distinguish between the cases $\delta > 0$ and $\delta = 0$:
\begin{equation}
\label{eq:sys_eps_1D}
  \Theta[{\bf q^\delta}] = \begin{pmatrix}
\mathbb L^{\rm 1D} + \mathbb D^{\rm 1D}_+[U_+] & {\underbar 0} & \mathbb D^{\rm 1D}_+[c_{+,h}]\\
\underbar 0 & \mathbb L^{\rm 1D} + \mathbb D^{\rm 1D}_-[U_-]  & \mathbb D^{\rm 1D}_-[c_{-,h}]\\
-\mathbb I  & \mathbb I & - \varepsilon \mathbb L^{\rm 1D} \\
\end{pmatrix}, \quad {\rm if } \, \delta > 0
\end{equation}
\begin{equation}
\label{eq:sys_0_1D}
  \Theta[{\bf q^0}] = \begin{pmatrix}
\mathbb L^{\rm 1D} & {\underbar 0} & \mathbb D^{\rm 1D}_+[c_{+,h}]\\
\underbar 0 & \mathbb L^{\rm 1D} & \mathbb D^{\rm 1D}_-[c_{-,h}]\\
-\mathbb I  & \mathbb I & - \varepsilon \mathbb L^{\rm 1D} \\
\end{pmatrix}, \quad {\rm if } \, \delta = 0
\end{equation}
with appropriate boundary conditions, defined in Eqs.~(\ref{eq:c_ghost}-\ref{eq:Phi_ghost}).


We apply IMEX method to \eqref{eq:Theta_eps}. Let us first set $\textbf{q}^1_E = \textbf{q}^n$, then the stage fluxes are calculated as
\begin{subequations}  
\label{eq_imex1}
\begin{align}
\label{eq_imex_qE}
    \textbf{q}_E^{i}& = \textbf{q}^n + \Delta t\,\sum_{j=1}^{\rm i-1}\widetilde a_{i,j}\Theta[\textbf{q}_E^j]\,\textbf{q}_I^j, \quad i = 1,\cdots,s \\
\label{eq_imex_qI}
  \textbf{q}_I^{i} &= \textbf{Q}^n + \Delta t\,\sum_{j=1}^{\rm i} a_{i,j}\Theta[\textbf{q}_E^j]\,\textbf{q}_I^j, \quad i = 1,\cdots,s 
\end{align}
\end{subequations}
and the numerical solution is finally updated with
\begin{align}
\label{eq_imex_q_bi}
   \textbf{q}^{n+1} &= \textbf{q}^n + \Delta t\sum_{i=1}^{\rm s} b(i) \Theta[\textbf{q}_E^i]\,\textbf{q}_I^i
\end{align}
where $\rm s$ is the number of stages  of the scheme, and $\Delta t>0$ the time step. We choose a two stages IMEX-RK methods \cite{ascher1997implicit,pareschi2003high,pareschi2005implicit}, with the double Butcher tableau of the form
\begin{equation}
\label{b_tableau}
    \begin{array}{c|cc}
        0 & 0 & 0  \\
        1/(2\gamma) & 1/(2\gamma) & 0 \\ \hline
         & 1-\gamma & \gamma 
    \end{array} \hspace{4cm}
        \begin{array}{c|cc}
        \gamma & \gamma & 0  \\
        1 & 1-\gamma & \gamma \\ \hline
         & 1-\gamma & \gamma 
    \end{array}
\end{equation}
This scheme is L-stable and stiffly accurate. We refer to this scheme as IMEX-SA(2,2,2).


\section{Numerical results in one dimension}
\label{section_accuracy_test_1D}
In this section, we test the accuracy of the numerical method. We choose
an exact solution $ {\bf q}^{\rm exa}$ and augment the system~\eqref{eq:Theta_eps} as:
\begin{equation}\label{eq:Theta_eps2}
	\td{ \textbf{q}^\delta }{t} = \Theta^\delta [ \textbf{q}^\delta ] + {\bf f}^{\rm exa}({\bf q}^{\rm exa}),
\end{equation}
choosing ${\bf f}^{\rm exa}({\bf q}^{\rm exa}) = [f_+(c_+^{\rm exa}),f_-(c_-^{\rm exa}),f(\Phi^{\rm exa})]$ in such a way that $ {\bf q} = {\bf q}^{\rm exa}$ is the exact solution, see for example \cite{roache2002code}. We choose the following exact solutions:
\begin{subequations}
\label{exact_rho}
\begin{align} 	
c_\pm^{\rm exa}(x,t) &= v_0 \left( \cos(t)^2 c_\pm^0(x) + \sin(t)^2 c_\pm^1(x) \right), \\ 
\Phi^{\rm exa}(x,t) & = \cos(t) \cos(2\pi x), \\ \nonumber
c_\pm^\beta(x) &= \frac{1}{\sqrt{2\pi}\sigma}\exp\left(-{\left(x-x_\pm^\beta\right)^2}/{\sigma}\right),\quad \beta = 0,1
\end{align}
\end{subequations}
where $x_+^0 = 0.45, x_-^0 = 0.5, x_+^1 = 0.5, x_-^1 = 0.55$, $M = 3$ (as we did in \cite{astuto2023multiscale}) and $v_0 = 10^{-4}$ {denotes, in the one-dimensional setting, the volume (length in 1D and surface in 2D) occupied by the solute.}  

We compute the ${L}^2$ norm of the relative error at $t=0.1$, as follows
\begin{align} 	\label{relativeerr}
{\rm error}({\bf q}) = \frac{|| {\bf q} - {\bf q}^{\rm exa}||}{||{\bf q}^{\rm exa}||},
\end{align}
for different values of $N$ and show the results in Fig.~\ref{fig:accuracy_1D}, for $\delta \neq 0$, and in Fig.~\ref{fig:accuracy_1D_MPNP}, for $\delta = 0$. {We compare the space and time accuracy of the method, showing that the observed convergence order is independent of $\delta$, see Fig.~\ref{fig:accuracy_1D} ($\delta = 10^{-2}$ in panels (a),(c),(e) and $\delta = 10^{-4}$ in panels (b),(d),(f)), while it depends on $\varepsilon$ ($\varepsilon = 10^{-1}$ in panels (a),(b), $\varepsilon = 10^{-4}$ in panels (c),(d), and $\varepsilon = 10^{-7}$ in panels (e),(f)). In particular, as $\varepsilon \to 0$, the order degrades from second order to approximately first order. In Section~\ref{sec:QNL}, we introduce a new formulation of the problem that is able to maintain the second order accuracy for smaller values of $\varepsilon$.
}

\begin{figure}
    \centering
\begin{minipage}[b]
		{.32\textwidth}
		\centering
	\begin{overpic}[abs,width=\textwidth,unit=1mm,scale=.25]{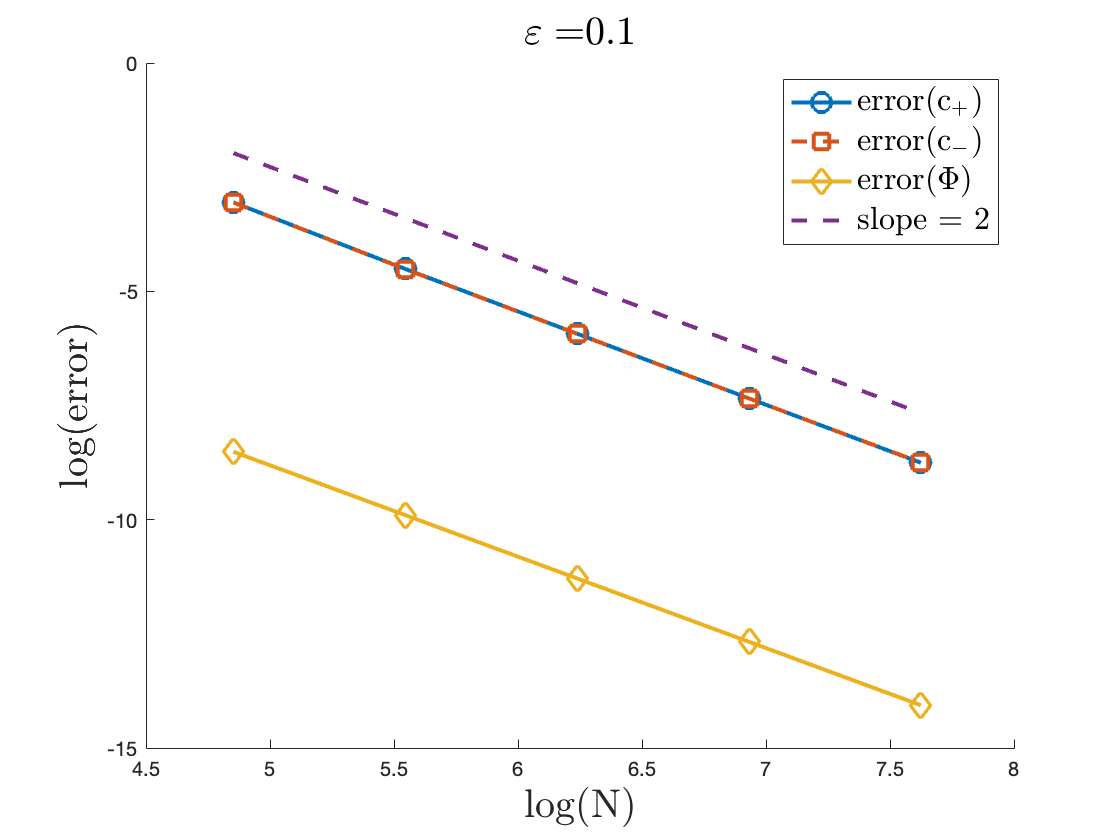}
\put(1,37){(a)}
\end{overpic}
\end{minipage}         
\begin{minipage}[b]
		{.32\textwidth}
		\centering
	\begin{overpic}[abs,width=\textwidth,unit=1mm,scale=.25]{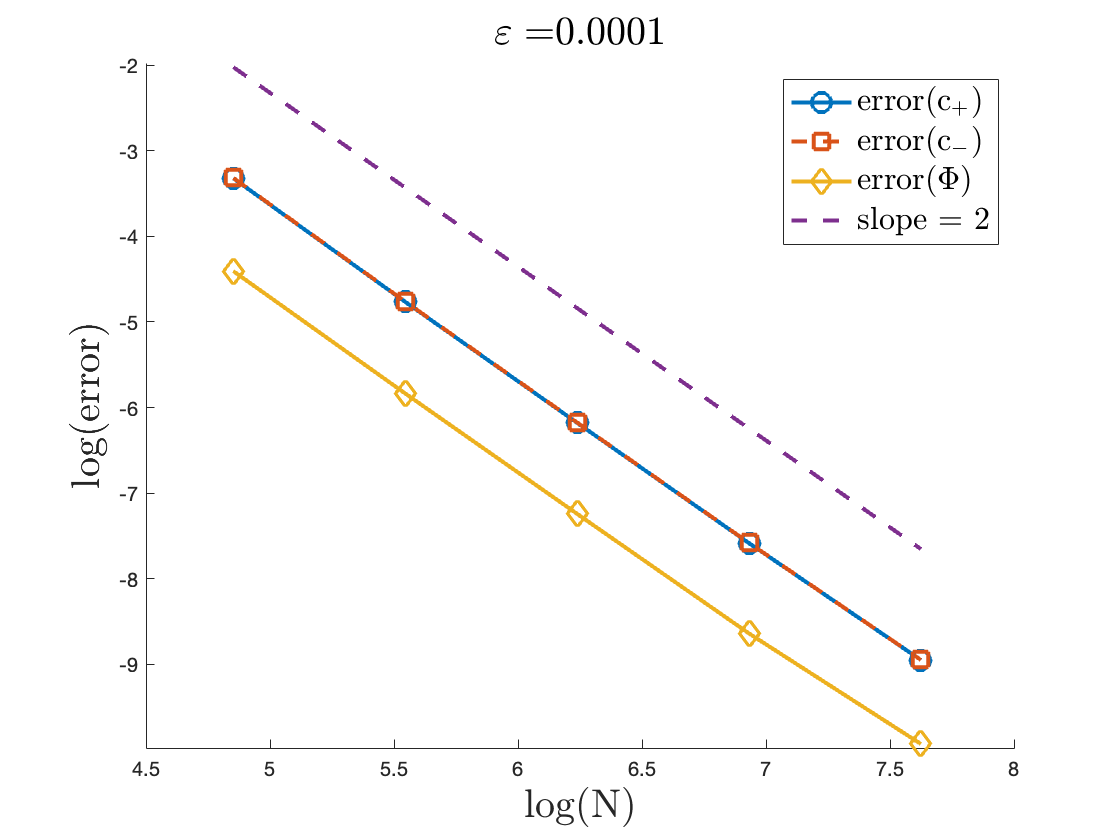}
\put(1,37){(b)}
\end{overpic}
\end{minipage}        
\begin{minipage}[b]
		{.32\textwidth}
		\centering
	\begin{overpic}[abs,width=\textwidth,unit=1mm,scale=.25]{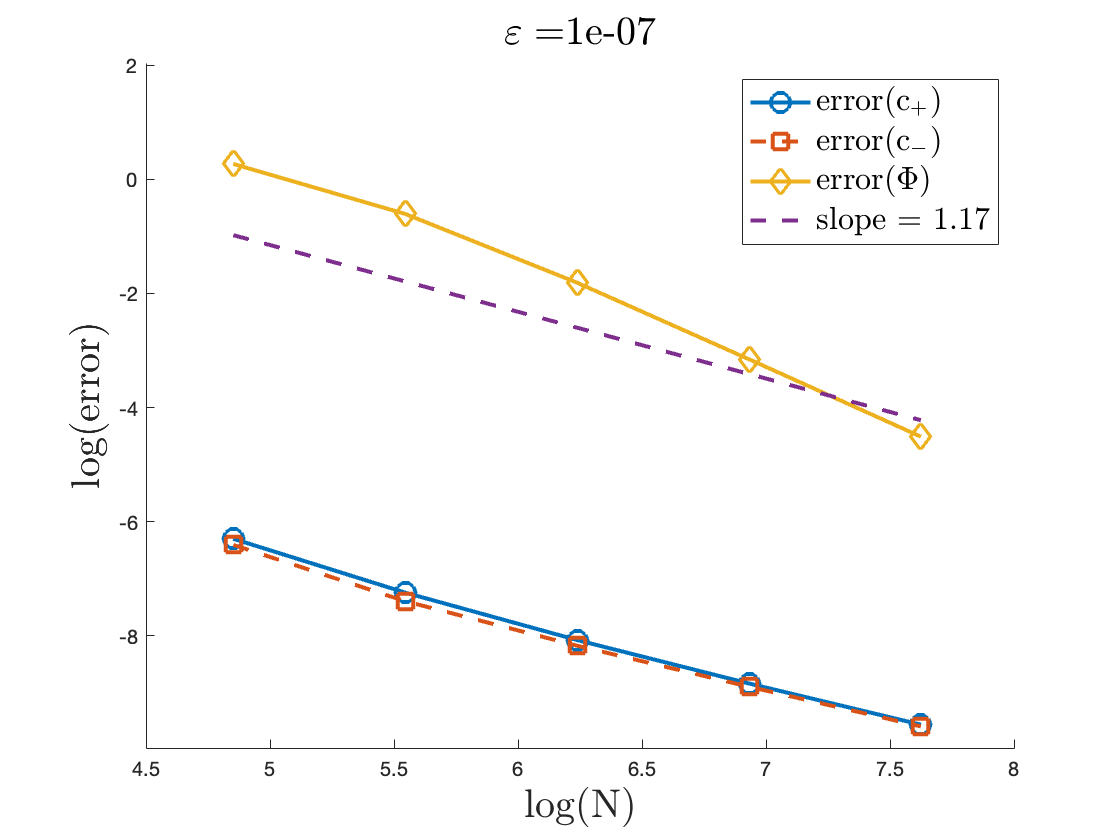}
\put(1,37){(c)}
\end{overpic}
\end{minipage}         
\caption{\textit{Space and time accuracy orders of the 0-model at final time $t = 0.1$, for different values of $\varepsilon$. Simulations details are provided in Section~\ref{section_accuracy_test_1D}.}}
\label{fig:accuracy_1D_MPNP}
\end{figure}

A qualitative comparison, of the full and multiscale models, is shown in Fig.~\ref{fig:comp_qualitative}, for different values of $\delta$. {As $\delta \to 0$, the solutions of the $\delta$-model and of the $0$-model overlap in the fluid domain $\Omega_{\rm f}^\delta$, as expected. In order to perform a quantitative comparison between the two models, in Fig.~\ref{fig:diff_cm_delta}, we calculate the difference between $c_-^\delta$ and $c_-^0$ in the bubble domain $\Omega^\delta_{\rm b}$ which we denote by ${\rm diff_m}$, defined as}

\begin{equation}
\label{eq_diffm}
\displaystyle {\rm diff_m} = \frac{\left|\int_{-\delta}^{\delta L}c_-^\delta(x,t)d x - M c_-^0(x=0,t) \right|}{v_0} 
\end{equation}

{Last but not least, in order to show the efficiency of the MPNP model, we highlight its advantages. Although the computational time of the MPNP model is comparable to that of the full system (see Table~\ref{comp_time}) for fixed parameters $\delta$ and $\Delta x$, the main benefit lies in its ability to provide a stable and accurate solution in the regime of realistic $\delta$ using a reasonable number of points. In contrast, the full model requires a much finer discretization to avoid numerical oscillations when considering realistic $\delta$ to describe the bubble surface ($\delta \approx 10^{-9}$, and reasonably $\Delta x \approx \delta/100$). Fig.~\ref{fig_comp_time} shows the numerical oscillations that appear near the origin for $\delta = 0.0125$ and $N_x = 500$ points in the $\delta$-model, whereas no oscillations are observed for the multiscale model with the same number of points.}

\begin{figure}[h]
    \centering
    \includegraphics[width = 0.4\textwidth]{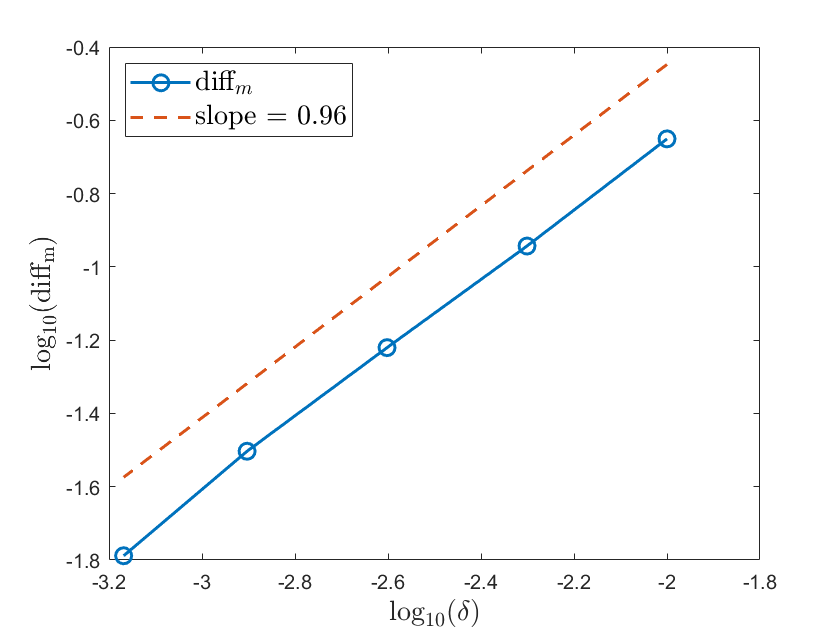}
    \caption{\textit{Difference between the solutions $c_-^\delta$ and $c_-^0$ of the $\delta$-model and 0-model, respectively. The quantity is computed is Eq.~\eqref{eq_diffm}.}}
\label{fig:diff_cm_delta}
\end{figure}

\begin{figure}
    \centering
\begin{minipage}[b]
		{.45\textwidth}
		\centering
	\begin{overpic}[abs,width=\textwidth,unit=1mm,scale=.25]{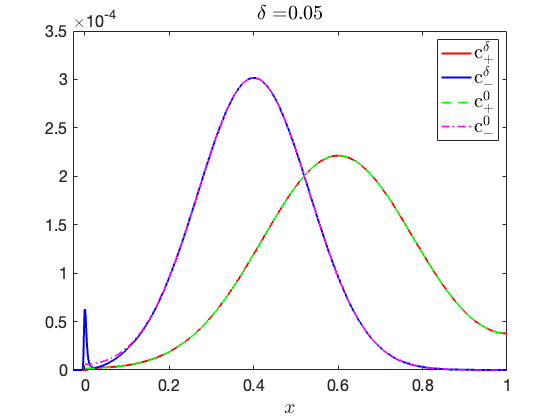}
\put(1,52){(a)}
\end{overpic}
\end{minipage}         
\begin{minipage}[b]
		{.45\textwidth}
		\centering
	\begin{overpic}[abs,width=\textwidth,unit=1mm,scale=.25]{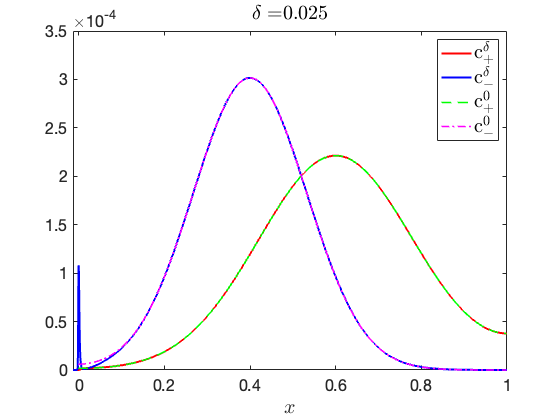}
\put(1,52){(b)}
\end{overpic}
\end{minipage}        
\begin{minipage}[b]
		{.45\textwidth}
		\centering
	\begin{overpic}[abs,width=\textwidth,unit=1mm,scale=.25]{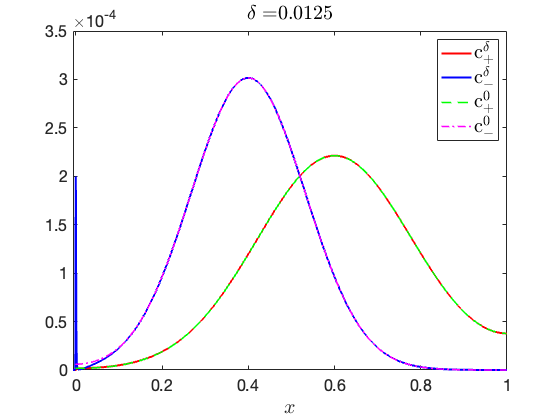}
\put(1,52){(c)}
\end{overpic}
\end{minipage}         
\begin{minipage}[b]
		{.45\textwidth}
		\centering
	\begin{overpic}[abs,width=\textwidth,unit=1mm,scale=.25]{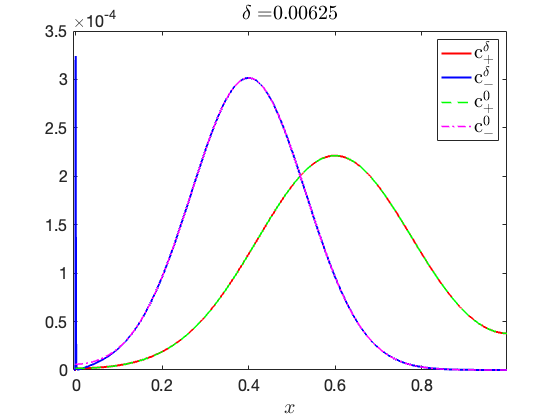}
\put(1,52){(d)}
\end{overpic}
\end{minipage}        
\begin{minipage}[b]
		{.45\textwidth}
		\centering
	\begin{overpic}[abs,width=\textwidth,unit=1mm,scale=.25]{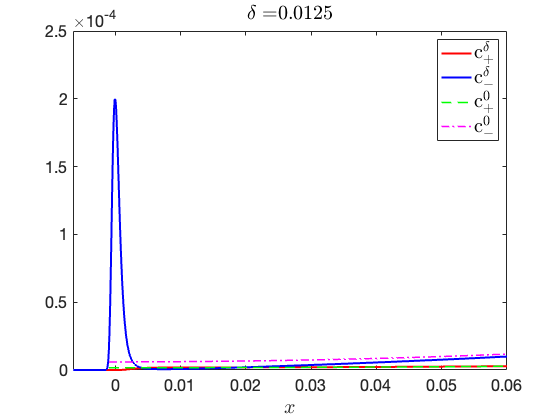}
\put(1,52){(e)}
\end{overpic}
\end{minipage}         
\begin{minipage}[b]
		{.45\textwidth}
		\centering
	\begin{overpic}[abs,width=\textwidth,unit=1mm,scale=.25]{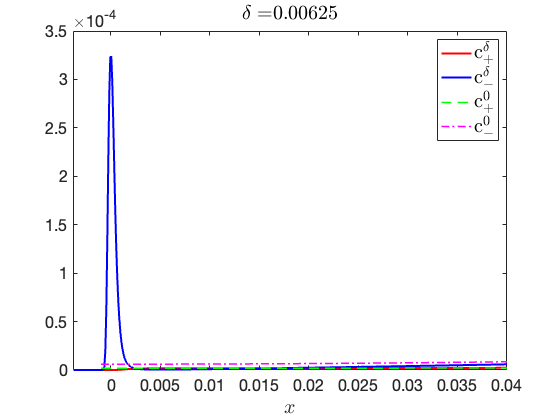}
\put(1,52){(f)}
\end{overpic}
\end{minipage}        
    \caption{\textit{Comparison between the solutions $c_\pm^\delta$ and $c_\pm^0$ of the $\delta$-model and the 0-model, respectively, for different values of $\delta$, {number of points $N_x = 10000$  and final time $t = 1.5$}. Zoom-in {of panels (b) and (c) in (e) and (d)}.}}
\label{fig:comp_qualitative}
\end{figure}

\begin{figure}
    \centering
\begin{minipage}[b]
		{.45\textwidth}
		\centering
	\begin{overpic}[abs,width=\textwidth,unit=1mm,scale=.25]{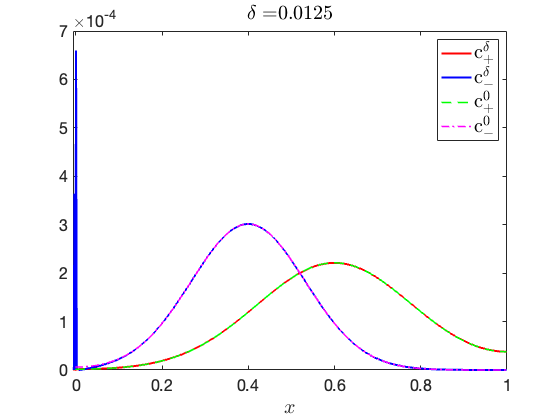}
\put(1,52){(a)}
\end{overpic}
\end{minipage}         
\begin{minipage}[b]
		{.45\textwidth}
		\centering
	\begin{overpic}[abs,width=\textwidth,unit=1mm,scale=.25]{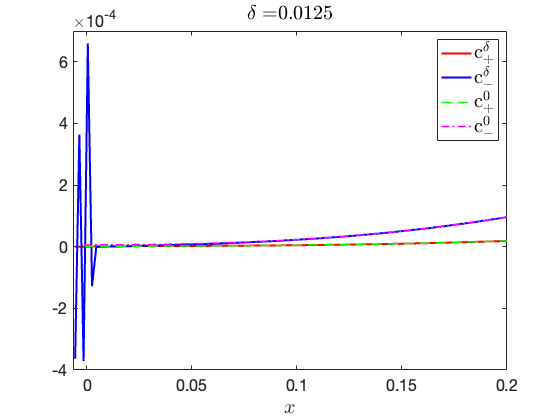}
\put(1,52){(b)}
\end{overpic}
\end{minipage}               
    \caption{\textit{{Comparison between the solutions $c_\pm^\delta$ and $c_\pm^0$ of the $\delta$-model and the 0-model, respectively, for $N_x=500$ number of points and final time $t = 1.5$. Zoom-in in panel (b).}}}
\label{fig_comp_time}
\end{figure}

\begin{table}
	\centering  
	\begin{tabular}{|c|c|c|}
		\hline
		 & $\delta$-model & 0-model \\ 
		\hline\hline 
        $ t_{\rm comp}$ & 0.0076582$s$ & 0.0070704$s$ \\
		\hline 
	\end{tabular}
	\caption{\textit{{Computational time for one time iteration and $N_x = 500$.}}}
	\label{comp_time}
\end{table}

\section{Quasi--Neutral Limit (QNL)}
\label{sec:QNL}
In this section, we deduce the Quasi-Neutral limit for the MPNP system in (\ref{system_multiscale_all}-\ref{system_multiscale_bc}). 

{Considering system~\eqref{system_multiscale_all}, we notice that the term $\varepsilon$ is very small and that is source of stiffness for the system. For example, in Figs.~\ref{fig:accuracy_1D}-\ref{fig:accuracy_1D_MPNP}, panels (c), we see that the order of the numerical scheme deteriorates to first order when considering $\varepsilon = 10^{-7}$. This justifies the introduction of Asymptotic-Preserving numerical schemes in Section~\ref{sec:time_discretization_2D} for small, but not negligible, values of $\varepsilon$. Another strategy, commonly used when interested in the limit case $\varepsilon \to 0$, is the approach based on the so called \textit{Quasi-Neutral Limit} (QNL), proposed in \cite{jungel}. In that regime, 
both species diffuse at the same rate with a common diffusivity that is intermediate between the ones of the two species.}

{To apply the quasi-neutral limit we need to perform the limit $\varepsilon\rightarrow 0$ in system~\eqref{system_multiscale_all}.}

We start introducing two new quantities, that we obtain from the sum and the difference of the two concentrations, as follows
\[ \mathcal C = { {c_+} + {c_-}}, \qquad \mathcal Q = \frac{{c_+ - c_-}}\varepsilon\]
where $\mathcal{Q}$  is proportional to the difference of the net charge density divided by  $\varepsilon$, and presumably remains finite in the Quasi-Neutral limit $\varepsilon\to 0$. In this way, we avoid the instability in the Poisson equation for the electrostatic potential caused by the strong degeneracy when $\varepsilon \approx 0$. 

We rewrite system~(\ref{system_multiscale_all}-\ref{system_multiscale_bc}), using the two new quantities $\mathcal{C}$ and $\mathcal{Q}$, obtaining
\begin{subequations}
\label{eq_CQP_system}
\begin{align}
\label{eq_CQP_system_C}
    \frac{\partial \mathcal C}{\partial t} & = \widetilde D \Delta \mathcal C + \varepsilon \widehat D \Delta \mathcal Q + \nabla \cdot \left(\left(\widehat D \mathcal C + \varepsilon\widetilde D  \mathcal Q \right)\nabla \Phi \right), \qquad {\rm in } \, \Omega \\ \label{eq_CQP_system_Q}
  \frac{\partial \mathcal Q}{\partial t} & =  \frac{\widehat D}{\varepsilon} \Delta \mathcal C + \widetilde D \Delta \mathcal Q+ \nabla \cdot \left(\left( \frac{\widetilde D}{\varepsilon} \mathcal C + \widehat D \mathcal Q \right)\nabla \Phi \right), \qquad {\rm in } \, \Omega \\ \label{eq_CQP_system_Phi}
   -\Delta \Phi & = \mathcal Q, \qquad {\rm in } \, \Omega
   \end{align}
\end{subequations}
with boundary conditions
\begin{subequations}
\label{eq_CQP_system_bc}
 \begin{align}
  \label{eq_CQP_system_bc_C}
 \frac{M}{2} \frac{\partial \mathcal C}{\partial t} - \varepsilon \frac{M}{2} \frac{\partial  \mathcal Q}{\partial t} & = \widetilde D \frac{\partial \mathcal C}{\partial n} + \varepsilon \widehat D \frac{\partial \mathcal Q}{\partial n}  + \left(\widehat D \mathcal C+ \varepsilon \widetilde D \mathcal Q \right) \frac{\partial \Phi}{\partial n}, \qquad {\rm on } \, \Gamma_\mathcal{B} \\
      -\frac{M}{2\varepsilon} \frac{\partial \mathcal C}{\partial t} +  \frac{M}{2} \frac{\partial  \mathcal Q}{\partial t} & = \frac{\widehat D}{\varepsilon} \frac{\partial \mathcal C}{\partial n} + \widetilde D \frac{\partial \mathcal Q}{\partial n}  + \left(\frac{\widetilde D}{\varepsilon} \mathcal C+  \widehat D \mathcal Q \right) \frac{\partial \Phi}{\partial n}, \qquad {\rm on } \, \Gamma_\mathcal{B} \\ 
      \frac{\partial \Phi}{\partial n} & = \frac{M}{2\varepsilon}\mathcal C - \frac{M}{2}\mathcal Q, \qquad {\rm on } \, \Gamma_\mathcal{B}\\ 
   0 & = \widetilde D \frac{\partial \mathcal C}{\partial n} + \varepsilon \widehat D \frac{\partial \mathcal Q}{\partial n} , \qquad {\rm on } \, \Gamma_S \\  
    0 & = \frac{\widehat D}{\varepsilon} \frac{\partial \mathcal C}{\partial n} +  \widetilde D \frac{\partial \mathcal Q}{\partial n}, \qquad {\rm on } \, \Gamma_S \\
\frac{\partial \Phi}{\partial n} & = 0 , \qquad {\rm on } \, \Gamma_S  
\end{align}
\end{subequations}
where $\widetilde D = \left({D_+ + D_-}\right)/{2}, \widehat D = \left({D_+ - D_-}\right)/{2}.$ 

As $\varepsilon \to 0,$ we obtain the Quasi-Neutral limit model
\begin{subequations}
\label{eq:QNL}
\begin{align}
\label{eq_CQP_system_lim}
    \frac{\partial \mathcal C^0}{\partial t} & = \frac{\widetilde D^{\, 2} - \widehat D^{\, 2}}{\widetilde D} \Delta \mathcal C^0 \qquad {\rm in } \, \Omega \\
  0 & = \widehat D  \Delta \mathcal C^0 + \widetilde D \nabla \cdot \left( \mathcal C^0 \nabla \Phi^0 \right) \qquad {\rm in } \, \Omega \\
   -\Delta \Phi^0 & = \mathcal Q^0 \qquad {\rm in } \, \Omega 
\end{align}
\end{subequations}
with boundary conditions
\begin{subequations}
\label{eq:QNL_bc}
\begin{align}
   \frac{M}{2} \frac{\partial \mathcal C^{0}}{\partial t}  & = \widetilde D \frac{\partial \mathcal C^{0}}{\partial n}  + \widehat D \mathcal C\frac{\partial \Phi^{0}}{\partial n} \qquad {\rm on } \, \Gamma_\mathcal{B} \\
      -\frac{M}{2} \frac{\partial \mathcal C^{0}}{\partial t}  & = \widehat D \frac{\partial \mathcal C^{0}}{\partial n}  + \widetilde D \mathcal C \frac{\partial \Phi^{0}}{\partial n} \qquad {\rm on } \, \Gamma_\mathcal{B} \\ 
    \frac{\partial \mathcal C^{0}}{\partial n} & = 0  \qquad {\rm on } \, \Gamma_S \\  
\frac{\partial \Phi^{0}}{\partial n} & = 0 \qquad {\rm on } \, \Gamma_S    \label{eq_CQP_system_lim_bc}
\end{align}
\end{subequations}
\begin{figure}
    \centering
\begin{minipage}[b]
		{.32\textwidth}
		\centering
	\begin{overpic}[abs,width=\textwidth,unit=1mm,scale=.25]{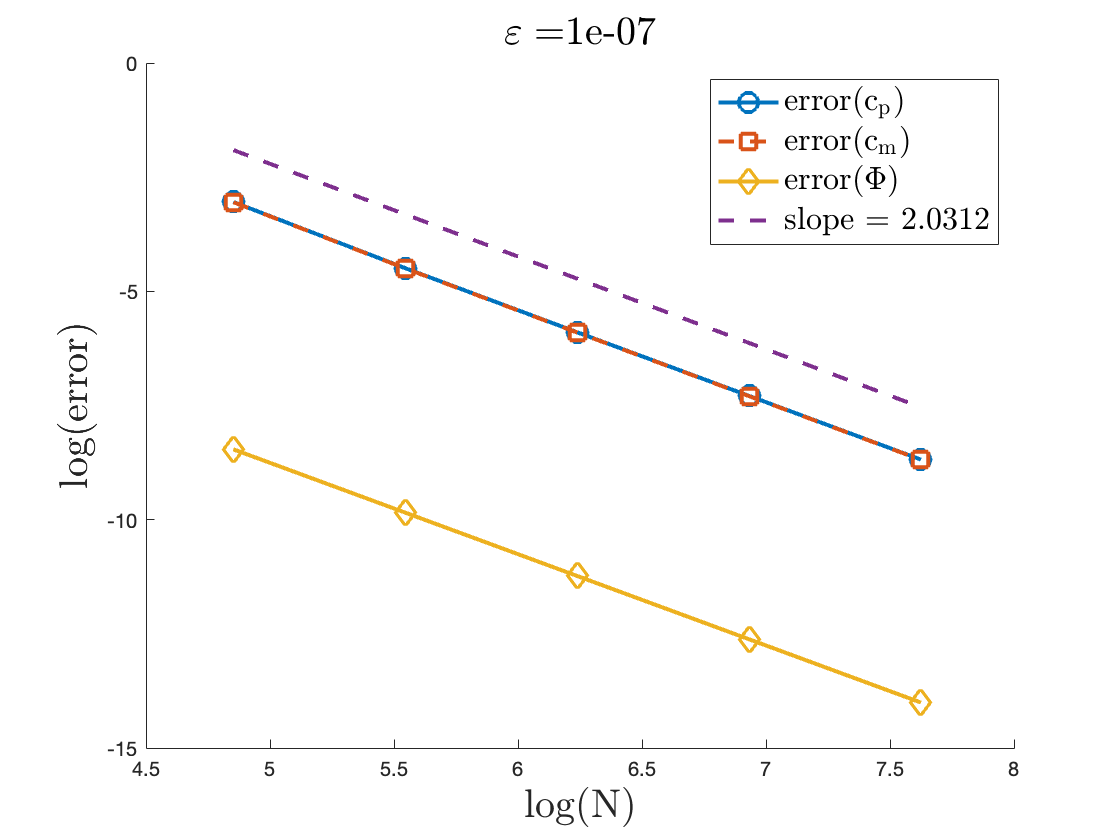}
\put(1,37){(a)}
\end{overpic}
\end{minipage}         
\begin{minipage}[b]
		{.32\textwidth}
		\centering
	\begin{overpic}[abs,width=\textwidth,unit=1mm,scale=.25]{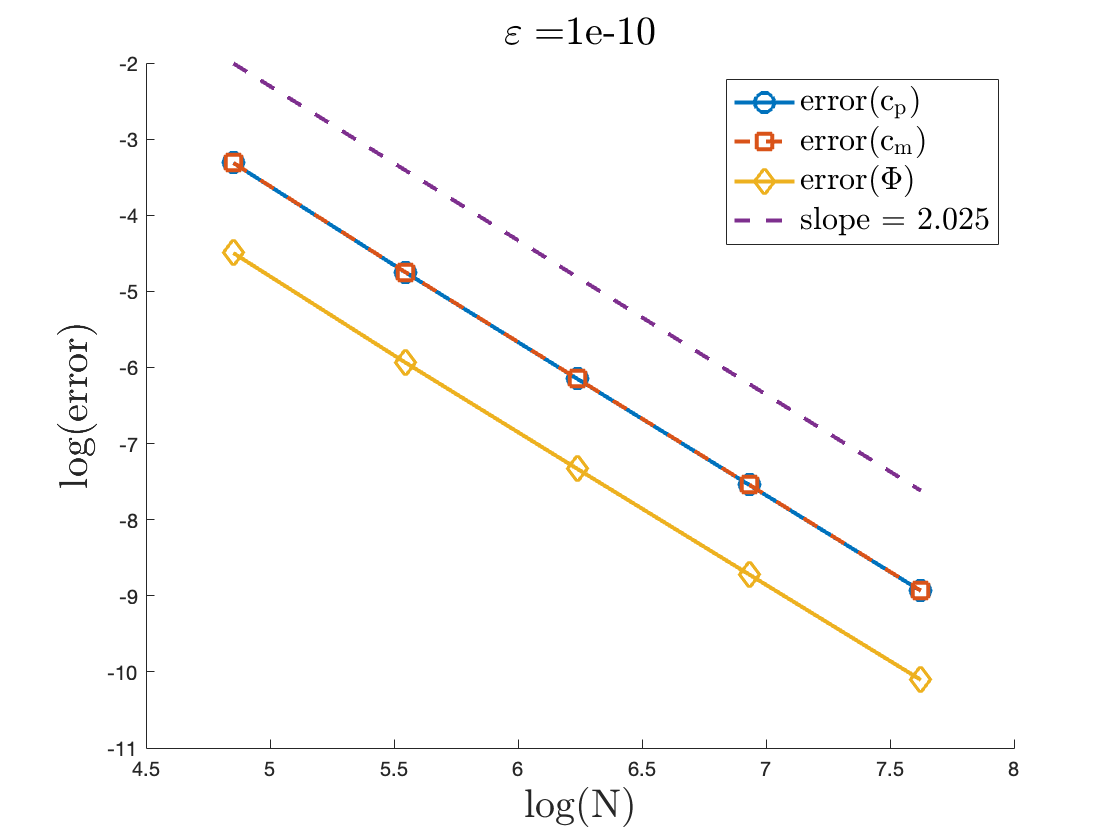}
\put(1,37){(b)}
\end{overpic}
\end{minipage}        
\begin{minipage}[b]
		{.32\textwidth}
		\centering
	\begin{overpic}[abs,width=\textwidth,unit=1mm,scale=.25]{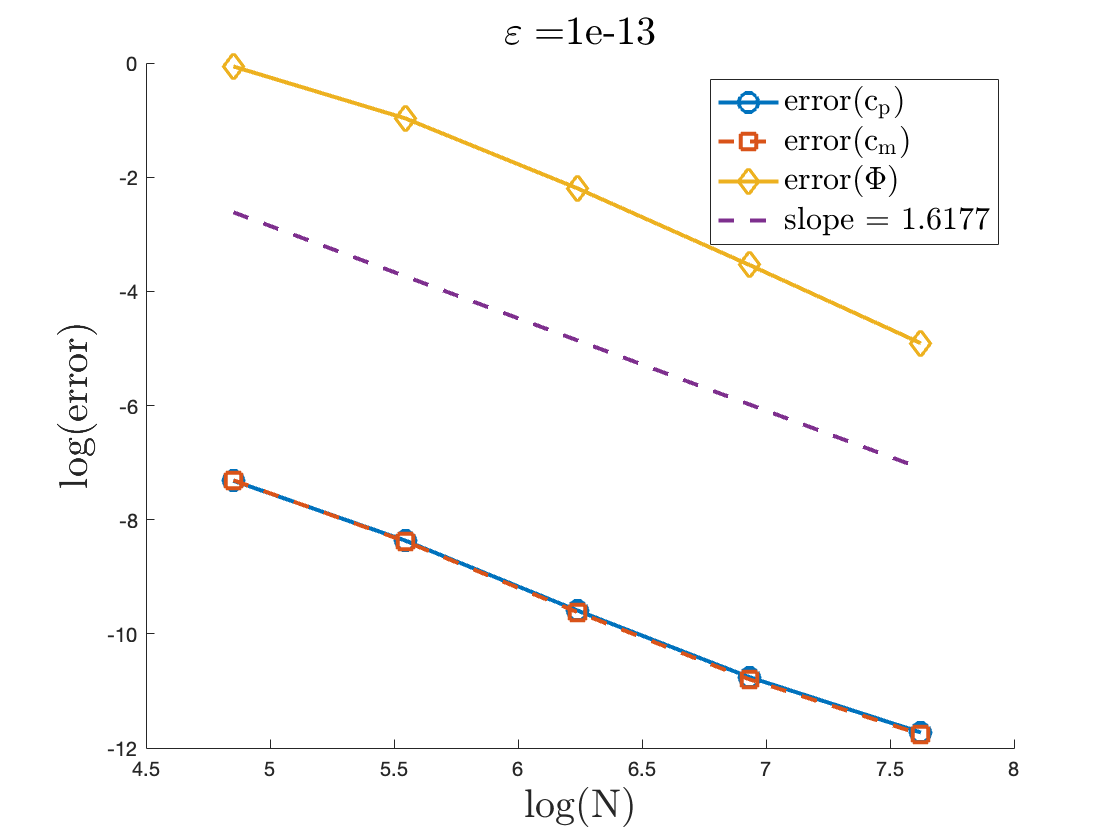}
\put(1,37){(c)}
\end{overpic}
\end{minipage}         
    \caption{\textit{Space and time accuracy orders of the QNL system (see Eqs.~(\ref{eq_CQP_system}-\ref{eq_CQP_system_bc})) at final time $t = 0.1$, for different values of $\varepsilon$. Simulations details are provided in Section~\ref{section_accuracy_test_1D}. In these tests $\Delta t = 0.1\Delta x$.}}
    \label{fig:accuracy_1D_MPNP_CQ}
\end{figure}

\subsection*{Numerical tests in one dimension}
To complete the one-dimensional tests, in Fig.~\ref{fig:accuracy_1D_MPNP_CQ} we show the accuracy in space and time, as we do in Section~\ref{section_accuracy_test_1D}. In practice, we choose the same exact solutions defined in Eq.~\eqref{exact_rho}, but this time our initial conditions are $\displaystyle \mathcal C^{\rm exa} = {c_+^{\rm exa}} + {c_-^{\rm exa}}$ and $\displaystyle \mathcal Q^{\rm exa} = \frac{c_+^{\rm exa} - c_-^{\rm exa}}{\varepsilon}$.

{The formulation proposed in Eqs.~(\ref{eq_CQP_system}-\ref{eq_CQP_system_bc}) allows us to consider smaller values of $\varepsilon$. While in Figs.~\ref{fig:accuracy_1D}-\ref{fig:accuracy_1D_MPNP}, panels (c), we observe a degradation of the accuracy order for $\varepsilon \leq 10^{-7}$, Fig.~\ref{fig:accuracy_1D_MPNP_CQ} demonstrates that second order accuracy is preserved even for smaller values, such as $\varepsilon = 10^{-10}$.}

\section{Two-dimensional space discretization and Asymptotic Preserving numerical scheme as $\protect \varepsilon \to 0$}
\label{sec_2D}
In this section, we describe the space and time discretization for the model~(\ref{eq_CQP_system}-\ref{eq_CQP_system_bc}). For the space discretization, we follow the strategy in \cite{astuto2024nodal}, a recently developed ghost-FEM method. The same space discretization has been extended to the numerical solution of biological network formation in a leaf-shaped domain in  \cite{astuto2024self}. Regarding the time discretization, we consider a second order Asymptotic Preserving numerical scheme
\cite[Section 1.3]{boscarino2024implicit}, \cite{jin1999efficient}.

\subsection{Variational formulation}
Here we consider the variational formulation of the system (\ref{eq_CQP_system}-\ref{eq_CQP_system_bc}). We introduce the two spaces $V$ and $W$ defined as
\begin{equation}
\label{eq:FEMspaces}
    V = \biggl\{ v\in H^1(\Omega) \biggr\}, \qquad W = \biggl\{ w\in H^1(\Omega) : \int_\Omega w\,{ \rm d}\Omega = 0\biggr\}.
\end{equation}
When imposing homogeneous Neumann boundary conditions in the Poisson equation, the potential
$\Phi$ is determined up to an additive constant, leading to a lack of uniqueness. To address this issue, in Eq.~\eqref{eq:FEMspaces} we impose an additional constraint in the space $W$, namely, that $\Phi$ has zero mean. 

Multiplying \eqref{eq_CQP_system_C} by a test function $v \in V$, and integrating over $\Omega$, we obtain
\begin{subequations}
\begin{align}
&\int_\Omega \frac{\partial \mathcal{C}}{\partial t} v\,{ \rm d}\Omega = 
\widetilde D \int_\Gamma \left(\nabla \mathcal C \cdot n \right) v \, { \rm d}\Gamma 
- \widetilde D \int_\Omega \nabla \mathcal{C} \cdot \nabla v \, { \rm d}\Omega 
+ \varepsilon \widehat D \int_\Gamma \left(\nabla \mathcal Q \cdot n \right)v  \, { \rm d}\Gamma 
- \varepsilon \widehat D \int_\Omega \nabla \mathcal Q \cdot \nabla v \, { \rm d}\Omega \\ & + \int_\Gamma \left( \left(  \widehat D \mathcal C  + \varepsilon \widetilde D \mathcal Q  \right)  \nabla \Phi \cdot n \right)  v \, { \rm d}\Gamma 
- \int_\Omega \left( \widehat D \mathcal C + \varepsilon \widetilde D \mathcal Q \right) \nabla \Phi \cdot \nabla v \, { \rm d}\Omega
\end{align}
\end{subequations}
Taking into account the boundary condition in Eq.~\eqref{eq_CQP_system_bc_C}, we have
\begin{subequations}
\begin{align} \nonumber
&\int_\Omega \frac{\partial \mathcal{C}}{\partial t} v\,{ \rm d}\Omega = 
- \widetilde D \int_\Omega \nabla \mathcal{C} \cdot \nabla v \, { \rm d}\Omega  
- \varepsilon \widehat D \int_\Omega \nabla \mathcal Q \cdot \nabla v \, { \rm d}\Omega  
- \int_\Omega \left( \widehat D \mathcal C + \varepsilon \widetilde D \mathcal Q \right) \nabla \Phi \cdot \nabla v \, { \rm d}\Omega  \\ & 
+ \int_{\Gamma_\mathcal{B}} \left( \frac{M}{2} \frac{\partial \mathcal C}{\partial t} - \varepsilon \frac{M}{2} \frac{\partial \mathcal Q}{\partial t}\right)  v \, { \rm d}\Gamma
\end{align}
\end{subequations} 
We adopt the same procedure on Eqs.~(\ref{eq_CQP_system_Q}-\ref{eq_CQP_system_Phi}), thus obtaining the variational formulation of our problem. 

\begin{pro}\label{pro:variational}
Find $\mathcal C(t), \mathcal Q(t) \in V$ and $\Phi(t) \in W$ for almost every $t\in(0,T)$, such that
\begin{subequations}
\label{eq_var_formulation}
\begin{align} 
&\left( \frac{\partial \mathcal{C}}{\partial t}, v\right)  = - \widetilde D \left( \nabla \mathcal{C},\nabla v \right)   - \varepsilon \widehat D \left( \nabla \mathcal Q, \nabla v \right)  - \left(\left( \widehat D \mathcal C + \varepsilon \widetilde D \mathcal Q \right) \nabla \Phi,\nabla v \right)   + \left( \frac{M}{2} \frac{\partial \mathcal C}{\partial t} - \varepsilon \frac{M}{2} \frac{\partial \mathcal Q}{\partial t}, v \right)_{L^2(\Gamma_\mathcal{B})}  \\
&\left( \frac{\partial \mathcal{Q}}{\partial t}, q \right)  = - \frac{\widehat D}{\varepsilon} \left( \nabla \mathcal{C},\nabla q  \right)  -  \widetilde D \left( \nabla \mathcal Q, \nabla q  \right)  - \left(\left( \frac{\widetilde D}{\varepsilon} \mathcal C +  \widehat D \mathcal Q \right) \nabla \Phi,\nabla q  \right)  + \left( -\frac{M}{2\varepsilon} \frac{\partial \mathcal C}{\partial t} + \frac{M}{2} \frac{\partial \mathcal Q}{\partial t}, q \right)_{L^2(\Gamma_\mathcal{B})} \\
&\left( \nabla \Phi,\nabla w\right)  + \left( \frac{M}{2\varepsilon}\mathcal C - \frac{M}{2}\mathcal Q, w\right)_{L^2(\Gamma_\mathcal{B})} = \left(\mathcal{Q},w\right) 
\end{align}
\end{subequations}
\end{pro}
where we denoted by $(\cdot,\cdot)$ the scalar product in $L^2(\Omega)$.
{The problem formulation for the Quasi-Neutral limit is obtained, analogously, from Eqs.~\eqref{eq:QNL}.
\begin{pro}\label{pro:variational0}
Find $\mathcal C^0(t), \mathcal Q^0(t) \in V$ and $\Phi^0(t) \in W$ for almost every $t\in(0,T)$, such that
\begin{subequations}
\label{eq_var_formulation2}
\begin{align} 
\label{eq_var_formulation2a}
\left( \frac{\partial \mathcal{C}^0}{\partial t}, v\right)  -  \frac{M}{2} \left(  \frac{\partial \mathcal C^0}{\partial t}, v \right)_{L^2(\Gamma_\mathcal{B})} & = - \widetilde D \left( \nabla \mathcal{C}^0,\nabla v \right)    -  \widehat D\left( \mathcal C^0  \nabla \Phi^0,\nabla v \right)   \\ \label{eq_var_formulation2b}
\frac{M}{2} \left(  \frac{\partial \mathcal C^0}{\partial t} , q \right)_{L^2(\Gamma_\mathcal{B})} & = -{\widehat D} \left( \nabla \mathcal{C}^0,\nabla q  \right)   - \widetilde D\left(\mathcal C^0 \nabla \Phi^0,\nabla q  \right)   \\
\left( \nabla \Phi^0,\nabla w\right) &= \left(\mathcal{Q}^0,w\right).
\end{align}
\end{subequations}
\end{pro}}

\subsection{Space discretization in two dimensions}
\label{sec:space_discretization_2D}
In this section, we adopt a two-dimensional space discretization based on finite elements method \cite{astuto2024nodal,astuto2024comparison}. For the sake of completeness, we provide the relevant details of the spatial discretization.

The domain is $\Omega = [0,1]^2\setminus \mathcal{B} \subset R$, with $\mathcal{B}$ a circle centered in $(x_c,y_c)$ and radius $R_{\mathcal{B}}$, and $R$ a rectangular region. The set of grid points will be denoted by $\mathcal N$, with $\# \mathcal N = (1+N)^2$, the active nodes (i.e.,\ internal $\mathcal{I}$ or ghost $\mathcal{G}$) by $\mathcal A = \mathcal{I}\cup\mathcal{G} \subset \mathcal N$, the set of inactive points by $\mathcal O \subset \mathcal N$, with $\mathcal O\cup\mathcal A = \mathcal N$ and $\mathcal O \cap \mathcal A = \emptyset$ and the set of cells by $\mathcal C$, with $\# \mathcal C = N^2$. Finally, we denote by $\Omega_c = R\setminus \Omega$ the outer region in $R$.

Following the approach shown in  \cite{Osher2002,Russo2000,book:72748, Sussman1994}, the domain $\Omega$ is implicitly defined by a level set function $\phi(x,y)$ that is negative inside $\Omega$, positive in $R\setminus \Omega$ and zero on the boundary ${\Gamma_\mathcal{B}}$:
\begin{equation}
	\Omega = \{(x,y): \phi(x,y) < 0\}, \qquad
	\Gamma_\mathcal{B} = \{(x,y): \phi(x,y) = 0\}.
\end{equation}
\begin{figure}[H]
    \centering
    \begin{minipage}{.49\textwidth}
\centering\begin{overpic}[abs,width=0.85\textwidth,unit=1mm,scale=.25]{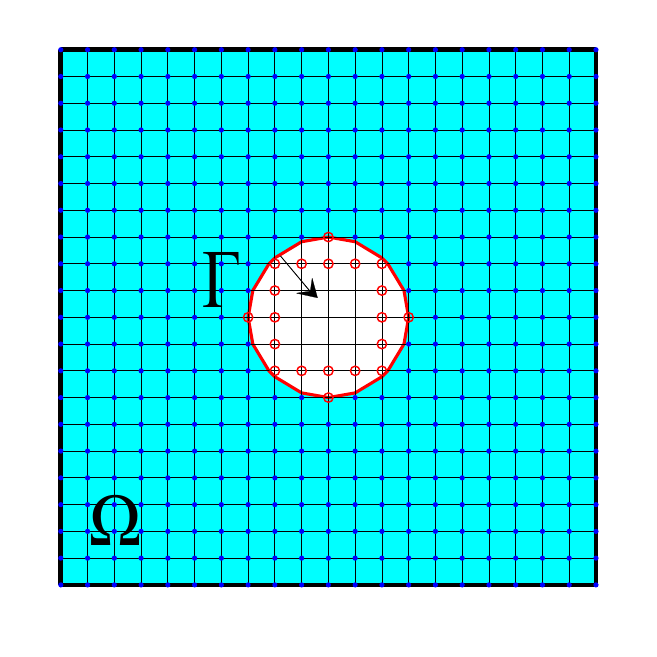}        
\put(-2,54){(a)}
\put(33,34){$\widehat n$}
\put(22.5,35){$\mathcal B$}
\end{overpic}
    \end{minipage}
    \begin{minipage}{.49\textwidth}
\centering\begin{overpic}[abs,width=0.85\textwidth,unit=1mm,scale=.25]{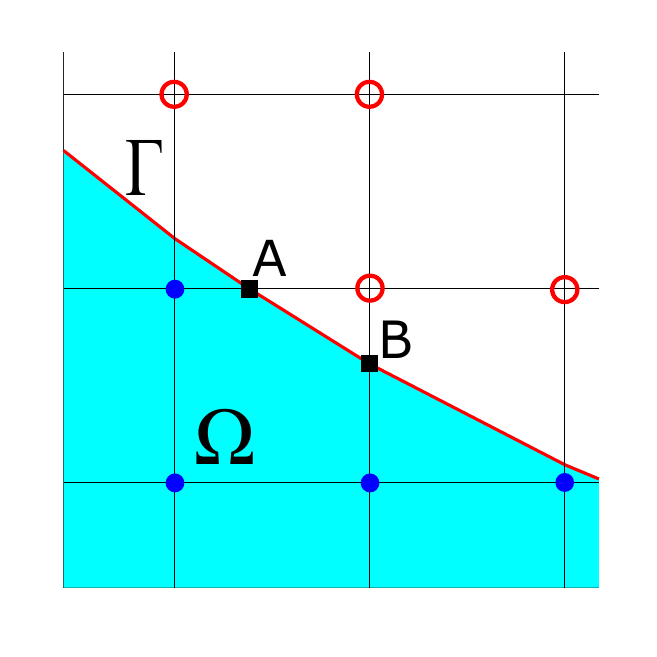}        
\put(-2,54){(b)}
\put(15,46){$\mathcal B$}
\end{overpic}
    \end{minipage}
\caption{\textit{Discretization of the computational domain. $\Omega$ is the {light blue} region inside the unit square $R$. (a): classification of the grid points: the blue points are the internal ones while the red circles denote the ghost points. (b): points of intersection between the grid and the boundary $\Gamma$ (see the definition of $A$ and $B$ in Algorithm \ref{alg_ab}).}}  
\label{fig:Domain2D}
\end{figure}

Here, we define the set of ghost points $\mathcal{G}$, which are grid points that belong to $\Omega_c$, { and are vertices of cut cells}, formally defined as
\begin{equation}
\notag
	(x,y) \in \mathcal{G} \iff (x,y) \in {\mathcal N}\cap \Omega_c  \text{ and } \{(x \pm h,y),(x,y\pm h), (x \pm h,y\pm h) \} \cap \mathcal I \neq \emptyset.
\end{equation}


The discrete spaces $V_h$ and $W_h$ are given by the piecewise bilinear functions which are continuous in $R$.
As a basis of $V_h$ and $W_h$, we choose the following functions:
\begin{equation}
    v_{i}(x,y) = \max\left\{
        \left(1-\frac{|x-x_i|}{h}\right),0
    \right\}\max\left\{\left(1-\frac{|y-y_i|}{h}\right),0
    \right\},
    \label{eq:V_h2}
\end{equation}
with $i = (i_1,i_2)$ an index that identifies a node on the grid. The generic element $u_h\in V_h$ (or, $w_h\in W_h$) {has } the following representation
\begin{equation}
    u_h(x,y) = \sum_{i\in\nodes}u_i v_i(x,y).
    \label{eq:u_h2}
\end{equation}

To solve the variational problem \eqref{eq_var_formulation}, we employ a finite-dimensional discretization. Specifically, the functions $\mathcal C, \mathcal Q \in V$ and $\Phi \in W$ are approximated by functions $\mathcal C_h, \mathcal Q_h \in V_h$ and $\Phi_h \in W_h$, respectively.  
To perform computations, the domain $\Omega$ is approximated by a polygonal domain $\Omega_h$. This approximation also extends to the boundary $\partial \Omega = \Gamma \cup \Gamma_\mathcal{B}$,  represented by $\Gamma_h$ and $\Gamma_{\mathcal{B},h}$, respectively. Consequently, the original integrals, defined over $\Omega$ and its boundaries $\Gamma$ and $\Gamma_\mathcal{B}$ are now evaluated over $\Omega_h$, $\Gamma_h$ and $\Gamma_{\mathcal{B},h}$, respectively.
\begin{pro}\label{pro:variational_2}
Find $\mathcal C_h, \mathcal Q_h \in V_h$ and $\Phi_h \in W_h$  such that, for almost every $t\in(0,T)$, it holds
\begin{subequations}
\label{eq_varform}
\begin{align}
&\left( \frac{\partial \mathcal{C}_h}{\partial t}, v_h \right)_{L^2\left(\Omega_h\right)} = - \widetilde D \left( \nabla \mathcal{C}_h,\nabla v_h  \right)_{L^2\left(\Omega_h\right)}  - \varepsilon \widehat D \left( \nabla \mathcal Q_h, \nabla v_h  \right)_{L^2\left(\Omega_h\right)} - \left(\left( \widehat D \mathcal C_h + \varepsilon \widetilde D \mathcal Q_h \right) \nabla \Phi_h,\nabla v_h  \right)_{L^2\left(\Omega_h\right)}  \\ & + \left( \frac{M}{2} \frac{\partial \mathcal C_h}{\partial t} - \varepsilon \frac{M}{2} \frac{\partial \mathcal Q_h}{\partial t}, v_h  \right)_{L^2(\Gamma_{\mathcal{B},h})} \\ 
&\varepsilon\left( \frac{\partial \mathcal{Q}_h}{\partial t}, q_h \right)_{L^2\left(\Omega_h\right)} = - {\widehat D} \left( \nabla \mathcal{C}_h,\nabla q_h  \right)_{L^2\left(\Omega_h\right)}  -  \varepsilon \widetilde D \left( \nabla \mathcal Q_h, \nabla q_h  \right)_{L^2\left(\Omega_h\right)} - \left(\left( {\widetilde D} \mathcal C_h + {\varepsilon} \widehat D \mathcal Q_h \right) \nabla \Phi_h,\nabla q_h  \right)_{L^2\left(\Omega_h\right)}  \\ & + \left( -\frac{M}{2} \frac{\partial \mathcal C_h}{\partial t} + \varepsilon\frac{M}{2} \frac{\partial \mathcal Q_h}{\partial t}, q_h  \right)_{L^2(\Gamma_{\mathcal{B},h})} \\
&\left( \nabla \Phi_h,\nabla w_h\right)_{L^2\left(\Omega_h\right)} + \left( \frac{M}{2\varepsilon}\mathcal C_h - \frac{M}{2}\mathcal Q_h, w_h\right)_{L^2\left(\Gamma_h^B\right)} = \left(\mathcal{Q}_h,w_h\right)_{L^2\left(\Omega_h\right)},
\end{align}
\end{subequations}
with appropriate initial conditions that we define in Section~\ref{section_numerics_2D}.
\end{pro}

\begin{figure}[H]
    \centering
    \begin{minipage}{.4\textwidth}
\begin{overpic}[abs,width=\textwidth,unit=1mm,scale=.25]{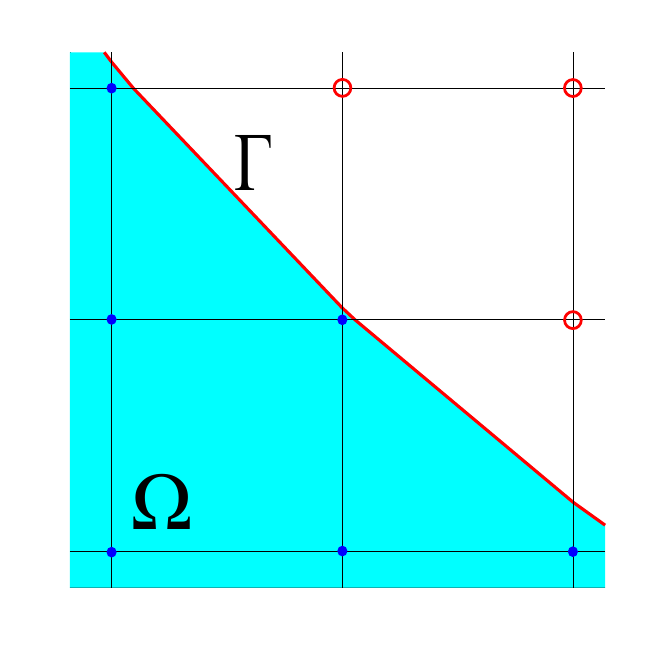}  
\put(-2,50){(a)}
\put(25.5,44.5){$\mathcal B$}
\end{overpic}
    \end{minipage}
    \begin{minipage}{.4\textwidth}
\begin{overpic}[abs,width=\textwidth,unit=1mm,scale=.25]{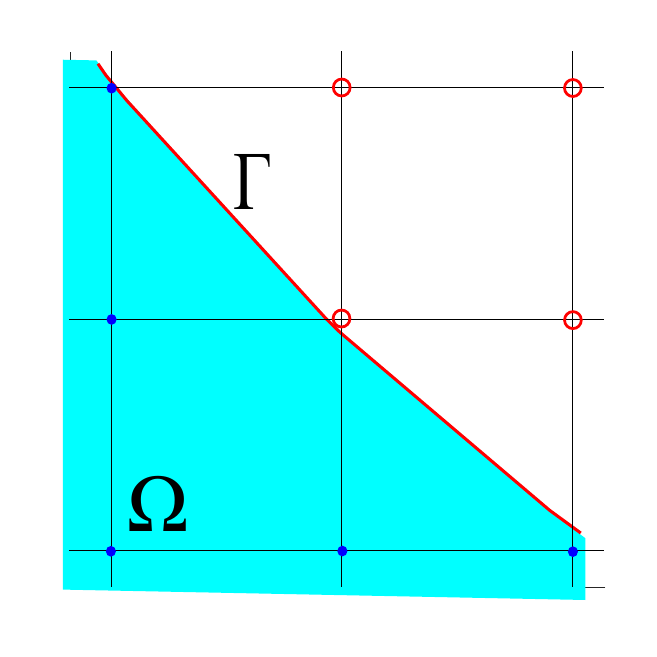}        
\put(-1,50){(b)}
\put(25,43){$\mathcal B$}
\end{overpic}
    \end{minipage}
\caption{\textit{Grid before and after snapping technique. (a): representation of the cell related to the internal point $P$ (blue points), whose distance from $\Gamma$ is less than $h^2$; (b): zoom-in of the shape of the domain, after the grid point $P$ has changed its classification, from internal to ghost point (red circles).
}}  
\label{fig:snapping}
\end{figure}

\begin{algorithm}
\caption{Computation of the intersection of the boundary with the grid (see Fig. \ref{fig:ref_cell})}\label{alg_ab}
\begin{algorithmic}
\State $k_4 = k_0$
\For{i = 0:3}
\If{$\phi(k_i)\phi(k_{i+1})<0$} 
   \State $\theta = \phi(k_i)/(\phi(k_i)-\phi(k_{i+1})) $ 
   \State $P = \theta k_{i+1}+(1-\theta)k_i$
   \If{$\phi(k_i)<0 $}
       \State ${\bf A}:=P$
       \Else
       \State ${\bf B}:=P$
   \EndIf
\EndIf 
\EndFor
\end{algorithmic}
\end{algorithm}

To compute the integrals shown in Problem~\ref{pro:variational_2}, we use exact quadrature rules. To explain our strategy, let us start considering the product between two test functions $v_i,v_j \in V_h$ restricted within the cell $K \in \mathcal C$, i.e., 
\begin{align}
\label{eq_int_discr2}
      (v_i,v_j)_\Omega = \sum_{K\in\mathcal C} \left( \left.v_i\right|_K , \left.v_j\right|_K \right),
\end{align}
where, for each $K$, $\left.\varphi_i\right|_K$ is the restriction of $\varphi_i$ in cell $K$, which is a bilinear function and takes value 1 in grid node $i$ and 0 in the other vertices of the cell $K$; see \cite{astuto2024nodal} for more details. 

We observe that the product of two elements in $Q_h$ is an element of $\mathbb{Q}_2(K)$, i.e., the set of bi-quadratic polynomials in $K$.
\begin{figure}[H]
    \centering
\begin{minipage}{.4\textwidth}
\centering
    \end{minipage}
\begin{minipage}{.49\textwidth}
\centering \begin{overpic}[abs,width=0.65\textwidth,unit=1mm,scale=.25]{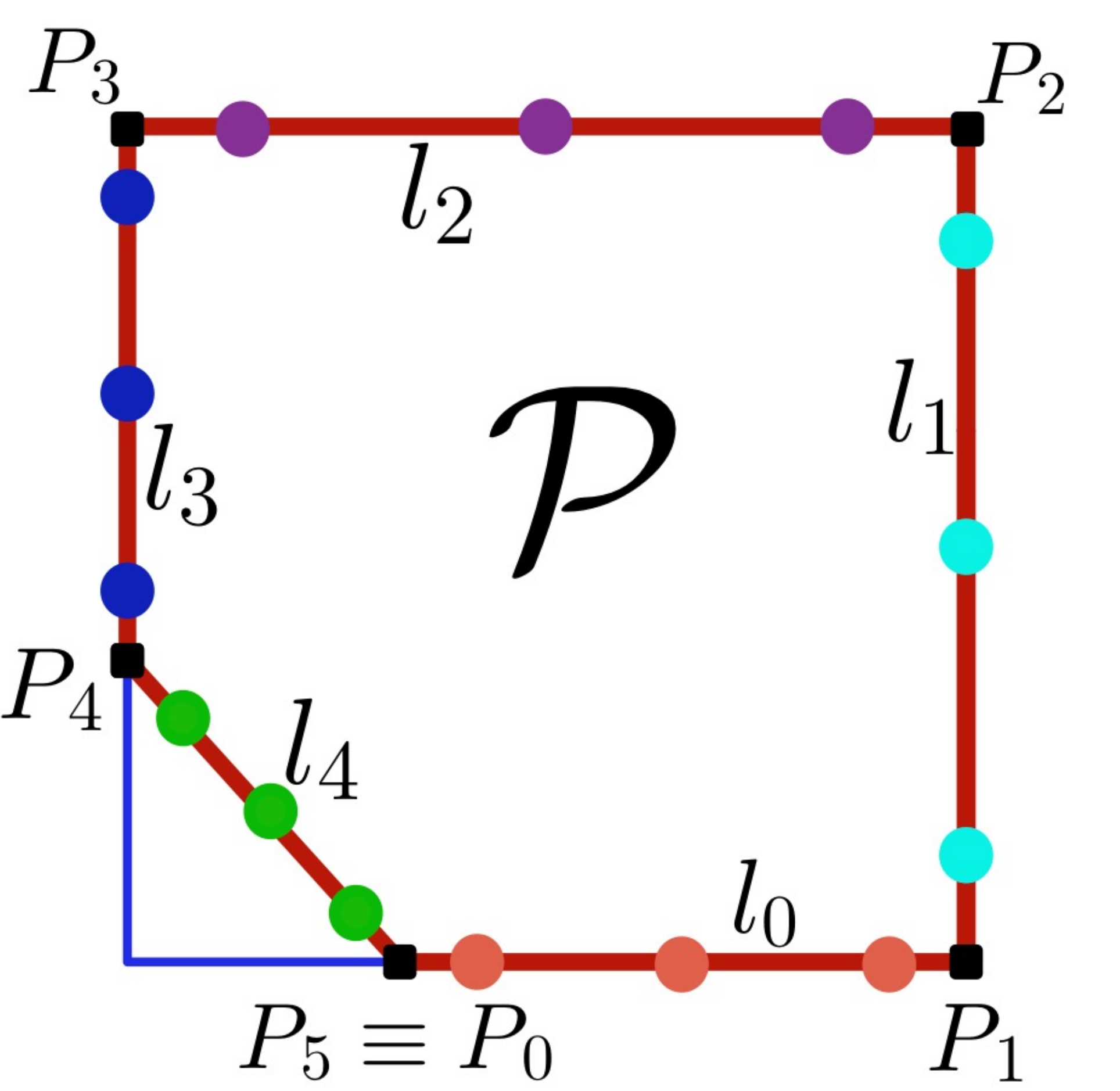}  
\end{overpic}
    \end{minipage}
\caption{\textit{Scheme of the three quadrature points (circles) for each edge $l_i, \, i = 0,\cdots,4$. The squared points represent the vertices $P_i,\, i = 0,\cdots,4,$ of the polygon $\mathcal P$.}}  
\label{fig:ref_cell}
\end{figure}

Let us consider a general integrable function $f(x,y)$ defined in $\Omega$, with $F(x,y) = \int f(x,y) dx$ (in our strategy we consider the primitive in $x$ direction; analogue results can be obtained integrating in $y$ direction). We now define a vector function $\textbf{F}= (F,0)^\top$ which has the property \(\nabla \cdot \textbf{F} = f\).
Thus, we have
\begin{align}
\label{eq_diverg}
    \int_K fdx\, dy = \int_K \nabla \cdot  \textbf{F} \,dx\, dy = \int_{\partial K} \textbf{F}\cdot \textbf{n}\, dl,
\end{align}
where we applied Gauss theorem, and $\textbf{n} = (n_x,n_y)$ is the outer normal vector to $K$.
If $\mathcal{P}$ is the generic polygon with $m$ edges $l_r, \, r = 0,\ldots,m-1$ (see Fig.~\ref{fig:ref_cell}) 
 we can express \eqref{eq_diverg} as
\begin{align}
\label{eq_gauss_interval}
   \int_\mathcal{P} \nabla \cdot  \textbf{F}\, dx\, dy = \int_{\partial \mathcal{P}} \textbf{F}\cdot \textbf{n}\, d l = \sum_{r=0}^{m-1} \int_{l_r} \textbf{F}\cdot \textbf{n}\, dl = \sum_{r=0}^{m-1} \int_{l_r}  F {n}_x\, dl = \sum_{r=0}^{m-1} \int_{l_r}  F dy. 
\end{align}
To evaluate the integral over the generic edge $l_r$, we  choose the three-point Gauss-Legendre quadrature rule, which is exact for polynomials in $\mathbb P_5(\mathbb R)$. Thus, we write
\begin{equation}
     \int_{l_r}  F dy = \sum_{s = 1}^3w_sF(\widehat x_{r,s},\widehat y_{r,s}) (y_{P_{r+1}}-y_{P_r})
\end{equation}
where $w_s$ 
and $\widehat x_{r,s}, \, s = 1,2,3$ are the weights and the nodes, respectively, of the considered quadrature rule, see Fig.~\ref{fig:ref_cell} (b).
Choosing $f=v_{k_\eta}v_{k_\mu}\in \mathbb{Q}_2(K)\subset \mathbb P_4(K)$, and making use of \eqref{eq_gauss_interval}, we write
\begin{equation}
     \left( \left.v_i\right|_K , \left.v_j\right|_K \right) = \sum_{r=0}^{m-1} \left(\sum_{s=1}^3 w_s\Phi_{ij}(\widehat x_{r,s},\widehat y_{r,s}) (y_{P_{r+1}}-y_{P_r})\right)|l_r|.
\end{equation}
where $v_{ij} = \int v_i v_j\, dx$, and the formula is exact because $v_{ij} \in \mathbb P_5(K), \, \forall i,j \in \mathcal N$.

When evaluating the integrals described above, we observe stability issues arising from the presence of cut cells near the boundary $\Gamma_\mathcal{B}$. The problem derives from the inability to control the size of the cut cells, which can become arbitrarily small. Consequently, this can lead to a loss of coercivity in the bilinear form. 
In Fig.~\ref{fig:snapping} (a), we see a case in which the stability of the numerical scheme fails, and in panel (b), we present an approach to address this issue. In other words, to avoid instability, 
we evaluate the level set function $\phi$ at the vertices of each cell: if the value is {smaller than a threshold} (that we choose proportional to a power of the length of the cell, i.e.,
{if $0<-\phi<\zeta h^\alpha$, for suitable chosen $\zeta$ and $\alpha$}) we disregard the respective cell {by setting the level set function equal to a small positive value, as illustrated in Algorithm~\ref{alg_snap}; see also \cite{astuto2024nodal,astuto2024comparison}. In our numerical results it corresponds to the machine epsilon.} Alternative techniques are employed to address the ill-conditioning caused by the presence of small cells; see, for instance,
\cite{burman2015cutfem,lehrenfeld2016high,ahlkrona2021cut}.

\begin{algorithm}[H]
\caption{Snapping back to grid}\label{alg_snap}
\begin{algorithmic}
\For{$k \in \mathcal N$} 
\If{$\phi(k)<0 \quad \& \quad |\phi(k)|<
\zeta h^\alpha$} 
   \State $\phi(k) := \texttt{eps} $ 
\EndIf 
\EndFor
\end{algorithmic}
\end{algorithm}

\subsection{Asymptotic Preserving time discretization in two dimensions}
\label{sec:time_discretization_2D}
{The concept of an Asymptotic-Preserving (AP) method has been introduced in \cite{jin1999efficient}. The commutative diagram below shows the main property of AP schemes.  In order to capture the limit using standard numerical methods, one typically needs to solve problems with increasingly small values of $\varepsilon$, while choosing a discretization parameter $h$ fine enough to resolve the corresponding small scales, which may be computationally very expensive. Conversely, an AP scheme becomes a consistent discretization of the limit problem $\mathcal{P}^0$ as $\varepsilon\to 0$, with no order relation between $h$ and $\varepsilon$.} 

\begin{figure}[H]
{\Large \[
  \begin{matrix}
\mathcal P_h^\varepsilon & {\longrightarrow} & \mathcal P^\varepsilon\\
\downarrow &  &  \downarrow\\
\mathcal P^0_h & \longrightarrow & \mathcal P^0 \\
\end{matrix}
\]
}

\caption{The AP diagram. $\mathcal P^\varepsilon$ is the original problem and $\mathcal P^\varepsilon_h$ its numerical approximation characterized by a discretization parameter $h$. The AP property corresponds to the request that $\mathcal P^0_h$ is consistent with $\mathcal{P}^0$ as $\varepsilon \to 0$, independently of $h$.}
\end{figure}
To reformulate system \eqref{eq_varform} using the computational matrices for the spatial derivatives, we express it as follows
\begin{subequations}
\label{eq:FEM_semidiscr}
\begin{align}
& \mathbb B[v_h]\frac{\partial \mathcal{C}_h}{\partial t} - \left( \frac{M}{2} \frac{\partial \mathcal C_h}{\partial t} - \varepsilon \frac{M}{2} \frac{\partial \mathcal Q_h}{\partial t}, v_h  \right)_{L^2\left(\Gamma_{\mathcal{B},h}\right)} = - \widetilde D \, \mathbb L[v_h] \, \mathcal{C}_h  - \varepsilon \widehat D\, \mathbb L[v_h]\, \mathcal Q_h - \mathbb H\left[ \widehat D \mathcal C_h + \varepsilon \widetilde D \mathcal Q_h,v_h \right]  \Phi_h  \\ 
&\varepsilon \mathbb B[q_h] \frac{\partial \mathcal{Q}_h}{\partial t} - \left( -\frac{M}{2} \frac{\partial \mathcal C_h}{\partial t} + \varepsilon\frac{M}{2} \frac{\partial \mathcal Q_h}{\partial t}, q_h  \right)_{L^2\left(\Gamma_{\mathcal{B},h}\right)} = - {\widehat D}\, \mathbb L[q_h]\, \,\mathcal{C}_h  -  \varepsilon \widetilde D \,\mathbb L[q_h]\, \, Q_h - \mathbb H \left[ {\widetilde D}\mathcal C_h + \varepsilon \widehat D \mathcal Q_h,q_h \right]  \Phi_h \\
& \mathbb L[w_h]\, \Phi_h +  \left( \frac{M}{2\varepsilon}\mathcal C_h - \frac{M}{2}\mathcal Q_h, w_h\right)_{L^2\left(\Gamma_{\mathcal{B},h}\right)} = \mathbb B[w_h] \, \mathcal{Q}_h 
\end{align}
\end{subequations}
where $\mathbb L[v_h] $ is the discrete operator that defines the stiffness matrix $\left( \nabla \bullet,\nabla v_h \right)_{L^2\left(\Omega_h\right)}$, such that $ \mathbb L[v_h] \, \mathcal C_h = \left( \nabla \mathbb \mathcal C_h,\nabla v_h\right)_{L^2\left(\Omega_h\right)}$; $ \mathbb H[\mathcal{C}_h,v_h]$ is the operator such that $\mathbb H \left[ \mathcal C_h,v_h \right]  \Phi_h = \left(  \mathcal C_h \nabla \Phi_h,\nabla v_h  \right)_{L^2\left(\Omega_h\right)}$; finally, we define the discrete operator $\mathbb B[v_h]$ for the mass matrix, such that $\mathbb B[v_h] \, \mathcal C_h = \left(\mathcal C_h,v_h\right)_{L^2\left(\Omega_h\right)}$.

Here we perform the limit for  $\varepsilon \to 0$ in \eqref{eq:FEM_semidiscr}, that becomes
\begin{subequations}
\begin{align} \label{eq:bL}
  \mathbb B[v_h]\frac{\partial \mathcal{C}^0_h}{\partial t} - \left( \frac{M}{2} \frac{\partial \mathcal C^0_h}{\partial t}, v_h  \right)_{L^2(\Gamma_{\mathcal{B},h})} & = - \widetilde D \, \mathbb L[v_h] \, \mathcal{C}^0_h  - {\widehat D } \mathbb H\left[ \mathcal C^0_h,v_h \right]  \Phi^0_h  \\ \label{eq:bH}
 \left( \frac{M}{2} \frac{\partial \mathcal C^0_h}{\partial t} , q_h  \right)_{L^2(\Gamma_{\mathcal{B},h})} & = - {\widehat D} \mathbb L[q_h] \,\mathcal{C}^0_h - { \widetilde D} \mathbb H \left[  \mathcal C^0_h,q_h \right]  \Phi^0_h \\
 \mathbb L[w_h] \, \Phi^0_h & = \mathbb B[w_h] \, \mathcal{Q}^0_h.
\end{align}
\end{subequations}
{Rearranging Eq.~\eqref{eq:bH}, we get 
\begin{equation} \mathbb H \left[  \mathcal C^0_h,q_h \right]  \Phi^0_h = - \frac{1}{\widetilde D} \left( \frac{M}{2} \frac{\partial \mathcal C^0_h}{\partial t} , q_h  \right)_{L^2(\Gamma_{\mathcal{B},h})}  - \frac{\widehat D}{\widetilde D}\mathbb L[q_h] \,\mathcal{C}^0_h 
\end{equation}
and we substitute it in Eq.~\eqref{eq:bL}, obtaining}
\begin{subequations}
\label{sys_qnl_discr}
\begin{align}
  \mathbb B[v_h]\frac{\partial \mathcal{C}^0_h}{\partial t} - \left( \frac{M}{2} {\left(1 + \frac{\widehat D}{\widetilde D} \right)}\frac{\partial \mathcal C^0_h}{\partial t}, v_h  \right)_{L^2(\Gamma_{\mathcal{B},h})} & = -\frac{\widetilde D^{\, 2} - \widehat D^{\, 2}}{\widetilde D} \, \mathbb L[v_h] \, \mathcal{C}^0_h  \\ 
\left( \frac{M}{2}  \frac{\partial \mathcal C^0_h}{\partial t} , q_h  \right)_{L^2(\Gamma_{\mathcal{B},h})}  & = -\widehat D\, \mathbb L[q_h] \, \mathcal C^0_h - \widetilde D \mathbb H \left[\mathcal C^0_h,q_h\right]\Phi^0_h \\
\mathbb L[w_h]\, \Phi^0_h & = \mathbb B[w_h]\, \mathcal{Q}^0_h,
\end{align}
\end{subequations}
that is the numerical scheme for {Eqs.~\eqref{eq_var_formulation2}. There, we rearrange Eq.~\eqref{eq_var_formulation2b} to get the expression for $\left(\mathcal C^0 \nabla \Phi^0,\nabla q  \right)$ and substitute it in Eq.~\eqref{eq_var_formulation2a}. The main idea is that applying the same procedure to the variational formulation and to the numerical scheme helps us to obtain the AP property.}


Now, we consider a time discretization for system \eqref{eq:FEM_semidiscr}. Let us define the final time $T$ and the time step as $\Delta t = T/N_{\rm ts},\, N_{\rm ts} \in \mathbb N,$ the number of times steps, denoting the nodes in time by $t^n = n\Delta t$ and $\mathcal{C}_h^n \approx \mathcal{C}_h(t^n), \, n = 0,\cdots,N_{\rm ts}$. A semi--implicit discretization is adopted to achieve second order of accuracy in time. In particular, we make use of implicit-explicit (IMEX) Runge--Kutta schemes \cite{pareschi2000implicit,boscarino2016high,astuto2023self}, which are multi-step methods based on $s$-stages. 

We rewrite Eqs.~\eqref{eq:FEM_semidiscr}  in a vectorial form
\begin{equation}
\label{eq:vectorial}
 \BB\td{\bf Q}{t} =\Theta [{\bf Q}] {\bf Q}
\end{equation}
where 
\begin{equation}
 \BB  = \begin{pmatrix}
 \mathbb B[v_h] & \underbar 0 & \underbar 0 \, \\
\underbar 0 & \varepsilon \mathbb B[q_h] &  \underbar 0 \,\\
\underbar 0 & \underbar 0 & \underbar 0 \, \\
\end{pmatrix} ,\qquad \Theta[{\bf Q}] = \begin{pmatrix}
- \widetilde D \mathbb L[v_h] & - \varepsilon \widehat D \mathbb L[v_h] & - \mathbb H\left[\widehat D \mathcal C_h + \widetilde D \varepsilon \mathcal Q_h,v_h\right]\\
- {\widehat D} \mathbb L[q_h] & - \varepsilon\widetilde D \mathbb L[q_h] &  - \mathbb H\left[{\widetilde D} \mathcal C_h + \widehat D \varepsilon \mathcal Q_h,q_h\right]\\
\underbar 0 & - \mathbb B[w_h] & \mathbb L[w_h] \\
\end{pmatrix} 
\end{equation}
and ${\bf Q} = [\mathcal C_h,\mathcal Q_h,\Phi_h]^\top$. 

To apply the IMEX method to \eqref{eq:vectorial}, we follow the strategy and the Butcher tableau seen in Section~\ref{sec:time_discretization_1D}. Let us first set $\textbf{Q}^1_E = \textbf{Q}^n$, then the stage fluxes are calculated as
\begin{subequations}  
\label{eq_imex}
\begin{align}
\label{eq_imex_QE}
    \BB  \textbf{Q}_E^{i}& =  \BB \textbf{Q}^n + \Delta t\,\sum_{j=1}^{ i-1}\widetilde a_{i,j}\Theta(\textbf{Q}_E^j)\textbf{Q}_I^j, \quad i = 1,\cdots,s \\
\label{eq_imex_QI}
     \BB \textbf{Q}_I^{i} &=  \BB  \textbf{Q}^n + \Delta t\,\sum_{j=1}^{ i}  a_{i,j}\Theta(\textbf{Q}_E^j)\textbf{Q}_I^j, \quad i = 1,\cdots,s 
\end{align}
\end{subequations}
and the numerical solution is finally updated with
\begin{align}
\label{eq_imex_Q_bi}
    \BB  \textbf{Q}^{n+1} &=  \BB  \textbf{Q}^n + \Delta t\sum_{i=1}^{\rm s}  b(i) \Theta(\textbf{Q}_E^i)\textbf{Q}_I^i.
\end{align}
The Butcher tableau that we employ is the one defined in Eq.~\eqref{b_tableau}.

Here we prove that the numerical scheme (\ref{eq_imex}--\ref{eq_imex_Q_bi}) that we design for the system (\ref{eq_CQP_system}-\ref{eq_CQP_system_bc}), it is also a numerical scheme for the limit model in (\ref{eq:QNL}-\ref{eq:QNL_bc}). Let us look at the first steps:
\begin{align*}
    \mathcal{C}_E^1 &= \mathcal C^n \\
    \mathcal{Q}_E^1 &= \mathcal Q^n \\
    \mathbb B[v_h] \, \mathcal{C}_I^1 &= \mathbb B[v_h] \, \mathcal{C}_E^1 - a_{11}\Delta t \left(  \widetilde D \mathbb L[v_h]\, \mathcal C_I^1 +  \varepsilon \widehat D \mathbb L[v_h]\, \mathcal Q_I^1 + \mathbb H \left[\widehat D \mathcal C_E^1 + \widetilde D \varepsilon \mathcal Q_E^1,v_h\right]\Phi_I^1  \right) \\
  \varepsilon \mathbb B[q_h] \, \mathcal{Q}_I^1 &= \varepsilon \mathbb B[q_h] \, \mathcal{Q}_E^1 - a_{11}\Delta t \left(  {\widehat D} \mathbb L[q_h]\, \mathcal C_I^1 + \varepsilon \widetilde D \mathbb L[q_h]\, \mathcal Q_I^1 + \mathbb H \left[\widetilde D  \mathcal C_E^1 + \varepsilon \widehat D \mathcal Q_E^1,q_h\right]\Phi_I^1 \right)  \\
    0 &= \mathbb B[w_h] \, \mathcal Q^1_I - \mathbb L[w_h]\, \Phi_I^1
\end{align*}
and perform the limit $\varepsilon \to 0$. Thus it becomes
\begin{align*}
    \mathcal{C}_E^1 &= \mathcal C^n \\
    \mathcal{Q}_E^1 &= \mathcal Q^n \\
    \mathbb B[v_h] \, \mathcal{C}_I^1 &= \mathbb B[v_h] \, \mathcal{C}_E^1 - a_{11}\Delta t \left( \widetilde D \mathbb L[v_h]\, \mathcal C_I^1 + \widehat D \mathbb H \left[ \mathcal C_E^1,v_h \right]\Phi_I^1 \right) \\
    0 &= \widehat D  \mathbb L[q_h]\, \mathcal C_I^1  + \widetilde D \mathbb H \left[\mathcal C_E^1,q_h \right]\Phi_I^1 \\
    0 &= \mathbb B[w_h] \, \mathcal Q_I^1 - \mathbb L[w_h]\, \Phi_I^1
\end{align*}
that can be rewritten as
\begin{align*}
    \mathcal{C}_E^1 &= \mathcal C^n \\
    \mathcal{Q}_E^1 &= \mathcal Q^n \\
    \mathbb B[v_h] \, \mathcal{C}_I^1 &= \mathbb B[v_h] \, \mathcal{C}_E^1 - a_{11}\Delta t \frac{\widetilde D^{\, 2} - \widehat D^{\, 2}}{\widetilde D} \mathbb L[v_h]\, \mathcal C_I^1  \\
    0 &= \widehat D  \mathbb L[q_h]\, \mathcal C_I^1  + \widetilde D \mathbb H \left[\mathcal C_E^1,q_h \right]\Phi_I^1 \\
    0 &= \mathbb B[w_h] \, \mathcal Q_I^1 - \mathbb L[w_h]\, \Phi_I^1.
\end{align*}

Analogously, we proceed with the second step and obtain 
\begin{align*}
   \mathbb B[v_h] \, \mathcal{C}_E^2 &= \mathbb B[v_h] \, \mathcal{C}_E^1 - \widetilde a_{21}\Delta t \frac{\widetilde D^{\, 2} - \widehat D^{\, 2}}{\widetilde D} \mathbb L[v_h]\, \mathcal C_I^1 \\
    \mathbb B[v_h]\, \mathcal{C}_I^2 &= \mathbb B[v_h] \, \mathcal{C}_E^1 - a_{21}\Delta t \frac{\widetilde D^{\, 2} - \widehat D^{\, 2}}{\widetilde D} \mathbb L[v_h]\, \mathcal C_I^1 - a_{22}\Delta t \frac{\widetilde D^{\, 2} - \widehat D^{\, 2}}{\widetilde D} \mathbb L[v_h]\, \mathcal C_I^2  \\
    0 &= \widehat D  \mathbb L[q_h]\, \mathcal C_I^2  + \widetilde D \mathbb H \left[\mathcal C_E^2,q_h \right]\Phi_I^2 \\
    0 &= \mathbb B[w_h] \, \mathcal Q_I^2 - \mathbb L[w_h]\, \Phi_I^2.
\end{align*}
that is the second order IMEX scheme applied to the limit model in (\ref{eq:QNL}-\ref{eq:QNL_bc}).

\section{Numerical results in two dimensions}
\label{section_numerics_2D}
In this section, we present the results obtained by applying the 
numerical schemes described in Sections~\ref{sec:space_discretization_2D}-\ref{sec:time_discretization_2D}
applied to problem \eqref{pro:variational}
in two space dimensions. 
We define the level-set function $\phi$ as the signed distance from the interface of the bubble $\mathcal{B}$, i.e.,
\[\phi=R_\mathcal{B}-\sqrt{(x-x_c)^2+(y-y_c)^2}
\]
where $(x_c,y_c) = (0.5,0.5)$ is the center of the bubble and $R_\mathcal{B} = 0.05$ its radius. In our numerical simulations, we set $M = 10^{-6}$. This value is calculated using $E = 10$ (see, e.g., \cite{everett1976adsorption}) and $\delta = 10^{-3}$ in Eq.~\eqref{eq:Vminus}.

\begin{figure}[h]
    \centering
\begin{minipage}[b]
		{.49\textwidth}
		\centering
	\begin{overpic}[abs,width=\textwidth,unit=1mm,scale=.25]{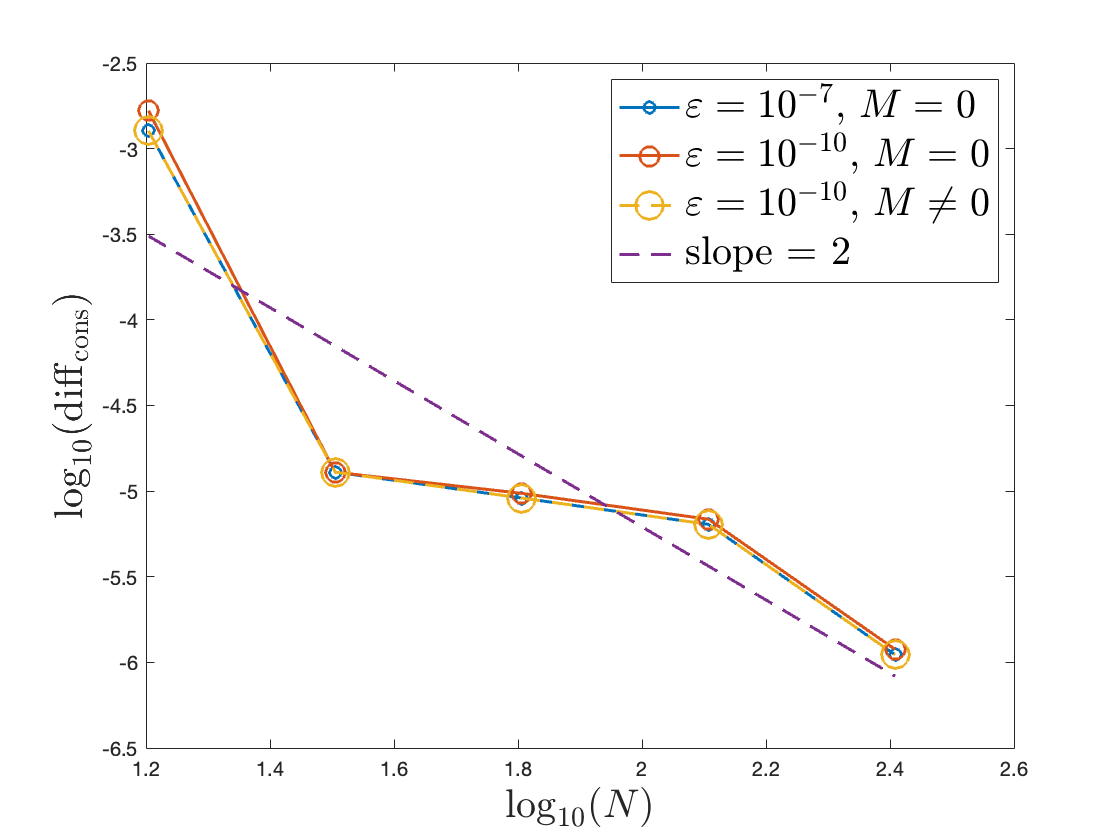}
\put(2,56){(a)}
\end{overpic}
\end{minipage}         
\begin{minipage}[b]
		{.49\textwidth}
		\centering
	\begin{overpic}[abs,width=\textwidth,unit=1mm,scale=.25]{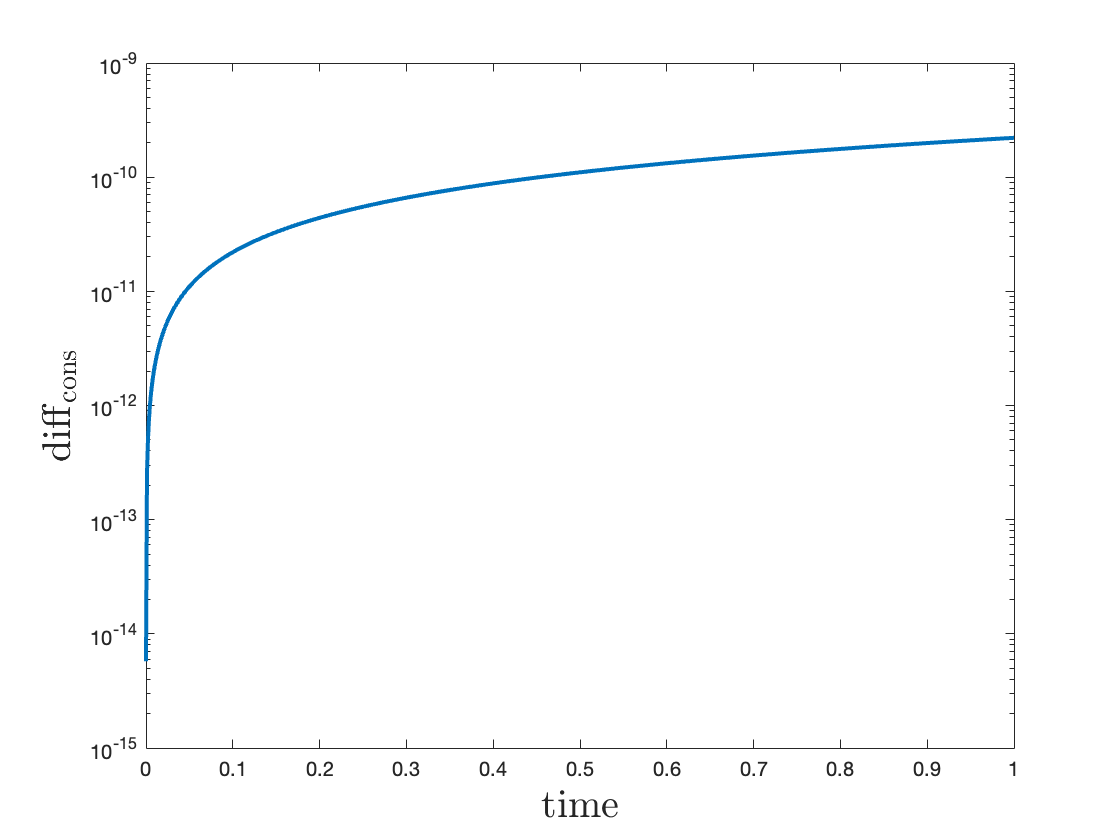}
\put(2,56){(b)}
\end{overpic}
\end{minipage}   
    \caption{\textit{Charge conservation test: we plot the difference ${\rm diff}_{\rm cons}$ defined in Eq.~\eqref{eq:diff_cons}, as a function of the number of cells of the space discretization (a). In panel (b), we show the same quantity in function of time, considering an explicit discretization in time. In (b) $\Delta t = 0.1h^2$.}}
    \label{fig:conservation}
\end{figure}


\begin{figure}[h]
    \centering
\begin{minipage}[b]
		{.49\textwidth}
		\centering
	\begin{overpic}[abs,width=\textwidth,unit=1mm,scale=.25]{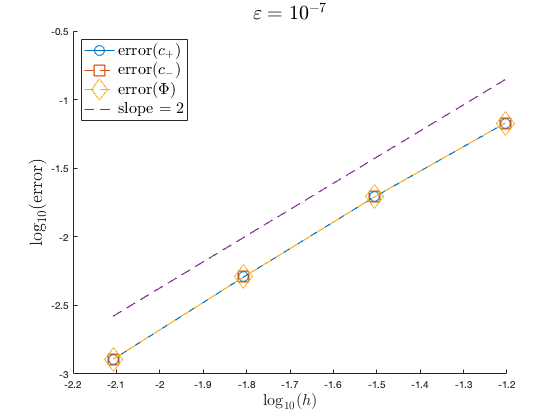}
\put(2,56){(a)}
\end{overpic}
\end{minipage}         
\begin{minipage}[b]
		{.49\textwidth}
		\centering
	\begin{overpic}[abs,width=\textwidth,unit=1mm,scale=.25]{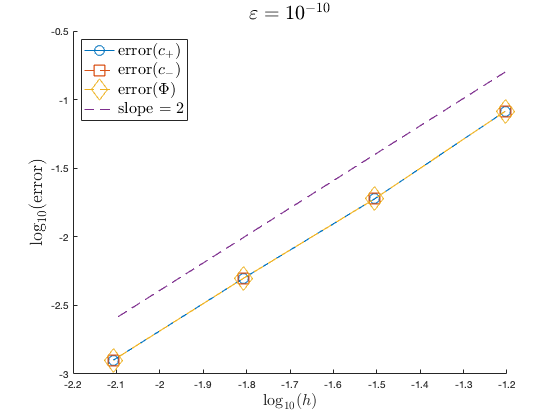}
\put(2,56){(b)}
\end{overpic}
\end{minipage}   
    \caption{\textit{Space and time accuracy of the QNL system in two dimensions (see Eqs.~(\ref{eq:vectorial}--\ref{eq_imex_Q_bi})) at final time $t = 0.3125$, for different values of $\varepsilon$. In this test the order is calculated with Richardson extrapolation technique. Simulation details are provided in Section~\ref{section_numerics_2D}.}}
    \label{fig:accuracy_2D}
\end{figure}
The initial conditions are given by
\begin{subequations}
    \begin{align}
        \mathcal{C}^{\rm in}(x,{y,}t=0) &= {c_+^{\rm in}(x,{y,}t=0)} + {c_-^{\rm in}(x,{y,}t=0)}\\
        \mathcal{Q}^{\rm in}(x,{y,}t=0) &=  \frac{c_+^{\rm in}(x,{y,}t=0) - c_-^{\rm in}(x,{y,}t=0)}{\varepsilon} 
    \end{align}
\end{subequations}
with
\begin{align} 	
\label{eq_initial_2D}
c_\pm^{\rm in}(x,{y,}t=0) &= {\frac{v_0}{2\pi\sigma^2}}\exp\left(-{\left((x-x^{\rm in}_\pm)^2 + (y-y^{\rm in}_\pm)^2 \right)}/{\sigma^2}\right) 
\end{align}
where $v_0$ denotes the total volume per unit surface. In our numerical tests $v_0 = 10^{-6}, \,x_-^{\rm in} = 0.5, \, x_+^{\rm in} = 0.4,\, y_-^{\rm in} = y_+^{\rm in} = 0.2,\, \sigma = 0.05$, the number of cells of the space discretization is $N = 100$ and $\Delta t = h$, unless otherwise specified. 

\paragraph{Charge conservation}
If zero flux boundary conditions are adopted on the external boundary, system 
\eqref{system_multiscale_all}, 
\eqref{system_multiscale_bc} conserves the total charges of both anions and cations, i.e. the sum of the bulk and surface integrals as follows
\[
    Q_+ = \int_\Omega c_+\, dx 
\quad {\rm and} \quad  
    Q_- = \int_\Omega c_-\, dx + M\int_{\Gamma_B} c_- \, d\Gamma.  
\]
However, the numerical scheme does not strictly conserve {the} total charge for two main reasons: first, the zero-flux condition is not exactly imposed on the boundary of the bubble; second, the IMEX time discretization introduces a coupling between variables at different time levels, which breaks exact conservation, even if the formulation is written in conservative form. A fully explicit (see Fig.~\ref{fig:conservation} (panel (b))) or fully implicit discretization would improve conservation {(but we do not consider this latter case; in this work, we focus on the semi-implicit strategy proposed by the IMEX discretization).}
We perform a test to check the lack of conservation of the method.
The results are illustrated in Fig.~\ref{fig:conservation}, where we show the following 
\begin{equation}
\label{eq:diff_cons}
    {\rm diff_{cons}} = \frac{|Q_-(t) - Q_-(t=0)|}{|Q_-(t=0)|}.
\end{equation} 
As it appears from Fig.~\ref{fig:conservation} (panel (a)), conservation error scales with the order of the scheme upon grid refinement. Moreover, the error does not vary significantly with respect to the parameters $\varepsilon$ and $M$.

In Fig.~\ref{fig:accuracy_2D}, we show the space and time accuracy of the numerical solutions, at final time $t = 0.3125$, for $\varepsilon = 10^{-7}$ and $\varepsilon = 10^{-10}$. In the absence of an exact solution, we apply the Richardson extrapolation technique (see, e.g., \cite{richardson1911ix}) to estimate the order of the method, and choose $\Delta t = h$. 
Moreover, to further investigate the accuracy order of the time discretization and to show that the numerical scheme maintains the same order for $\varepsilon \to 0$, in Fig.~\ref{fig:error_time} we show the $L^2$-norm of the error as a function of $\Delta t$ (and fixed $h = 10^{-2}$), at final time $t = 0.1$ and for different values of $\varepsilon$. The results confirm that the scheme maintains second-order accuracy with respect to $\varepsilon$, since the observed convergence rate remains consistent in all the values tested, including $\varepsilon = 0$. {This shows an improvement respect to the results obtained in one dimension in Figs.~\ref{fig:accuracy_1D}-\ref{fig:accuracy_1D_MPNP}, panels (c), where a degradation in the accuracy order is showed before the introduction of the new system for the variables $\mathcal C$ and $\mathcal Q$ in Section~\ref{sec:QNL}.} Although the value of the error for $\varepsilon = 10^{-11}$ is a few orders of magnitude higher, this confirms that the scheme is asymptotic preserving, since it retains its second-order accuracy even in the singular limit $\varepsilon \to 0$, where quasi-neutrality occurs. We believe that the increase in the absolute error for very small values of $\varepsilon$ is due to the conditioning number of the linear system that we solve, which becomes increasingly ill-conditioned as $\varepsilon$ decreases.

\begin{figure}[H]
    \centering
\begin{minipage}[b]
		{.49\textwidth}
		\centering
	\begin{overpic}[abs,width=\textwidth,unit=1mm,scale=.25]{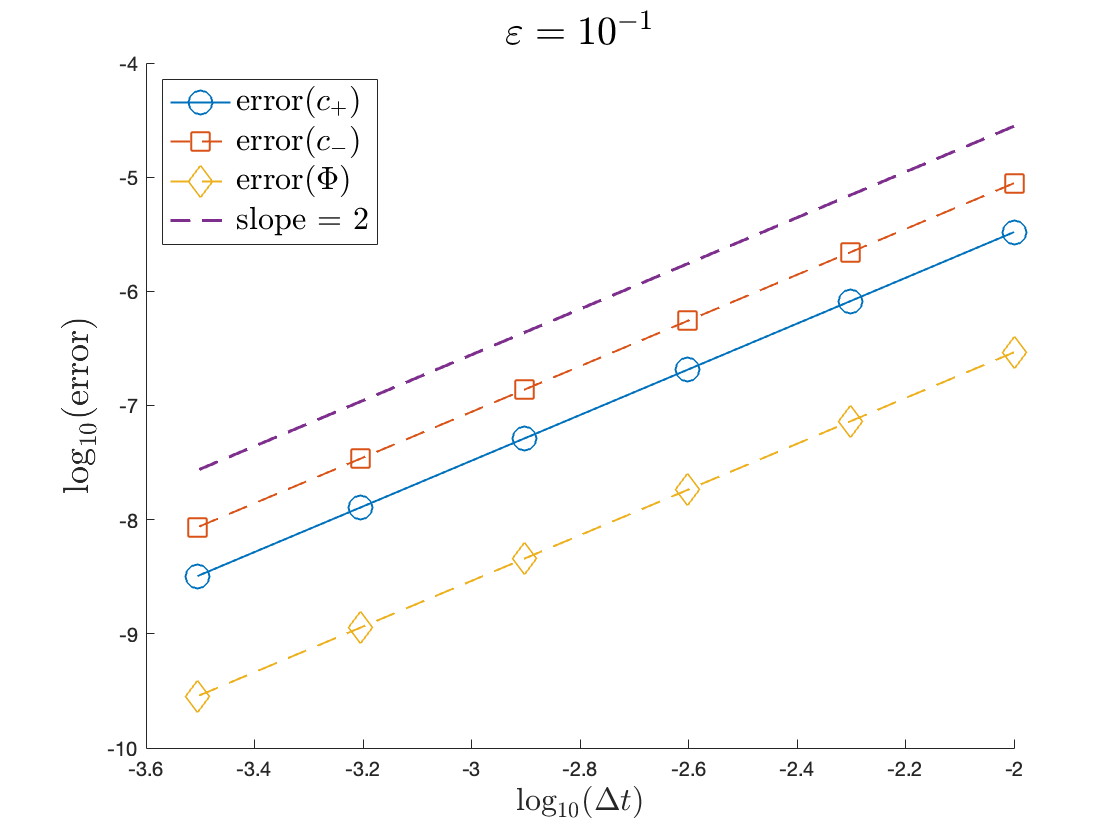}
\put(2,56){(a)}
\end{overpic}
\end{minipage}         
\begin{minipage}[b]
		{.49\textwidth}
		\centering
	\begin{overpic}[abs,width=\textwidth,unit=1mm,scale=.25]{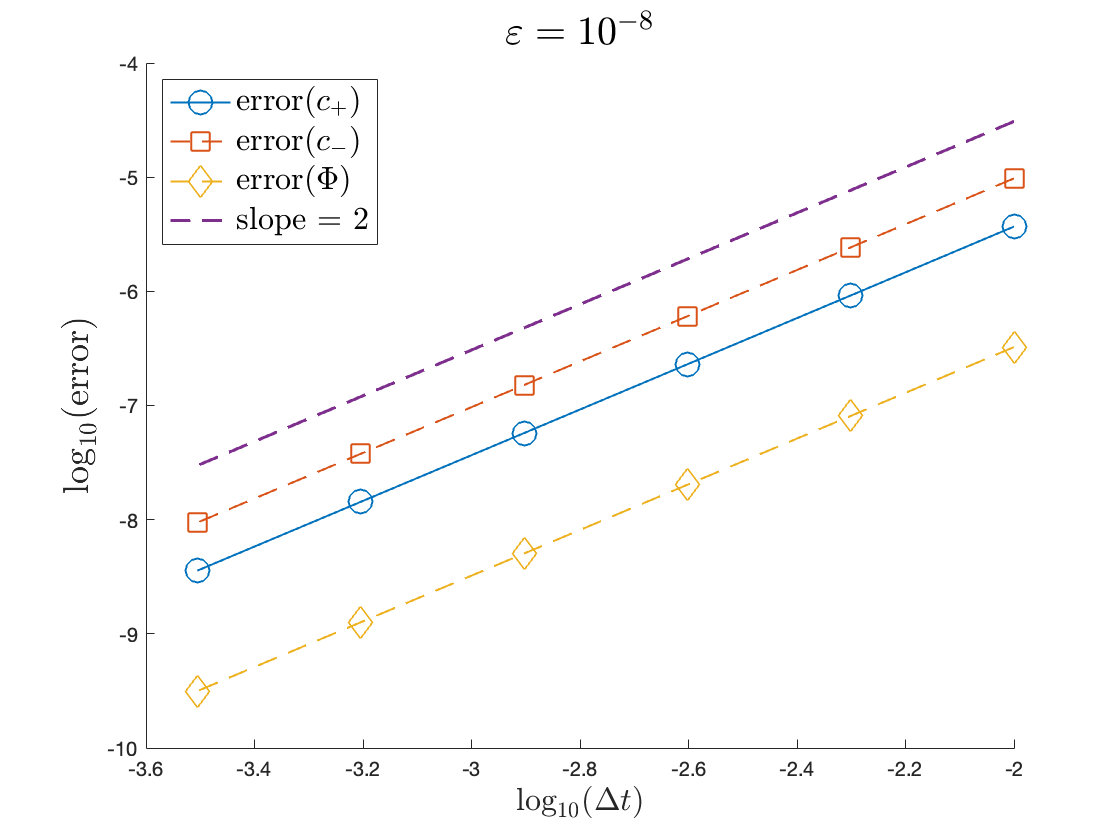}
\put(2,56){(b)}
\end{overpic}
\end{minipage}   
\begin{minipage}[b]
		{.49\textwidth}
		\centering
	\begin{overpic}[abs,width=\textwidth,unit=1mm,scale=.25]{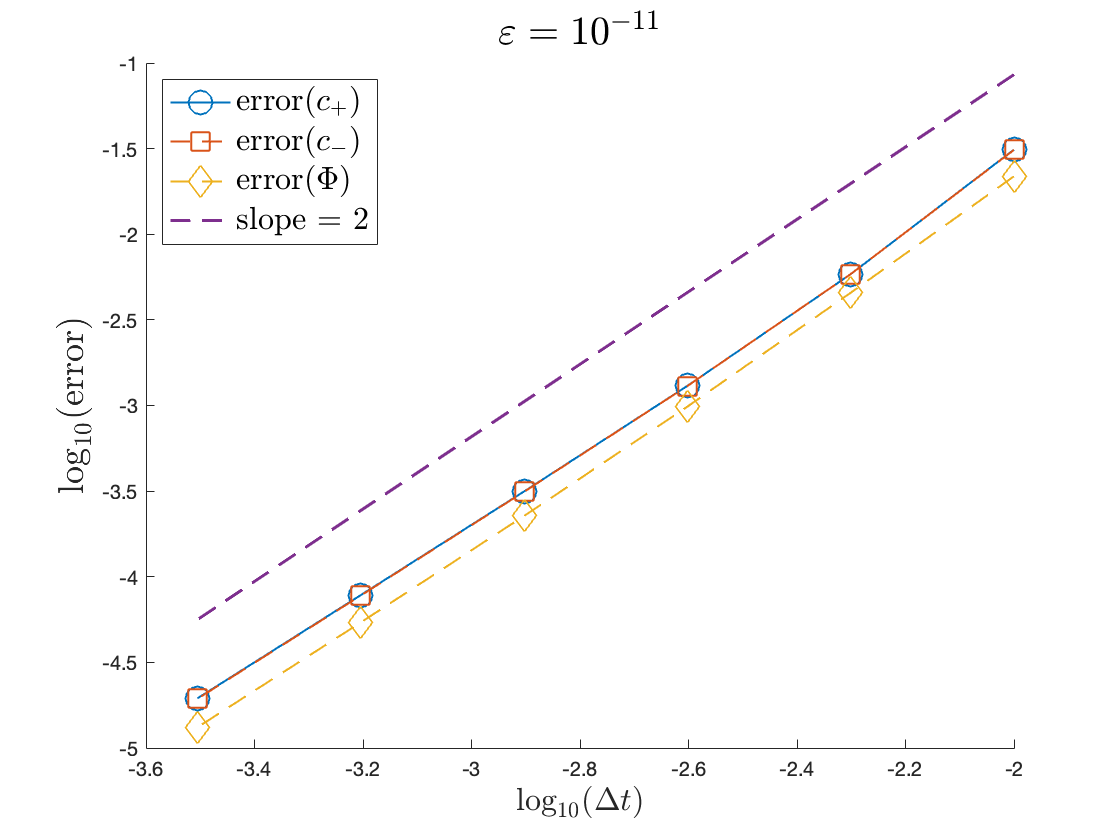}
\put(2,56){(c)}
\end{overpic}
\end{minipage}  
\begin{minipage}[b]
		{.49\textwidth}
		\centering
	\begin{overpic}[abs,width=\textwidth,unit=1mm,scale=.25]{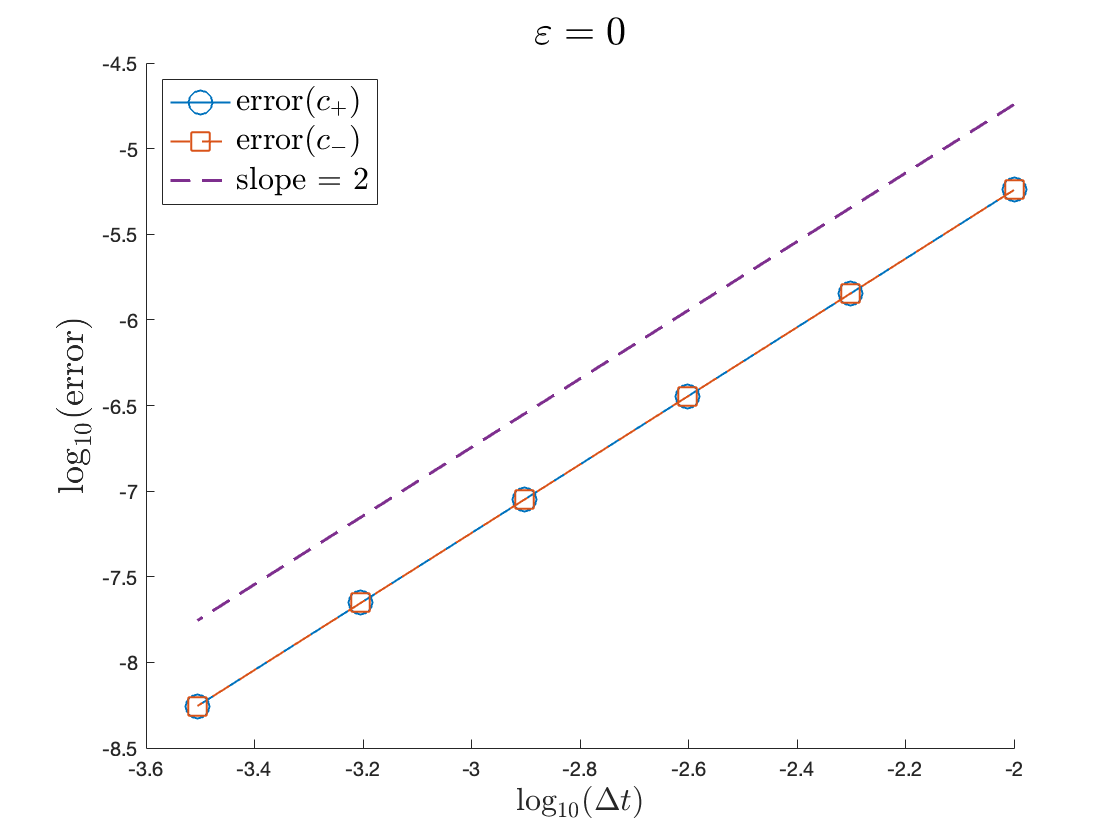}
\put(2,56){(d)}
\end{overpic}
\end{minipage}   
    \caption{\textit{Time accuracy orders of the QNL system in two dimensions (see Eqs.~(\ref{eq:vectorial}--\ref{eq_imex_Q_bi})) at final time $t = 0.1$, for different values of $\varepsilon$. Simulation details are provided in Section~\ref{section_numerics_2D}. In panel (d) the error in $\Phi$ is not plotted, however for vanishing $\varepsilon$, the value of the potential $\Phi$ remains constant within machine precision.}}
    \label{fig:error_time}
\end{figure}



\begin{figure}[H]
    \centering
\begin{minipage}[b]
		{.49\textwidth}
		\centering
	\begin{overpic}[abs,width=\textwidth,unit=1mm,scale=.25]{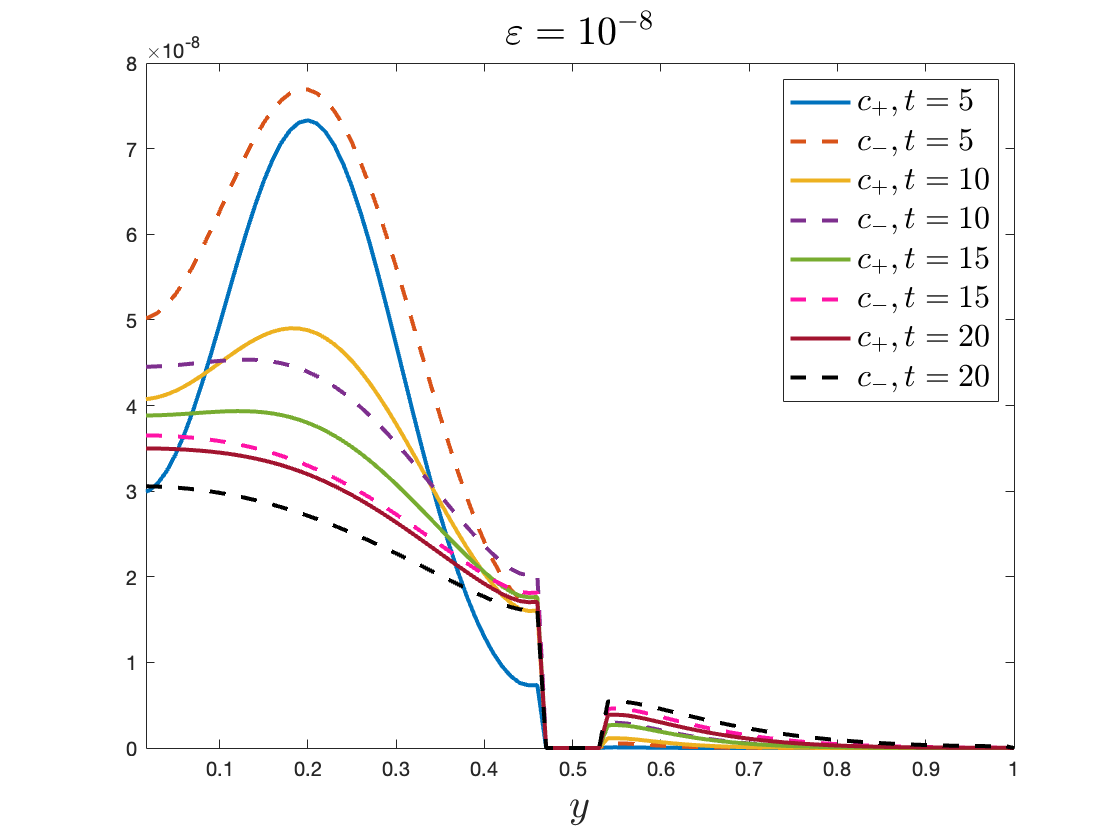}
\put(3,57){(a)}
\end{overpic}
\end{minipage}         
\begin{minipage}[b]
		{.49\textwidth}
		\centering
	\begin{overpic}[abs,width=\textwidth,unit=1mm,scale=.25]{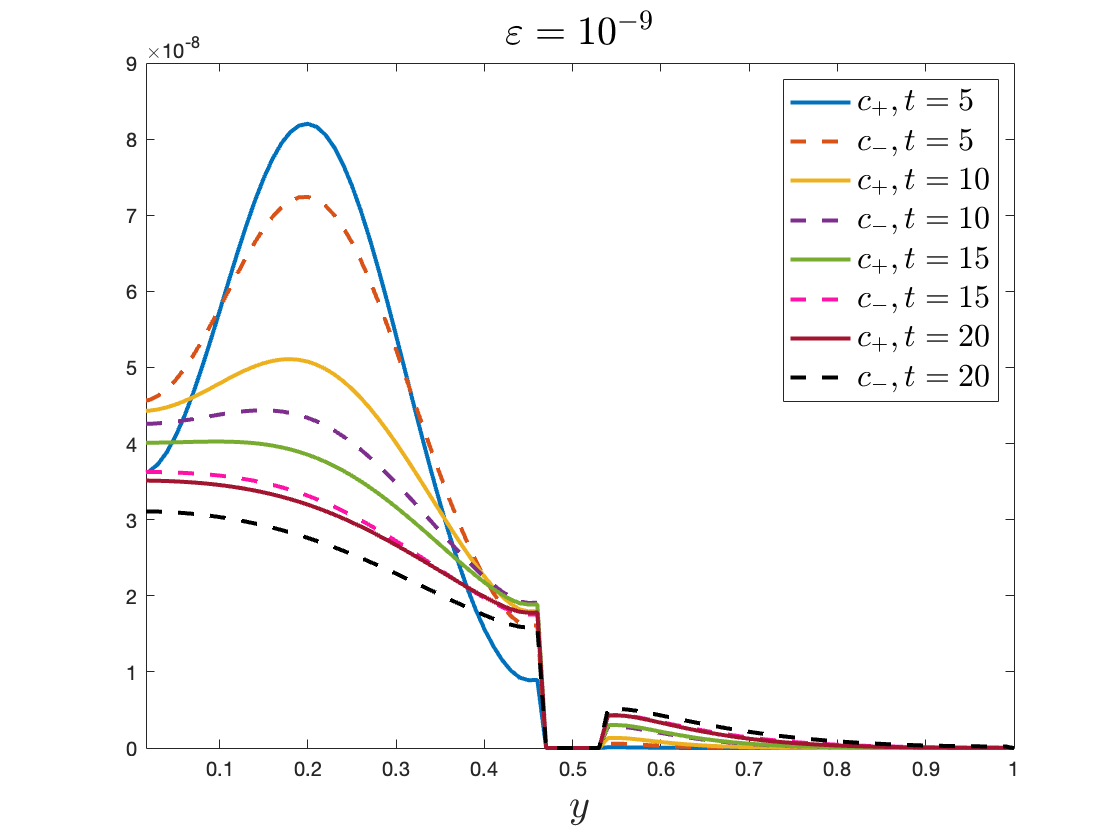}
\put(3,57){(b)}
\end{overpic}
\end{minipage}   
\begin{minipage}[b]
		{.49\textwidth}
		\centering
	\begin{overpic}[abs,width=\textwidth,unit=1mm,scale=.25]{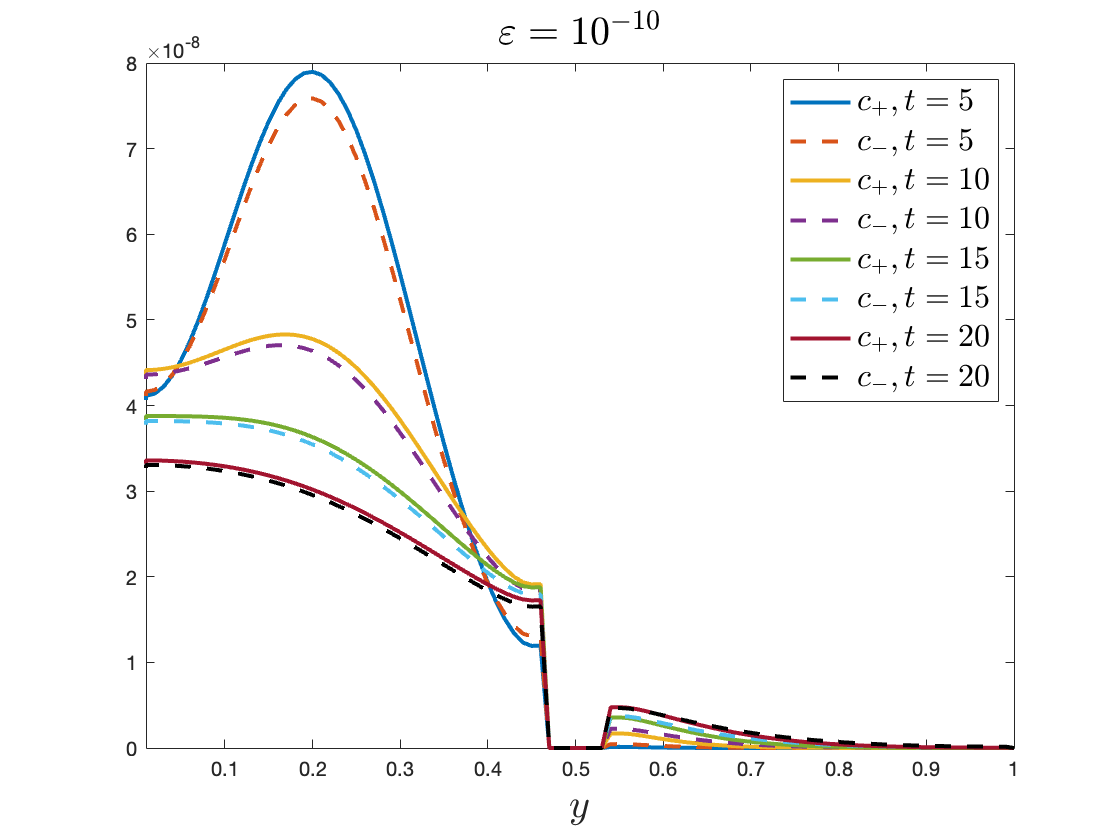}
\put(3,57){(c)}
\end{overpic}
\end{minipage}  
\begin{minipage}[b]
		{.49\textwidth}
		\centering
	\begin{overpic}[abs,width=\textwidth,unit=1mm,scale=.25]{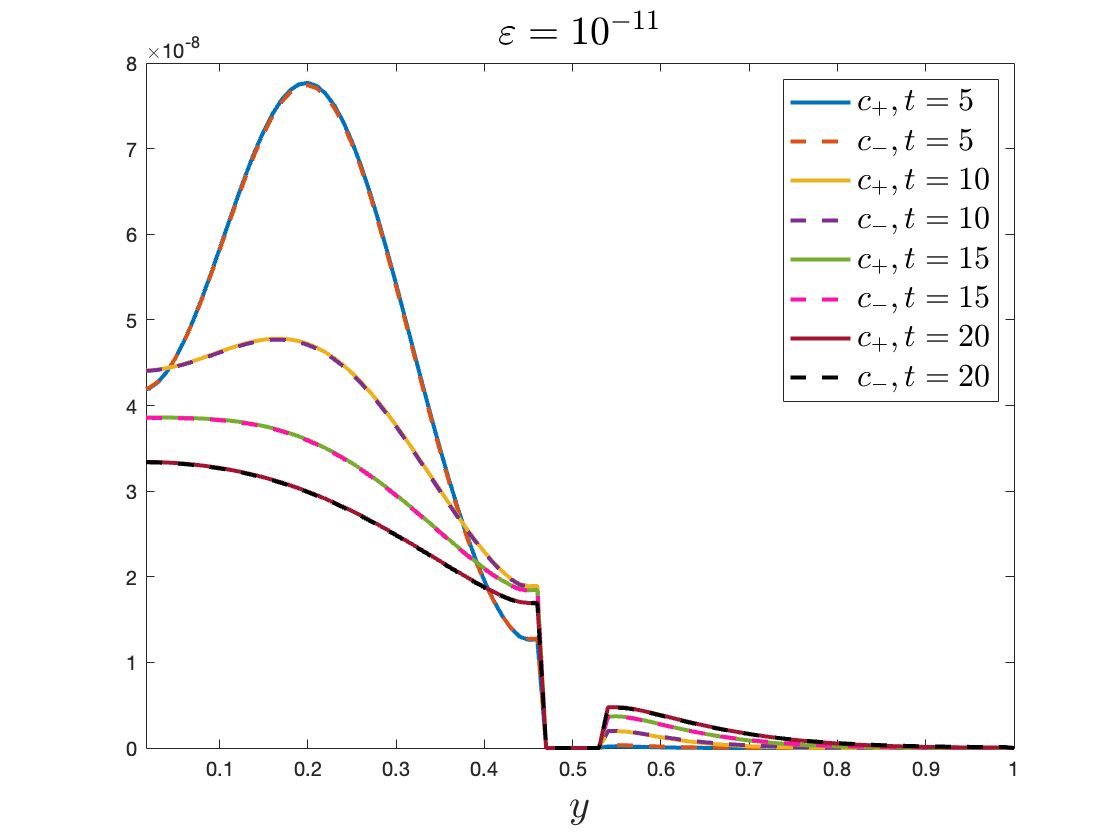}
\put(3,57){(d)}
\end{overpic}
\end{minipage}   
\begin{minipage}[b]
		{.49\textwidth}
		\centering
	\begin{overpic}[abs,width=\textwidth,unit=1mm,scale=.25]{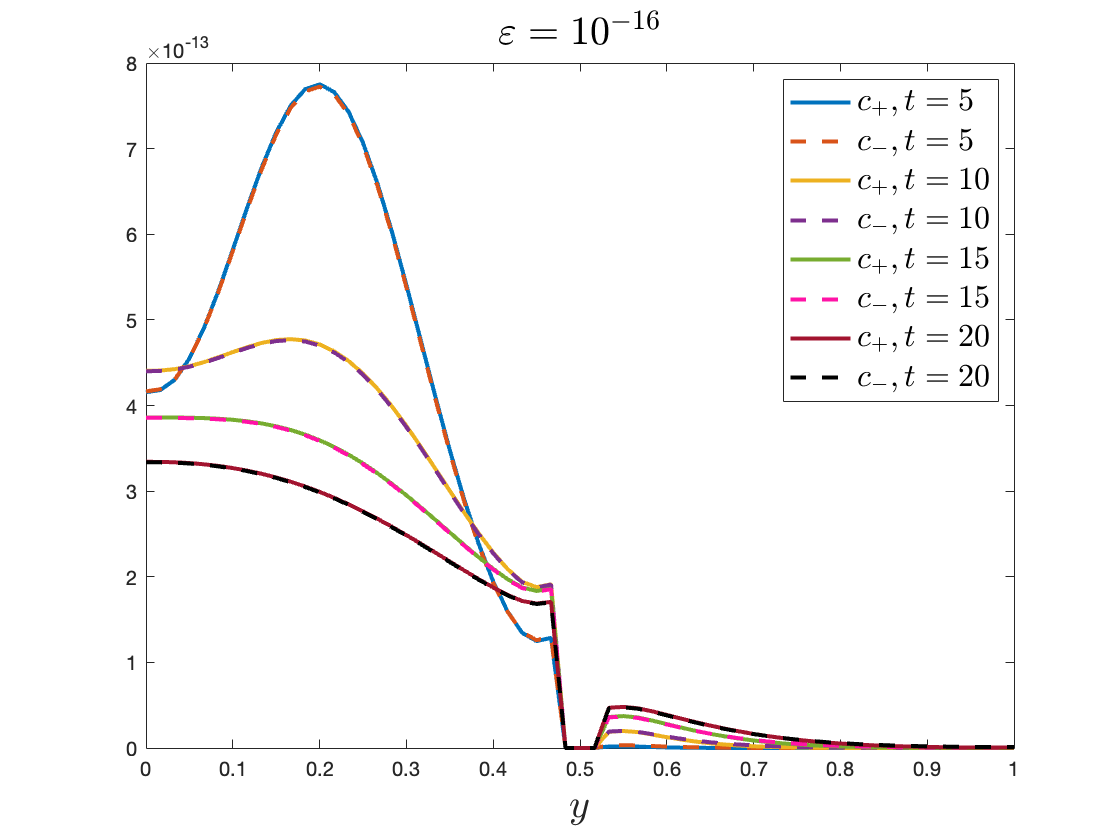}
\put(3,57){(e)}
\end{overpic}
\end{minipage}   
    \caption{\textit{Profiles of the concentrations $c_+$ (solid lines) and $c_-$ (dashed lines) in two dimensions, at $x = 0.5$. We show the solutions at different times $t = 5,10,15,20$, and for different $\varepsilon$. In panels (a)-(d) the initial volume is $v_0 = 10^{-6}$ while in panel (e) $v_0 = 10^{-11}$ and $\Delta t = 0.01 h$.}}
    \label{fig:profile_time_concentration}
\end{figure}

In Fig.~\ref{fig:profile_time_concentration}, we show the profiles of the ion concentrations $c_+$ (solid lines) and $c_-$ (dashed lines) {in $y$-direction and for a fixed $x=0.5$}, for different times and values of $\varepsilon$ {($\in \{ 10^{-8},10^{-9},10^{-10},10^{-11}\}$ and initial volume $v_0 = 10^{-6}$ in panels (a)-(d) and $\varepsilon = 10^{-16}$ and $v_0 = 10^{-11}$ in panel (e))}. We observe that the two concentrations align together when $\varepsilon\to 0$, in agreement with the Quasi-Neutral limit. We remark this aspect in Fig.~\ref{fig:difference_eps_concentration}, where we show the absolute value of the difference between the numerical density of ions (normalized respect to the initial volume $v_0$), i.e. 
\begin{equation}
\label{eq:diff}
\frac{\left|n_+ - n_-\right|}{v_0}
\end{equation}
for different values of $\varepsilon$. Finally, in Fig.~\ref{fig:bubble}, we display the anion concentration $c_-$ at different times. In the same plots, we highlight in red the values of the concentration at the boundary of the bubble $\Gamma_{\mathcal{B},h}$.

To validate the code, we study the loss of positivity of the solution. From the drift-diffusion equations, we expect the concentrations to remain positive at each time step. The reason for the investigation is the presence of the cut elements close to the bubble boundary, as we saw in \cite{astuto2024self}, and the small parameter $\varepsilon$ that can affect the coercivity of the iteration matrix. We investigate the loss of positivity in Fig.~\ref{fig:positivity} for $\varepsilon = 10^{-8}$ and $\varepsilon = 10^{-11}$. The {\tt semilogy} plots show the minimum of the concentration of anions $c_-$ (blue line) that, being negative for $t<2$ in panel (a) and for $t<4$ in panel (b), the corresponding curve is not visible due to the use of logarithmic scaling. For this reason, we also plot the absolute value of the minimum of the solution , i.e. $|{\rm min}(c_-)|$ (red dashed line). In this way, we are able to see the values of the solution for all times. For $\varepsilon = 10^{-8}$, the negative values are close to the zero machine. For $\varepsilon = 10^{-11}$ we realize that we need to investigate further because, for a few time steps, the minimum of the solution is $\approx -10^{-9}$, holding that
\( {\rm min}(c_-)\approx -10^{-2}{\rm max}(c_-).\) We consider different initial conditions for $c_\pm$ and the results show a better scenario. In Fig.~\ref{fig:positivity} and in Fig.~\ref{fig:positivity_xin} panel (a), the initial conditions are defined in Eq.~\eqref{eq_initial_2D} with $x^{\rm in}_+ = 0.4, x^{\rm in}_- = 0.5$; in Fig.~\ref{fig:positivity_xin} panel (b), $x^{\rm in}_+ = 0.45, x^{\rm in}_- = 0.5$, and the solution becomes positive earlier. We deduce that the presence of cut elements do not influence the loss positivity of the solutions, since they remain positive when $c_\pm$ is in the proximity of $\Gamma_\mathcal{B}$ (i.e., for $t>4$ the solution becomes not negligible in the proximity of the boundary and the minimum of the solutions is positive). We attribute the loss of positivity to the presence of $\varepsilon \to 0$. At the initial stage of the time evolution, the two concentrations experience a strong mutual attraction, which the time step fails to accurately capture.

\begin{figure}[H]
    \centering
\begin{minipage}[b]
		{.49\textwidth}
		\centering
	\begin{overpic}[abs,width=\textwidth,unit=1mm,scale=.25]{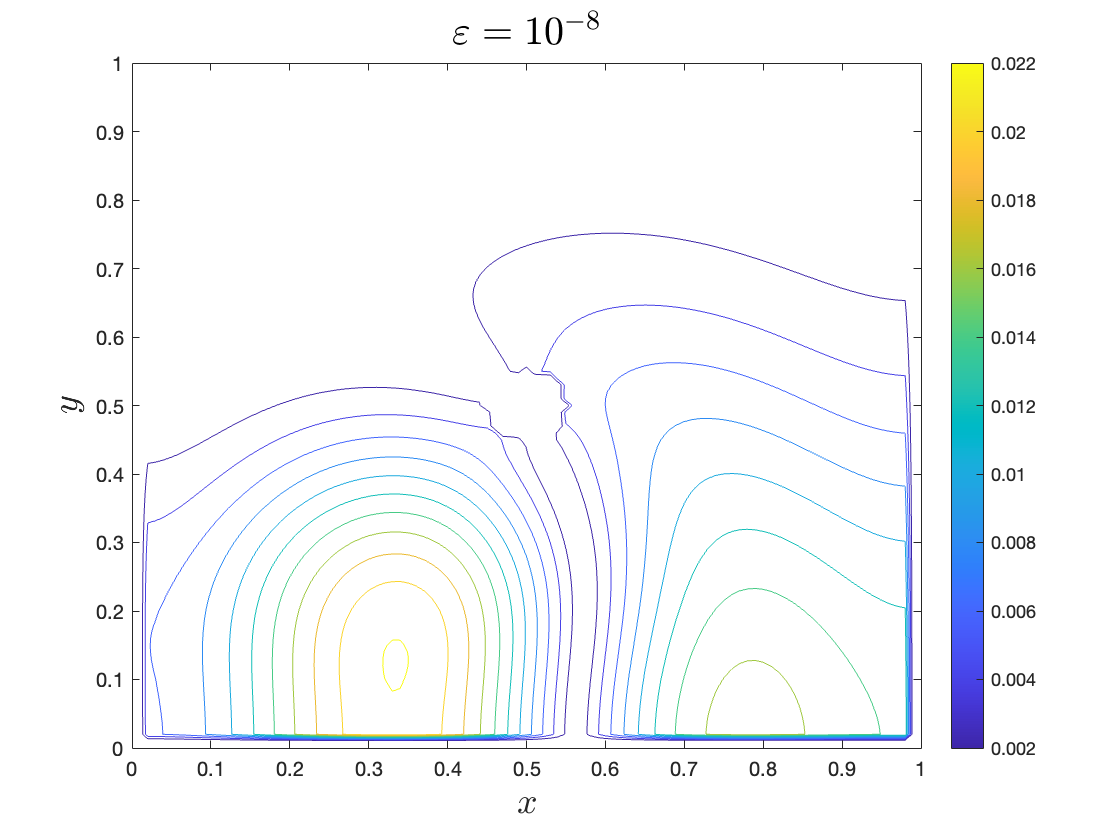}
\put(2,56){(a)}
\end{overpic}
\end{minipage}         
\begin{minipage}[b]
		{.49\textwidth}
		\centering
	\begin{overpic}[abs,width=\textwidth,unit=1mm,scale=.25]{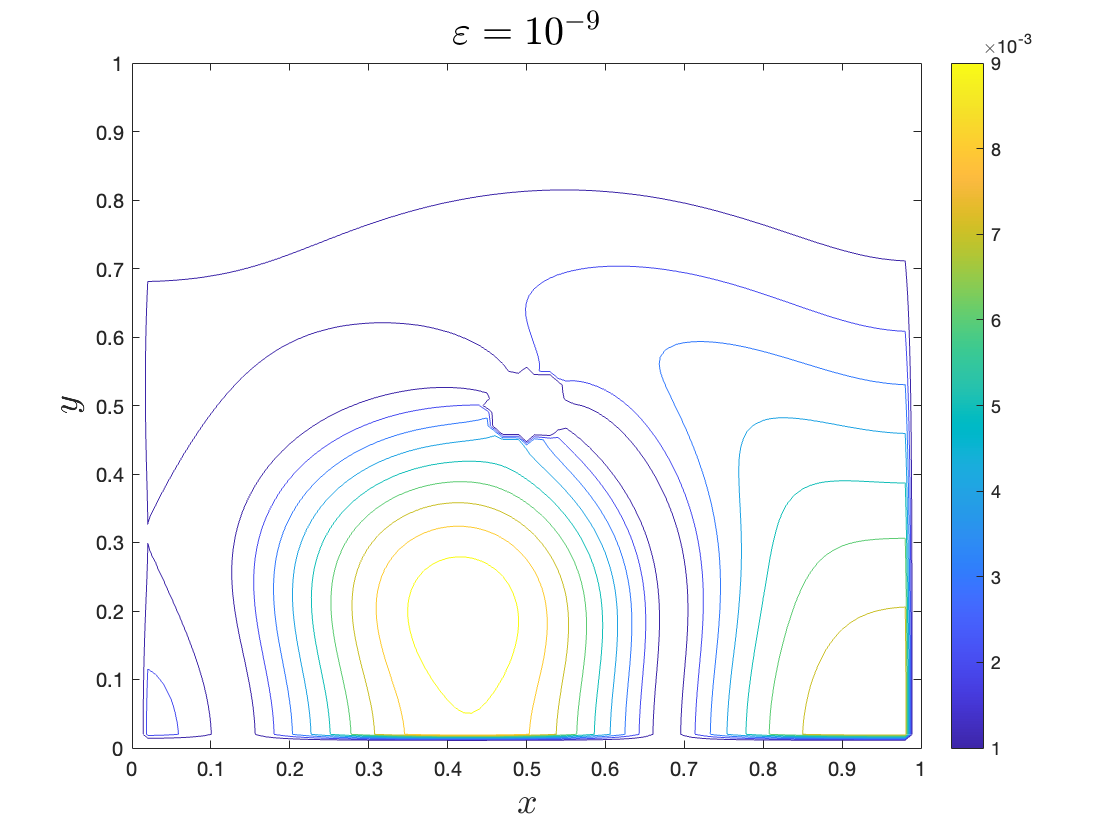}
\put(2,56){(b)}
\end{overpic}
\end{minipage}        
\begin{minipage}[b]
		{.49\textwidth}
		\centering
	\begin{overpic}[abs,width=\textwidth,unit=1mm,scale=.25]{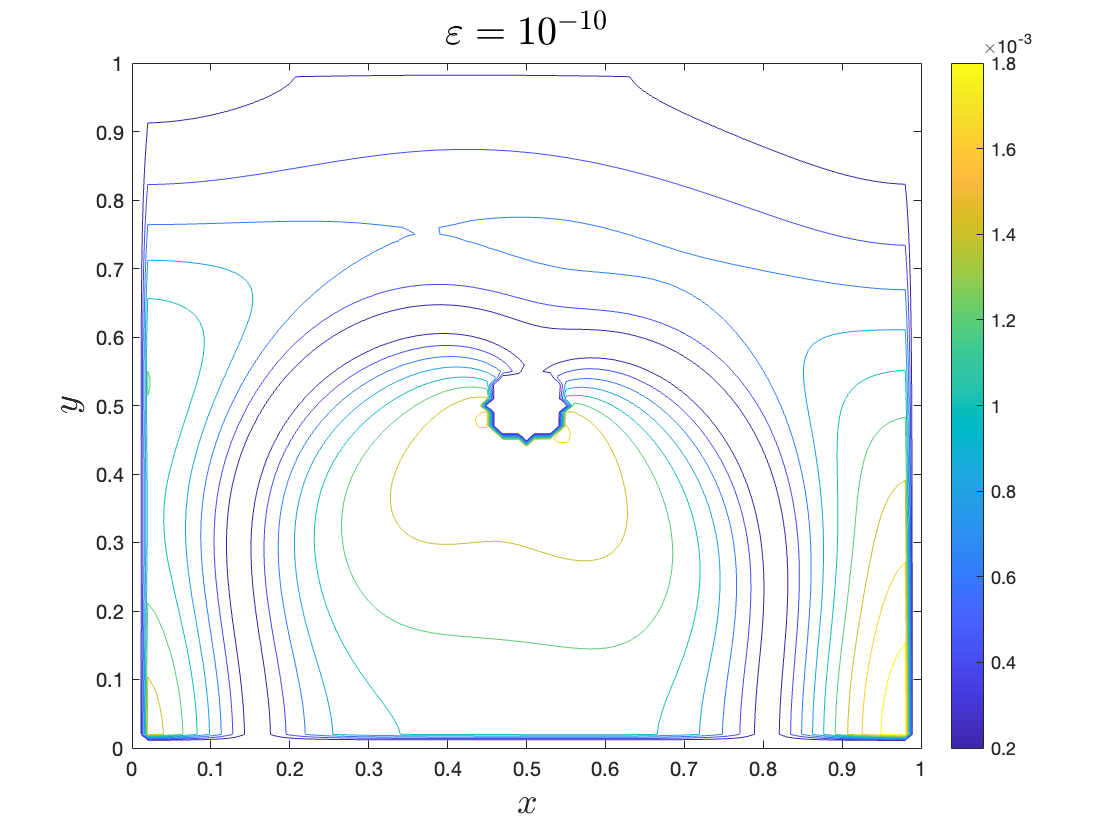}
\put(2,56){(c)}
\end{overpic}
\end{minipage}    
\begin{minipage}[b]
		{.49\textwidth}
		\centering
	\begin{overpic}[abs,width=\textwidth,unit=1mm,scale=.25]{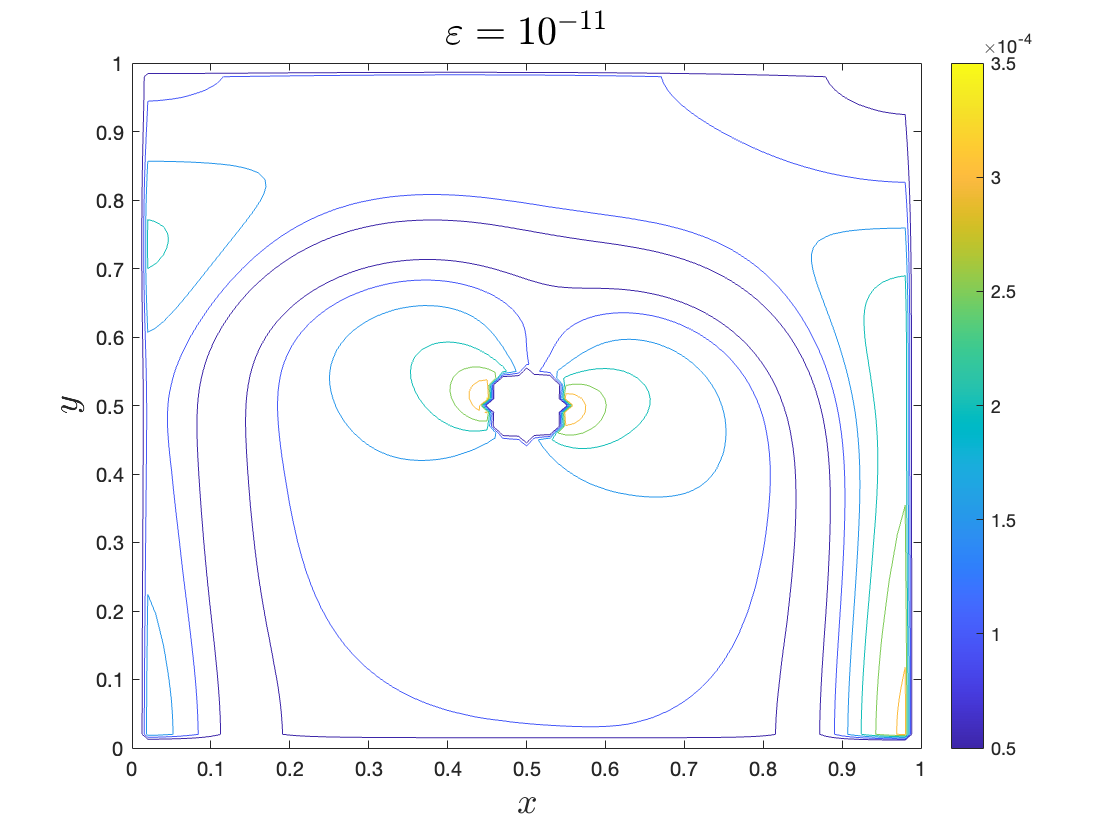}
\put(2,56){(d)}
\end{overpic}
\end{minipage}     
    \caption{\textit{Contour plot of the difference defined in Eq.~\eqref{eq:diff}, for different values of $\varepsilon$.}}
\label{fig:difference_eps_concentration}
\end{figure}

\begin{figure}[h]
    \centering
\begin{minipage}[b]
		{.49\textwidth}
		\centering
	\begin{overpic}[abs,width=\textwidth,unit=1mm,scale=.25]{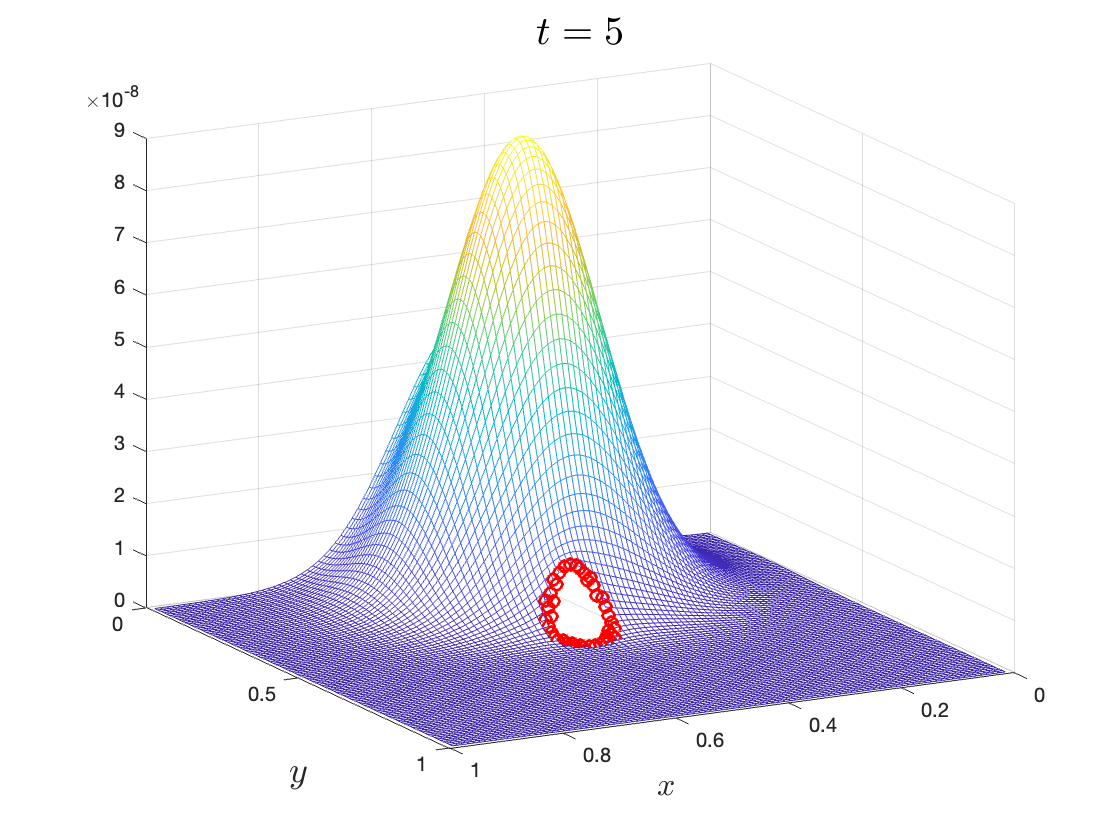}
\put(2,56){(a)}
\end{overpic}
\end{minipage}         
\begin{minipage}[b]
		{.49\textwidth}
		\centering
	\begin{overpic}[abs,width=\textwidth,unit=1mm,scale=.25]{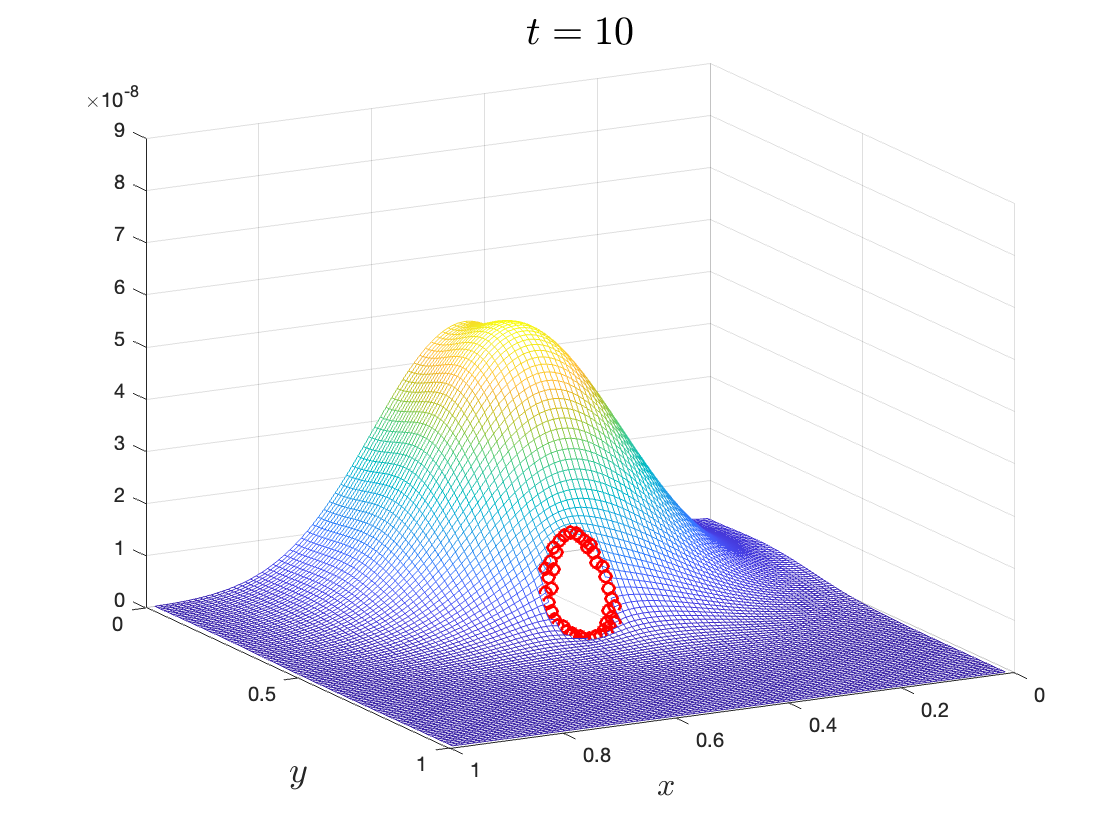}
\put(2,56){(b)}
\end{overpic}
\end{minipage}   
\begin{minipage}[b]
		{.49\textwidth}
		\centering
	\begin{overpic}[abs,width=\textwidth,unit=1mm,scale=.25]{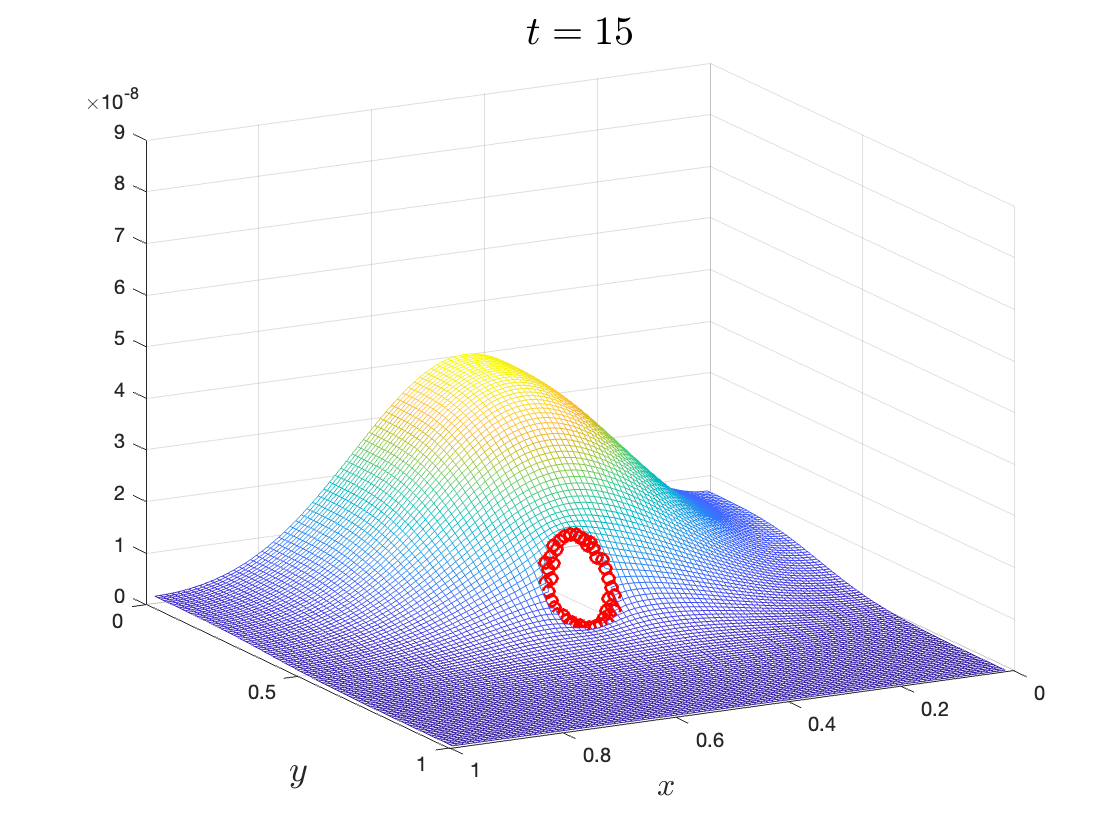}
\put(2,56){(c)}
\end{overpic}
\end{minipage}  
\begin{minipage}[b]
		{.49\textwidth}
		\centering
	\begin{overpic}[abs,width=\textwidth,unit=1mm,scale=.25]{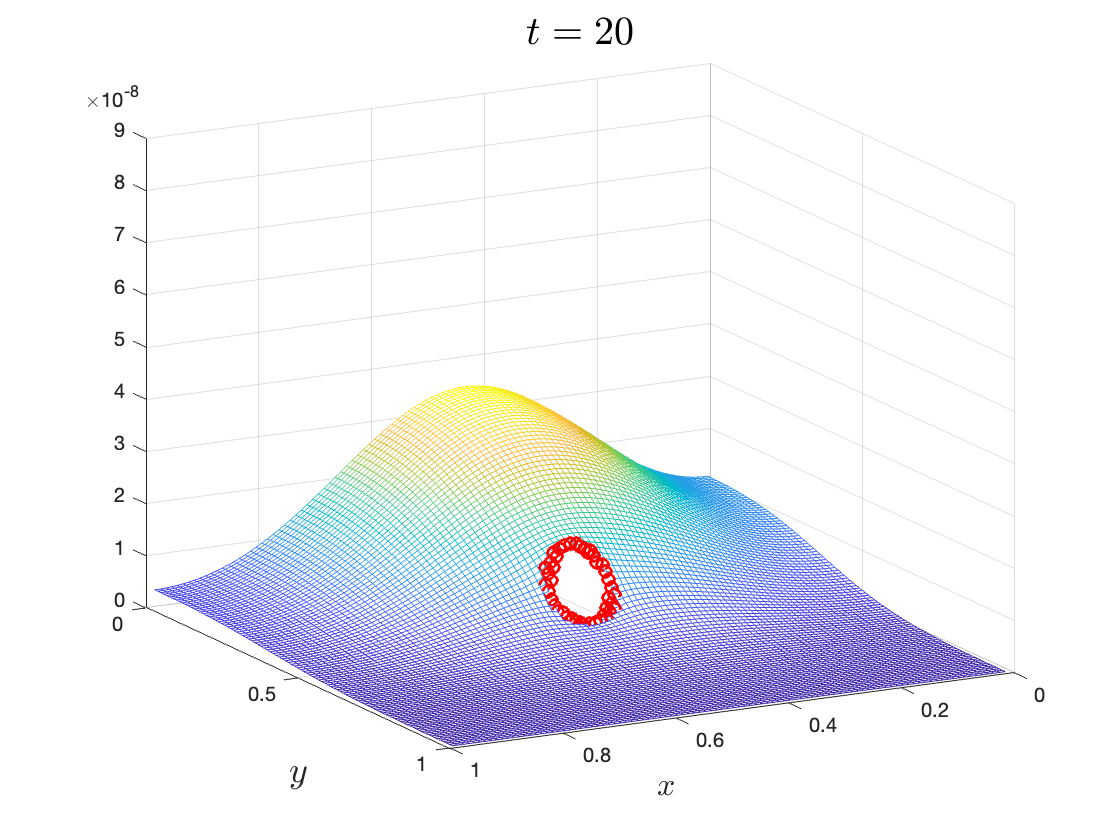}
\put(2,56){(d)}
\end{overpic}
\end{minipage}   
    \caption{\textit{Time evolution of anion concentration  $c_-$ at different times $t=5, 10, 15$, and $20$. We mark in red the concentration values at the boundary of the bubble $\Gamma_{\mathcal{B},h}$.}}
    \label{fig:bubble}
\end{figure}

\begin{figure}[h]
    \centering
\begin{minipage}[b]
		{.49\textwidth}
		\centering
	\begin{overpic}[abs,width=\textwidth,unit=1mm,scale=.25]{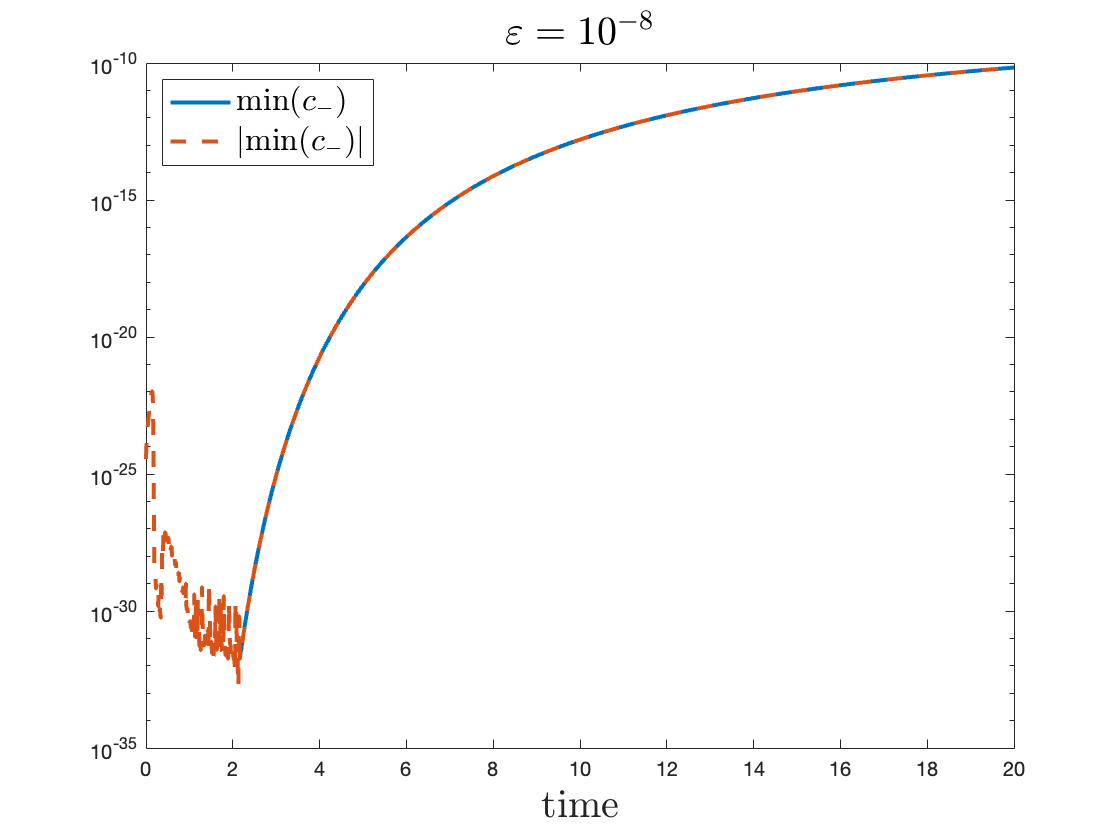}
\put(2,56){(a)}
\end{overpic}
\end{minipage}         
\begin{minipage}[b]
		{.49\textwidth}
		\centering
	\begin{overpic}[abs,width=\textwidth,unit=1mm,scale=.25]{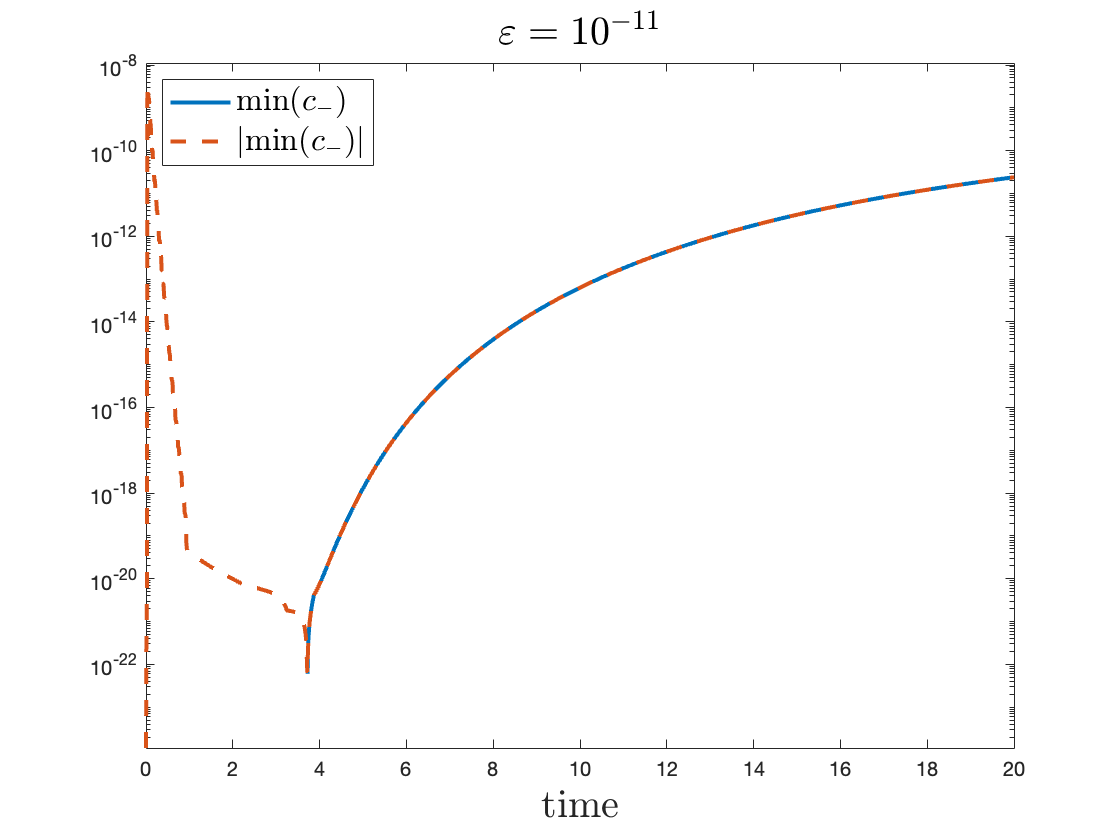}
\put(2,56){(b)}
\end{overpic}
\end{minipage}         
    \caption{\textit{Positivity of the solution for different values of $\varepsilon$.}}
    \label{fig:positivity}
\end{figure}

\begin{figure}[h]
    \centering
\begin{minipage}[b]
		{.49\textwidth}
		\centering
	\begin{overpic}[abs,width=\textwidth,unit=1mm,scale=.25]{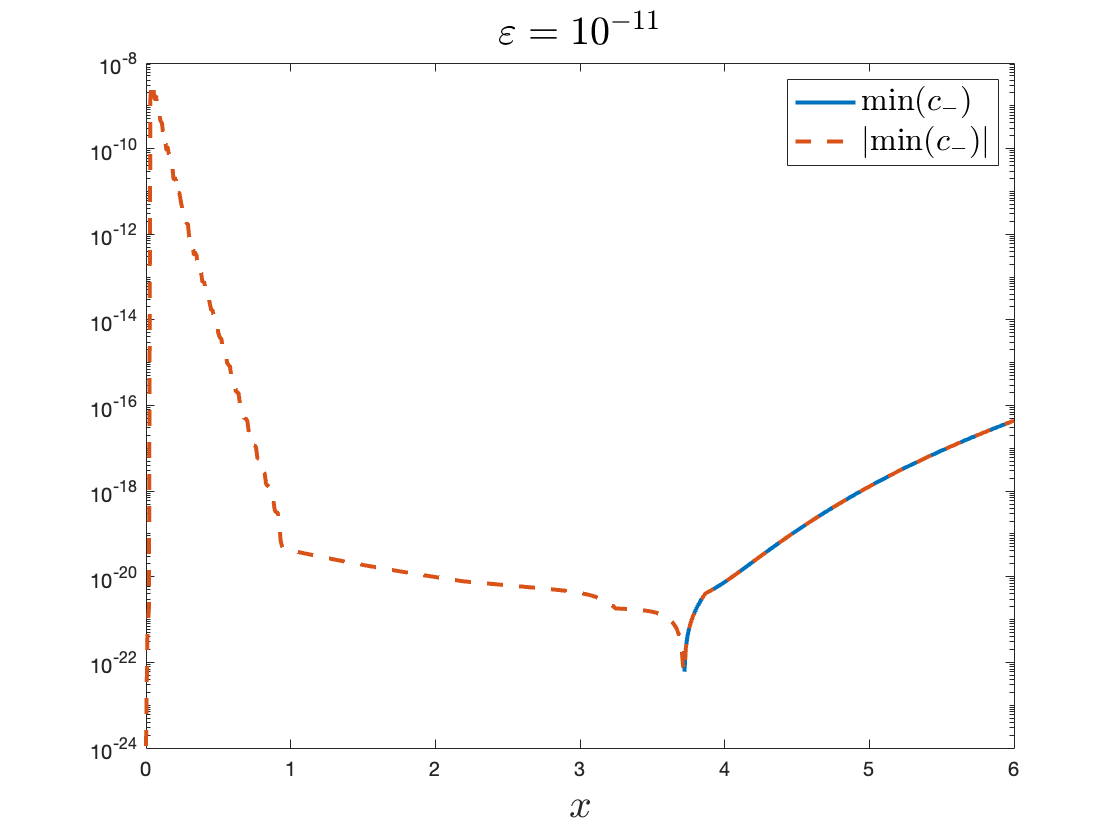}
\put(2,56){(a)}
\end{overpic}
\end{minipage}         
\begin{minipage}[b]
		{.49\textwidth}
		\centering
	\begin{overpic}[abs,width=\textwidth,unit=1mm,scale=.25]{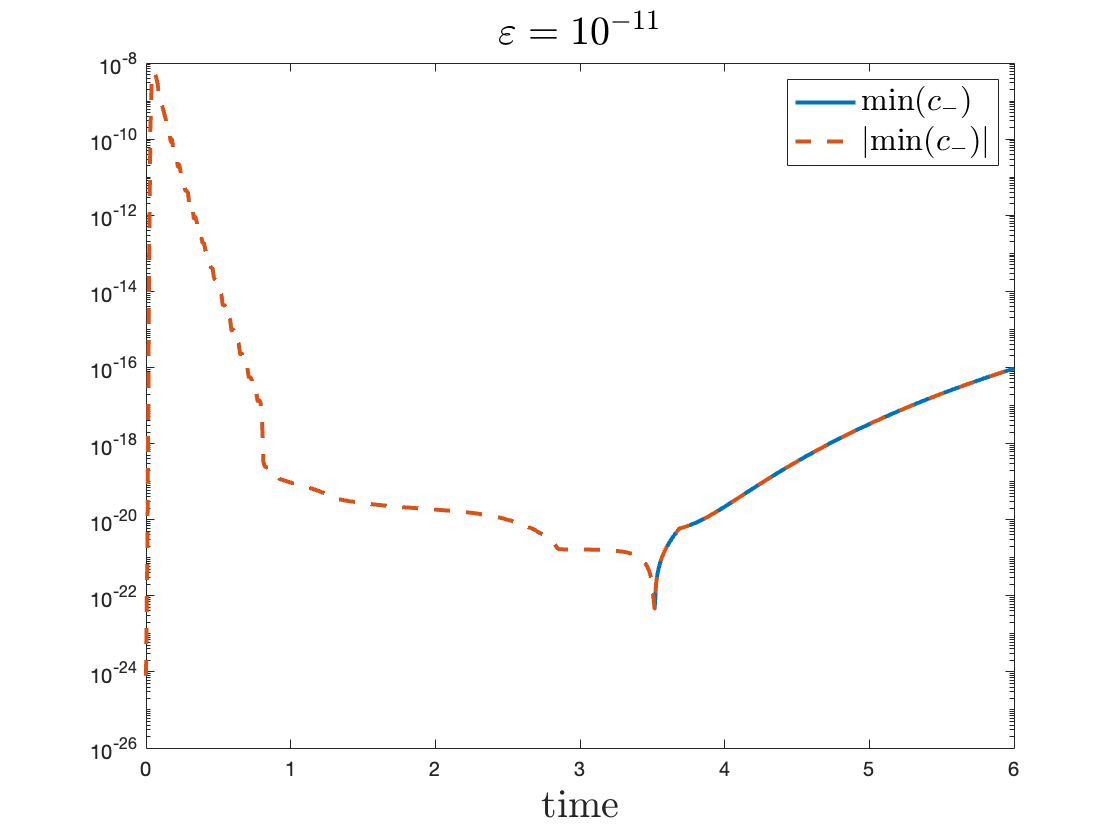}
\put(2,56){(b)}
\end{overpic}
\end{minipage}         
    \caption{\textit{Positivity of the solution for different initial conditions {defined in Eq.~\eqref{eq_initial_2D}, with $x^{\rm in}_+ = 0.4, x^{\rm in}_- = 0.5$ in (a), and $x^{\rm in}_+ = 0.45, x^{\rm in}_- = 0.5$ in (b).}}}
    \label{fig:positivity_xin}
\end{figure}

\section{Conclusions}
\label{sec:conclusions}
In this work, we present a multiscale model for a two-species Poisson-Nernst-Planck (PNP) system that describes the correlated dynamics of positive and negative ions in the presence of the trap. The model is derived from a system of two drift-diffusion equations, where the drift terms account for the gradient of a potential representing the effect of the bubble. The proposed Multiscale PNP model (MPNP) relies on the assumption that the potential range is much smaller than the relevant macroscopic length scales, such as the radius of a spherical trap. Building on our previous work \cite{astuto2023multiscale}, we show that the anion concentration follows a Boltzmann-type distribution. This leads to a significant simplification of the system, resulting in the substitution of the small scale interaction with an evolutionary time dependent boundary condition. 
	
The MPNP model is then solved in time with a second order IMEX scheme. The model is carefully numerically validated in one dimension against a detailed numerical solution of a fully resolved model with a potential of width $\delta$.
We show that the new MPNP model asymptotically coincides with the full model in the limit as the potential length $\delta \to 0$.

A key contribution of this study lies in the accurate treatment of the Coulomb interaction between ions in regimes where the Debye length is small but not negligible. While the Quasi-Neutral limit provides a simplified model in the asymptotic regime $\varepsilon = 0$, the case $\varepsilon \ll 1$ introduces significant numerical challenges: the system becomes stiff and the condition number of the discretized matrix increases, leading to potential loss of accuracy and increase of the computational time to solve the linear system. To address these issues, we develop a second-order Asymptotic Preserving (AP) scheme that ensures uniform accuracy across a wide range of Debye lengths. Finally, we validated the code by examining the loss of positivity in the solution. The results suggest that the issue is not caused by cut elements near the bubble boundary but is instead related to the small parameter $\varepsilon$, which affects the dynamics in the first part of the time evolution. {A further analysis of the AP numerical schemes is showed in \cite{astuto2025standard}, where a third order accurate numerical scheme is proposed. Thanks to the adaptability of the IMEX technique, we obtain a high order time discretization for the problem. Moreover, in that paper, we show that the proposed numerical scheme is robust and independent of the choice of the initial data.}

A natural extension of this work would be to include saturation effects, which become relevant for non-negligible ion concentrations near the bubble surface. Such effects, partially addressed in \cite{astuto2023multiscale}, require the development of nonlinear boundary conditions. 

An additional objective of future investigation is the improvement of the conservation properties of the scheme{, and a gradient flow structure for the multiscale model for one carriers is already under investigation in \cite{astuto2026gradient}.}

\begin{table}
	\centering  
	\begin{tabular}{|c|c|c|c|c|c|}
		\hline
		Symbol & value & Symbol & value & Symbol & value\\ 
		\hline\hline 
		$D_0$ & $10^{-9}m^2 s^{-1}$ & $D_+/D_0$ & $1.5$  & $D_-/D_0$ & $0.5$\\
		\hline 
		$\epsilon_0$ & $8.8541\times10^{-12}\,F m^{-1}$ & $\epsilon_r$ & $78$ & $\rho$ & $10^{3}\,Kg\, m^{-3}$ \\
		\hline		
		$m_0$ & $10^{-3}\,$ Kg$\, mol^{-1}$ & { $m_{H_2O}$ } & { $18$ } &$q$ & $1,602\times 10^{-19}C$\\
		\hline
		  $k_B$ & $1.38 \times 10^{-23}J/K$ & T & $300K$ & $N_A$  & $6.022\times10^{23}mol^{-1}$ \\
		\hline 
	\end{tabular}
	\caption{\textit{Parameters involved.}}
	\label{table_parameters}
\end{table}

\section*{Acknowledgments}
{The work has been supported by the Spoke 10 Future AI Research (FAIR) of the Italian Research Center funded by the Ministry of University and Research as part of the National Recovery and Resilience Plan (PNRR). \\ 
The authors have also been supported also by Italian Ministerial grant PRIN 2022 'Efficient numerical schemes and optimal control methods for time-dependent partial differential equations', No. 2022N9BM3N - Finanziato dall'Unione europea - Next Generation EU.} \\
{The work has also been supported by the Italian Ministerial grant PRIN 2022 PNRR 'FIN4GEO: Forward and Inverse Numerical Modeling of hydrothermal systems in volcanic regions with application to geothermal energy exploitation', No. P2022BNB97 - Finanziato dall'Unione europea - Next Generation EU. } \\
The authors have also been supported also by Italian Ministerial grant PRIN 2022 'Advanced numerical methods for time dependent parametric partial differential equations with applications', No. 2022KA3JBA - Finanziato dall'Unione europea - Next Generation EU. \\
{The authors are members of the Gruppo Nazionale Calcolo Scientifico-Istituto Nazionale di Alta Matematica (GNCS-INdAM).}
We acknowledge the CINECA award under the ISCRA initiative, for the availability of high-performance computing resources.

\bibliographystyle{plain}
\bibliography{bibliography.bib}

\end{document}